\documentclass{pspum-l}

\usepackage{amssymb}

\usepackage[cmtip,all]{xy}

\usepackage{hyperref}

\newtheorem{theorem}{Theorem}[section]
\newtheorem{lemma}[theorem]{Lemma}
\newtheorem{prop}[theorem]{Proposition}
\newtheorem{cor}[theorem]{Corollary}
\newtheorem{conjecture}[theorem]{Conjecture}

\theoremstyle{definition}
\newtheorem{definition}[theorem]{Definition}
\newtheorem{example}[theorem]{Example}
\newtheorem{examples}[theorem]{Examples}

\theoremstyle{remark}
\newtheorem{remark}[theorem]{Remark}

\numberwithin{equation}{section}

\newcommand\quash[1]{}
\newcommand\til[1]{\widetilde{#1}}

\newcommand\Nat{\mathbb{N}}
\newcommand\Z{\mathbb{Z}}
\newcommand\Q{\mathbb{Q}}
\newcommand\oQ{\overline{\mathbb{Q}}}
\newcommand\R{\mathbb{R}}
\newcommand\C{\mathbb{C}}
\newcommand\Fi{\mathbb{F}}
\newcommand\Ham{\mathbb{H}}
\newcommand\Ade{\mathbb{A}}
\newcommand\Af{\mathbb{A}_f}
\newcommand\Se{\mathbb{S}}

\newcommand\Acal{\mathcal{A}}
\newcommand\Dcal{\mathcal{D}}
\newcommand\Gf{\mathcal{G}}
\newcommand\Hf{\mathcal{H}}
\newcommand\Lf{\mathcal{L}}
\newcommand\Mcal{\mathcal{M}}
\newcommand\Of{\mathcal{O}}
\newcommand\Pf{\mathcal{P}}
\newcommand\Xf{\mathcal{X}}

\newcommand\ggoth{\mathfrak{g}}
\newcommand\hgoth{\mathfrak{h}}
\newcommand\kgoth{\mathfrak{k}}
\newcommand\Sgoth{\mathfrak{S}}

\newcommand\Aut{\mathrm{Aut}}
\newcommand\End{\mathrm{End}}
\renewcommand\H{\mathrm{H}}
\newcommand\Hom{\mathrm{Hom}}

\newcommand\GL{\mathrm{GL}}
\newcommand\GSp{\mathrm{GSp}}
\newcommand\GSpin{\mathrm{GSpin}}
\newcommand\GO{\mathrm{GO}}
\newcommand\GU{\mathrm{GU}}
\newcommand\Orth{\mathrm{O}}
\newcommand\PGL{\mathrm{PGL}}
\newcommand\PSp{\mathrm{PSp}}
\newcommand\PSU{\mathrm{PSU}}
\newcommand\SL{\mathrm{SL}}
\newcommand\SO{\mathrm{SO}}
\newcommand\Sp{\mathrm{Sp}}
\newcommand\SU{\mathrm{SU}}
\newcommand\U{\mathrm{U}}

\newcommand\lan\langle
\newcommand\ran\rangle

\newcommand\ab{\mathrm{ab}}
\newcommand\ad{\mathrm{ad}}
\newcommand\Ad{\mathrm{Ad}}
\DeclareMathOperator\card{card}
\DeclareMathOperator\Cent{Cent}
\newcommand\der{\mathrm{der}}
\DeclareMathOperator\diag{diag}
\newcommand\disc{\mathrm{disc}}
\newcommand\et{\mathrm{\acute{e}t}} 
\newcommand\flnom[1]{\stackrel{#1}{\ra}}
\newcommand\Frob{\mathrm{Frob}}
\DeclareMathOperator\Gal{Gal}
\newcommand\id{\mathbf{id}}
\DeclareMathOperator\im{Im}
\DeclareMathOperator\inv{inv}
\newcommand\iso{\stackrel{\sim}{\ra}}
\DeclareMathOperator\Ker{Ker}
\DeclareMathOperator\Lie{Lie}
\newcommand\op{\mathrm{op}}
\newcommand\Pic{\mathbf{Pic}}
\newcommand\Proj{\mathbb{P}}
\newcommand\ra\rightarrow
\DeclareMathOperator\Res{Res}

\newcommand\Sch{\mathrm{Sch}}
\DeclareMathOperator\Spec{Spec}
\newcommand\sous\backslash

\DeclareMathOperator\Tr{Tr}
\DeclareMathOperator\vol{vol}

\makeatletter
\DeclareRobustCommand{\Udots}{%
  \vcenter{\offinterlineskip
    \halign{%
      \hbox to .8em{##}\cr
      \hfil.\cr\noalign{\kern.2ex}
      \hfil.\hfil\cr\noalign{\kern.2ex}
      .\hfil\cr}%
  }%
}
\makeatother

\newcommand\ungras{1\mkern -5mu\mathrm{l}}

\begin{document}

\title{Shimura varieties (version of \today)}


\author{Sophie Morel}
\address{UMPA UMR CNRS 5669\\
ENS Lyon Site Monod\\
46 Allée d'Italie\\
69364 Lyon Cedex 07}
\email{sophie.morel@ens-lyon.fr}

\thanks{I thank the organizers of the IHES summer school (Pierre-Henri
Chaudouard, Wee Teck Gan, Tasho Kaletha and Yannis Sakellaridis) for the
enormous amount of work they put into the school, and the audience for
putting up with my lectures. I also thank Mikhail Borovoi and Yujie
Xu for remarks and corrections on an earlier version of these notes, as well
as the anonymous referee for their careful reading of the texts and numerous
comments and corrections.}

\subjclass[2020]{Primary }

\date{today}

\begin{abstract}
Lecture 1 gives an introduction to Shimura varieties over the complex
numbers (defined here as
a special type of locally symmetric spaces) and to the general theory
of canonical models; it also discusses in more detail the example of the
Siegel modular varieties. Lecture 2 presents some families of Shimura
varieties (PEL type, Hodge type, abelian type) and the results that are
known about their canonical and integral models.
Finally, lecture 3 discusses the cohomology of Shimura varieties,
concentrating mostly on the compact non-endoscopic case.

\end{abstract}

\maketitle


\bibliographystyle{amsplain}

\section{Lecture 1: locally symmetric spaces and Shimura varieties}

\subsection{Locally symmetric spaces}
\label{subsection_lss}

Let $G$ be a semi-simple algebraic group over $\Q$, for example
$\SL_n$, $\Sp_{2n}$, or a special orthogonal or special unitary group. 
Locally symmetric spaces for $G(\R)$ are ``nice enough'' spaces
whose cohomology is
related to automorphic representations of $G$.
A good reference for locally symmetric spaces is the introductory
paper \cite{Ji} by Ji.

To simplify the presentation, we will
assume here that $G(\R)$ is connected.
Let $K_\infty$ be a maximal compact subgroup of $G(\R)$, and let
$X=G(\R)/K_\infty$.
If $\Gamma$ is a discrete subgroup of $G(\R)$ such that
$\Gamma\sous G(\R)$ (or equivalently $\Gamma\sous X$) is compact
and that $\Gamma$ acts properly and freely on $X$,
\footnote{This holds for example if $\Gamma$ is torsion free.}
then there is a classical connection between the cohomology of
$\Gamma\sous X$ and automorphic representations of $G(\R)$, called
\emph{Matsushima's formula} (see Matsushima's paper \cite{Mat}).
We will state a modern reformulation in Lecture~3 (see
Theorem~\ref{thm_Matsushima}),
but roughly it relates the Betti numbers of $\Gamma\sous X$ and the
multiplicities of representations of $G(\R)$ in $L^2(\Gamma\sous G(\R))$.

In fact, Matsushima's paper deals with semi-simple real Lie groups. Here,
we have an algebraic group defined over $\Q$, so we have a particularly
nice way to produce discrete subgroups of $G(\R)$.
Remember that a subgroup $\Gamma$ of $G(\Q)$ is called an
\emph{arithmetic subgroup} if there exists 
a closed embedding of algebraic groups $G\subset\GL_N$ such that, setting
$G(\Z)=G(\Q)\cap\GL_N(\Z)$, we have that $\Gamma\cap G(\Z)$ is
of finite index in $\Gamma$ and in $G(\Z)$.
\footnote{We can check that this
definition does not depend on the embedding $G\subset\GL_N$,
see \cite[Proposition~4.2]{Ji}.}
If $\Gamma$ is small enough, then it acts properly and freely on $X$
(\cite[Proposition~5.5]{Ji}), so
the quotient $\Gamma\sous X$ is a real analytic manifold. Also,
the quotient $\Gamma\sous G(\R)$ is compact if and only if $G$ is anisotropic
(over $\Q$), which means that $G$ has no nontrivial parabolic subgroup
defined over $\Q$ (\cite[Theorem~5.10]{Ji}).
If $\Gamma\sous X$ is not compact but $G(\R)$ has a discrete
series, then there is an extension of
Matsushima's formula, due to Borel and Casselman in \cite{BC}, 
that involves $L^2$ cohomology of $\Gamma\sous X$; 
see~\ref{sssection_BC}.

We actually would like to see automorphic representations of $G(\Ade)$
(not just $G(\R)$) in the cohomology of our spaces, so we will use
adelic versions of $\Gamma\sous X$. Let $K$ be an open compact subgroup
of $G(\Af)$, where $\Af$ is the ring of finite ad\`eles; 
for example, if we have chosen an embedding $G\subset\GL_N$,
then we could take $K=G(\Af)\cap\Ker(\GL_N(\widehat{\Z})\ra\GL_N(\Z/n\Z))$,
for some positive integer $n$ (these are called \emph{principal congruence
subgroups}). Let
\[M_K=G(\Q)\sous X\times G(\Af)/K,\]
where the group $K$ acts by right translations on the factor
$G(\Af)$, and the group $G(\Q)$ acts by left translations on both factors
simultaneously. Choose a system of representatives $(x_i)_{i\in I}$ of the
finite
\footnote{See Lemmas~5.11 and~5.12 of \cite{Milne-SV}.} 
quotient $G(\Q)\sous G(\Af)/K$, 
and set $\Gamma_i=G(\Q)\cap
x_i Kx_i^{-1}$ for every $i\in I$. Then the $\Gamma_i$ are arithmetic subgroups
of $G(\Q)$, and we have
\[M_K=\coprod_{i\in I}\Gamma_i\sous X,\]
so $M_K$ is a real analytic manifold if $K$ is small enough.
But now we have an action of $G(\Af)$ on the projective system
$(M_K)_{K\subset G(\Af)}$, so we get an action on $\varinjlim_K\H^*(M_K)$,
where $\H^*$ is any ``reasonable'' cohomology theory, for example
Betti cohomology.
If $G$ is anisotropic over $\Q$, then Matsushima's result can be
reformulated to give a description of this action in terms
of irreducible representations of $G(\Ade)$ appearing in
$L^2(G(\Q)\sous G(\Ade))$; this is the version we give in
Theorem~\ref{thm_Matsushima}.
There is also a version of the
Borel-Casselman generalization (see~\ref{sssection_BC}).

There is another way to think about the action of
$G(\Af)$ on $(M_K)_{K\subset G(\Af)}$, which does not involve a limit
on $K$. Fix a Haar measure on $G(\Af)$ such that open compact subgroups
of $G(\Af)$ have rational volume (this is possible because these groups
are all commensurable); then every open subset of $G(\Af)$ has rational
volume, because it is a finite union of translates of open compact
subgroups.
The \emph{Hecke algebra} of $G$ is the space $\Hf_G$ of locally
constant functions with compact support
from $G(\Af)$ to $\Q$; if $f,g\in\Hf_G$, then the convolution
product $f*g$ still has rational values by the choice of Haar measure, so
convolution defines a multiplication on $\Hf_G$. For every open compact
subgroup $K$ of $G(\Af)$, the \emph{Hecke algebra at level $K$} is the
subalgebra $\Hf_{G,K}$ of bi-$K$-invariant functions in $\Hf_G$; we have
$\Hf_G=\bigcup_K\Hf_{G,K}$.

Fix $K$ small enough. Then $\H^*(M_K)$ is the set of
$K$-invariant vectors in $\varinjlim_{K'\subset G(\Af)}\H^*(M_{K'})$, so it
has an action of $\Hf_{G,K}$.
\footnote{In fact, we can recover the action of $G(\Af)$ on
$\varinjlim_{K\subset G(\Af)}\H^*(M_{K})$ from the action of $\Hf_{G,K}$ on
$M_K$ for every $K$ small enough.}
We can describe this action using Hecke correspondences:
\label{def_Hecke_correspondences}
let $g\in G(\Af)$ and let $K'$ be an open compact subgroup of $G(\Af)$ such
$K'\subset K\cap gKg^{-1}$. Then we have a \emph{Hecke correspondence}
$(T_1,T_g):M_{K'}\ra M_K^2$ sending the class of $(x,h)$ in $M_{K'}$ to that
of $((x,h),(x,hg))$ in $M_K\times M_K$,
and $T_1,T_g$ are both finite covering maps if
$K$ is small enough. 
Then the action of $\frac{1}{\vol(K')}\ungras_{KgK}$ on $\H^*(M_K)$ is given
by pulling back
cohomology classes along $T_1$, then pushing them forward along
$T_g$ (see Section~16 of~\cite{K-NP}).

We can also ask whether there is more structure on the spaces $\Gamma\sous X$
(or $M_K$). For example, suppose that
$G=\SL_2$ and $K_\infty=\SO(2)$. Then, for $\Gamma$ an arithmetic subgroup
of $\SL_2(\Z)$, the space $\Gamma\sous X$ is a modular curve, so it is
the set of complex points of an algebraic variety defined over a number
field $F$, and we can use the commuting actions of Hecke correspondences and
of the absolute Galois group of $F$ on the \'etale cohomology 
of this variety to construct
some instance of the global Langlands correspondence for
$\GL_2$. See for example Deligne's Bourbaki seminar \cite{Deligne-FML}.

In order to generalize this picture, 
as outlined in Langlands's papers \cite{Langlands-J} and
\cite{Langlands-M},
we first need to know when
the spaces $\Gamma\sous X$ or $M_K$ are
the set of $\C$-points of an algebraic variety, and whether this algebraic
variety is defined over a number field. 
As we will see in the section on canonical models (\ref{subsection_canonical_models}),
another advantage of $M_K$ over $\Gamma\sous X$
is that, when the answer to both questions above is ``yes'', then the $M_K$
for $K$ varying tend to all be defined over the same field, while this
is not the case for the $\Gamma\sous X$.

\begin{remark}
The first step is to check whether $\Gamma\sous X$ has the structure
of a complex manifold, and
there are obvious obstructions to that.
For example,
if $G=\SL_3$ and $K_\infty=\SO(3)$, then $\Gamma\sous X$ is $5$-dimensional
as a real manifold, so it cannot have the structure of a complex
manifold. 

\end{remark}

Choose a $G(\R)$-invariant Riemannian metric
on $X=G(\R)/K_\infty$ 
(such a metric is unique up to rescaling on each irreducible factor).
Then $X$ is a \emph{symmetric
space}, that is, a Riemannian manifold such that:
\begin{itemize}
\item[(a)] The group of isometries of $X$ acts transitively on $X$;
\item[(b)] For every $p\in X$, there exists a symmetry $s_p$ of $X$
(i.e. an involutive isometry) such that $p$ is an isolated fixed point
of $s_p$.

\end{itemize}
Moreover, the symmetric space $X$ is of \emph{noncompact type}, that is, it
has negative curvature. For $\Gamma$ a small enough arithmetic subgroup
of $G(\Q)$, the Riemannian manifold $\Gamma\sous X$ is a
\emph{locally symmetric space}; in particular, it does not satisfy
condition (a) anymore, and it satisfies a variant of condition (b)
where we only ask for the symmetry to be defined in a neighborhood of the
point. 
See Ji's notes \cite{Ji} for a review of locally
symmetric spaces.

We say that $X$ is a \emph{Hermitian symmetric domain} if it admits a
$G(\R)$-invariant Hermitian metric. See Section~1 of Milne's notes
\cite{Milne-SV} for a review of Hermitian symmetric domains. 

\begin{example}\label{ex_def_h_d}
Let $d$ be a positive integer. 
The \emph{Siegel upper half space} $\hgoth_d^+$ is the set of
symmetric $d\times d$ complex matrices $Y\in M_d(\C)$ such
$\im(Y)$ is positive definite;
if $d=1$, then this is just the usual upper half plane.
We claim that $\hgoth_d^+$ is a Hermitian symmetric
domain. For the proofs of the basic properties of $\hgoth_d^+$ that we state
below, a good reference is Siegel's paper \cite{Siegel}.

We first need to see $\hgoth_d^+$ as a symmetric space.
Let $\Sp_{2d}$ be the symplectic group of the symplectic form
with matrix $\begin{pmatrix} 0 & I_d \\ -I_d & 0\end{pmatrix}$,
where $I_d\in\GL_d(\Z)$ is the identity matrix.
For every commutative ring $R$, we have
\[\Sp_{2d}(R)=\left\{g\in\GL_{2d}(R)\mid {}^t g
\begin{pmatrix} 0 & I_d\\ -I_d & 0\end{pmatrix}g=
\begin{pmatrix} 0 & I_d\\ -I_d & 0\end{pmatrix}
\right\}.\]
Note that $\Sp_2=\SL_2$. We make $\Sp_{2d}(\R)$ act on $\hgoth_d^+$
by the following formula:
\[\begin{pmatrix} A & B \\ C & D\end{pmatrix}\cdot Y=
(A Y+B)(C Y+D)^{-1},\]
where $A,B,C,D$ are $d\times d$ matrices such that $\begin{pmatrix}
A & B \\ C & D\end{pmatrix}\in\Sp_{2d}(\R)$ (see page 9 of \cite{Siegel}).
Then this action is transitive (see page 10 of \cite{Siegel}).
Let $K_\infty$ be the stabilizer in
$\Sp_{2d}(\R)$ of $i I_d\in\hgoth_d^+$. Then $K_\infty=O(2d)\cap\Sp_{2d}(\R)$
(this is easy to check directly), 
so it is a maximal
compact subgroup of $\Sp_{2d}(\R)$, 
\footnote{In fact, we have an isomorphism $\U(d)\iso K_\infty$ sending
$X+iY\in\U(d)$ (with $X,Y\in\GL_d(\R)$) to $\begin{pmatrix}X & Y \\
-Y & X\end{pmatrix}$.}
and we have $\hgoth_d^+\simeq
\Sp_{2d}(\R)/K_\infty$ as real analytic manifolds.

Also, the space $\hgoth_d^+$ is an open subset of the complex
vector space of symmetric matrices in $M_d(\C)$, so it has an
obvious structure of complex manifold.
It remains to construct a
$\Sp_{2d}(\R)$-invariant Hermitian metric on $\hgoth_d^+$. 
Let $\Dcal_d$ be the set of symmetric matrices $A\in M_d(\C)$ such
that $I_d-A^*A$ is positive definite, this is a bounded domain in
the complex vector space of symmetric matrices in $M_d(\C)$, hence
is equipped with a canonical Hermitian metric called the
\emph{Bergman metric}, which has negative curvature
(see for example \cite[Theorem~1.3]{Milne-SV});
in particular, this metric is invariant by all holomorphic automorphisms
of $\Dcal_d$.
Now note that we have an isomorphism $\hgoth_d^+\iso\Dcal_d$ sending
$X\in\hgoth_d^+$ to $(i I_d-X)(i I_d+X)^{-1}$ (whose inverse sends
$A\in\Dcal_d$ to $i(I_d-A)(I_d+A)^{-1}$), see
pages 8-9 of \cite{Siegel}.
We can give a formula for the resulting Hermitian metric on $\hgoth_d^+$:
up to a positive scalar, it is given by
\[ds^2=\Tr(\im(Y)^{-2}dY\im(Y)^{-1}d\overline{Y})\]
(see formula (28) on page 17 of \cite{Siegel}).

The isomorphism $\hgoth_d^+\iso\Dcal_d$ is called a
\emph{bounded realization} of $\hgoth_d^+$.

\end{example}

We now state the classification of Hermitian symmetric domains
in terms of real algebraic groups. Let $U(1)=\{z\in\C\mid|z|=1\}$.

\begin{theorem}{\cite[Theorem~1.21]{Milne-SV}}\label{thm_HSD}
Suppose that $G(\R)$ is adjoint and connected.
The locally symmetric space $X$ is a Hermitian symmetric domain if and only
if there exists a morphism of real Lie groups $u:\U(1)\ra G(\R)$ such
that:
\begin{itemize}
\item[(a)] The only characters of $\U(1)$ that appear in its
representation on $\Lie(G(\R))$ via $\Ad\circ u$ 
are $1$, $z$ and $z^{-1}$;
\item[(b)] Conjugation by $u(-1)$ is a Cartan involution of
$G(\R)$, which means that the real Lie subgroup
$\{g\in G(\C)\mid g=u(-1)\overline{g}u(-1)^{-1}\}$ of $G(\C)$
is compact;
\item[(c)] The image of $u(-1)$ in every simple factor of $G(\R)$ is
nontrivial.

\end{itemize}

Moreover, we can choose $u$ such that $K_\infty$ is the centralizer
of $u$ in $G(\R)$, which means that $X$ is isomorphic to the set of conjugates
of $u$ by elements of $G(\R)$.

\end{theorem}

We explain the construction of the morphism $u$.
Suppose that $X$ is a Hermitian symmetric domain, and let $p\in X$.
For every $z\in\C$ with $|z|=1$, multiplication by $z$ on
$T_p X$ preserves the Hermitian metric and sectional curvatures, so
there exists a unique isometry $u_p(z)$ of $X$ fixing $p$ and such that
$T_p u_p(z)$ is multiplication by $z$. The uniqueness implies that
$u_p(z)u_p(z')=u_p(zz')$ if $|z|=|z'|=1$, so we get a morphism of
groups from $\U(1)$ to the group of isometries of $X$, which is equal
to $G(\R)^0_\ad=G(\R)$.

\begin{remark}
If we don't assume that $G(\R)$ is adjoint, then $u$ goes from $\U(1)$ to
$G(\R)_\ad$, as we see in the
sketch above, and it will not necessarily lift to a morphism from
$\U(1)$ to $G(\R)$.

\end{remark}

\begin{example}
\begin{itemize}
\item[(1)] If $G=\Sp_{2d}$, let $h:\C^\times\ra G(\R)$ be defined by
\[h(a+ib)=\begin{pmatrix}aI_d & -b I_d \\ b I_d & a I_d\end{pmatrix}.\]
Then we can take $u:\U(1)\ra\PSp_{2d}(\R)$ given by
$u(z)=h(\sqrt{z})$. Note that $u$ does not lift to a morphism from
$\U(1)$ into $G(\R)$.

\item[(2)] If $G=\PGL_n$ with $n\geq 3$, then the centralizer of a
character $u:\U(1)\ra G(\R)$ cannot be a maximal compact subgroup of
$G(\R)$ (exercise), so the locally symmetric space of maximal compact subgroups
of $G(\R)$ is not Hermitian.

\end{itemize}
\end{example}

Theorem~\ref{thm_HSD} puts some pretty strong restrictions on the
root systems of the simple factors of $G(\R)$, see
Theorem~1.25 of \cite{Milne-SV} and the table following it. In particular,
the type $A$ simple factors of $G(\R)$ must be of the form
$\PSU(p,q)$, and $G(\R)$ can have no simple factor of type
$E_8$, $F_4$ or $G_2$.

The natural next step would be to wonder for which Hermitian symmetric
domains $X$ the quotients $\Gamma\sous X$ are algebraic varieties, but in
fact it turns out that the answer is ``for all of them'', as was proved
by Baily and Borel. Borel later proved that this structure of algebraic variety
is unique.

\begin{theorem}[Baily-Borel and Borel, \cite{BB} and \cite{Bor-SM}]\label{thm_BB}
Suppose that $X=G(\R)/K_\infty$ is a Hermitian symmetric domain. Then, for
any torsion free arithmetic subgroup $\Gamma$ of $G(\Q)$, the quotient
$\Gamma\sous X$ has a unique structure of algebraic variety over $\C$
compatible with its structure of complex analytic manifold.

\end{theorem}

The very rough idea is that the sheaf of automorphic forms on
$\Gamma\sous X$ of sufficiently high weight will define an embedding
of $\Gamma\sous X$ into a projective space. Borel's uniqueness theorem
actually says that, if $V$ is a quasi-projective complex algebraic variety,
then any holomorphic map $V\ra X$ is regular; see the discussion below
Theorem~3.14 of~\cite{Milne-SV} for more details on its proof.

Remember that we did not just want the locally symmetric spaces
$\Gamma\sous X$ to be algebraic varieties, we also wanted them to be
defined over a number field, and we would ideally like the number field in
question to only depend on $G$ and $K_\infty$. For this, it will actually
be easier to work with reductive groups instead of semi-simple groups.
As a motivation for this, and for the definition of Shimura varieties, we
now spend some more time on the case of the symplectic group.

\subsection{The Siegel modular variety}
\label{subsection_Siegel_modular_variety}:

See the end of this subsection (\ref{ref_sss_background}
on page \pageref{ref_sss_background})
for some background on abelian schemes.

\subsubsection{Symplectic groups}

We use the notation of Example~\ref{ex_def_h_d}.
It is well-known that $\hgoth_1^+$ parametrizes elliptic curves over $\C$:
an element $\tau\in\hgoth_1^+$ is sent to the elliptic curve $E_\tau=
\C/(\Z+\Z\tau)$,
and $E_\tau\simeq E_{\tau'}$ if and
only if $\tau,\tau'\in\hgoth_1^+$
in the same orbit for the action of $\SL_2(\Z)=\Sp_2(\Z)$; so we can
recover $\tau$ from $E_\tau$ and the data of a symplectic isomorphism
$\H_1(E_\tau,\Z)\simeq\Z^2$ where $\Z^2$ is equipped with the standard
symplectic form. 
We want to give a similar picture for higher-dimensional abelian varieties;
in fact, the analogy works best if we consider abelian varieties with
a principal polarization (Definition~\ref{def_polarization}).

We first introduce some notation and
recall the definition of the (general) symplectic group as a group scheme
over $\Z$. 
If $R$ is a commutative ring, we denote by $\psi_R$ the perfect symplectic
pairing on $R^{2d}$ with matrix $\begin{pmatrix} 0 & I_d \\ -I_d & 0\end{pmatrix}$.
So we have
\[\psi_R((x_1,\ldots,x_d,y_1,\ldots,y_d),(x'_1,\ldots,x'_d,y'_1,\ldots,y'_d))=
\sum_{i=1}^d x_i y'_i-\sum_{i=1}^d x'_i y_i.\]
The \emph{general symplectic group} $\GSp_{2d}$ is
the reductive group scheme over $\Z$ whose
points over a commutative ring $R$ are given by:
\[\GSp_{2d}(R)=\{g\in\GL_{2d}(R)\mid\exists c(g)\in R^\times,\ 
\psi_R(g\cdot,g\cdot)=c(g)\psi_R\}.\]
The scalar $c(g)$ is called the \emph{multiplier} of $g\in\GSp_{2d}(R)$.
Sending $g$ to $c(g)$ defines a morphism of group schemes $c:\GSp_{2d}\ra
\GL_1$, whose kernel $\Sp_{2d}$ is called the \emph{symplectic group}.

\begin{example}
We have $\GSp_2=\GL_2$, $c=\det$ and $\Sp_2=\SL_2$.

\end{example}

\subsubsection{Complex abelian varieties}
\label{sssection_cav}

Let $A$ be complex abelian variety of dimension $d$; we identify $A$
and its set of complex points. Then
$A$ is a connected complex Lie group of dimension $d$, so we have
$A\simeq\Lie(A)/\Lambda$, with $\Lie(A)\simeq\C^d$
the universal cover of
$A$ and $\Lambda=\pi_1(A)=\H_1(A,\Z)\simeq\Z^{2d}$ a lattice in the underlying
$\R$-vector space. Let $A^\vee$ be the dual abelian variety, i.e. the
space of degree $0$ line bundles on $A$
(Definition~\ref{def_dual_av}).
We can identify
$\Lie(A^\vee)$ with the space of semi-linear forms on $\Lie(A)$ and
$\H_1(A^\vee,\Z)$ with the subspace $\Lambda^\vee$ of forms whose imaginary part
takes integer values on $\Lambda$ (see \cite{Mumford}, \S9). 
For every positive integer $n$, we have $A[n]=\frac{1}{n}\Lambda/\Lambda$ and
$A^\vee[n]=\frac{1}{n}\Lambda^\vee/\Lambda^\vee$, and the canonical
pairing $A[n]\times A^\vee[n]\ra\mu_n(\C)$ is given by
\[(v,u) \mapsto e^{-2i\pi n\im(u(v))}\]
(see \cite{Mumford}, \S24).
We then have a bijection
between the set of polarizations on $A$ (Definition~\ref{def_polarization})
and the set of positive definite
Hermitian forms\footnote{We take Hermitian forms to be semi-linear in the
first variable and linear in the second variable.}
$H$ on $\Lie(A)$ such that the symplectic form
$\im(H)$ takes integer values on $\Lambda$; given such a form $H$, the
corresponding isogeny $\lambda_H$ from $A$ to $A^\vee$ is given on $\C$-points
by:
\[\left\{\begin{array}{rcl}
\lambda_H:\Lie(A)/\Lambda &\ra&\Lie(A^\vee)/\Lambda^\vee\\
w & \mapsto & (v\mapsto H(v,w))\end{array}\right.\]
It follows that
the Weil pairing (see Remark~\ref{rmk_def_Weil_pairing}(2))
corresponding to $\lambda_H$ is the map
\begin{equation}\label{eq_Weil_pairing}
\left\{\begin{array}{rcl}A[n]\times A[n] & \ra & \mu_n(\C) \\
(v,w) & \mapsto & e^{-2i\pi n\im(H(v,w))}.
\end{array}\right.
\end{equation}
As $v,w\in\frac{1}{n}\Lambda$, we have $\im(H(v,w))\in\frac{1}{n^2}\Z$,
so the formula does define a $n$th root of $1$ in $\C$.

In particular, the polarization $\lambda_H$ is principal if and only
if $\Lambda$ is self-dual with respect to the symplectic form $\im(H)$,
that is,
\[\Lambda=\{w\in\Lie(A)\mid\forall v\in\Lambda,\ \im(H(v,w))\in\Z\}.\]
In that case, the symplectic $\Z$-module $(\Lambda,\im(H))$ is isomorphic
to $\Z^{2d}$ with the form $\psi_\Z$.

Let $\til{\Mcal}_d$ be the set of isomorphism classes
of triples $(A,\lambda,\eta_\Z)$, where
$A$ is a $d$-dimensional
complex abelian variety, $\lambda$ is a principal
polarization on $A$ and $\eta_\Z$ is an isomorphism of symplectic spaces
from $\H_1(A,\Z)$ to $(\Z^{2d},\psi_\Z)$.
We have an action of $\Sp_{2d}(\Z)$ on $\til{\Mcal}_d$:
if $c=(A,\lambda,\eta_\Z)\in\til{\Mcal}_d$ and $x\in\Sp_{2g}(\Z)$, set
$x\cdot c=(A,\lambda,x\circ\eta_\Z)$.

Let $(A,\lambda,\eta_\Z)\in
\til{\Mcal}_d$. Then $\Lambda=\H_1(A,\Z)$ is a lattice in the real
vector space $\Lie(A)$, we have $A=\Lie(A)/\Lambda$ and, if
$H_\lambda$ is the Hermitian form associated to the polarization $\lambda$, then
the symplectic form $\im(H_\lambda)_{\mid\Lambda}$ is sent to the form $\psi_\Z$ on
$\Z^{2d}$ by the isomorphism $\eta_\Z:\Lambda\iso\Z^{2d}$. 
So, if we see $\R^{2d}$ as a complex
vector space via the isomorphisms (of real vector spaces)
$\Lie(A)\iso\Lambda\otimes_\Z\R\iso\R^{2d}$, then the Hermitian form on
$\R^{2d}$ corresponding to $H_\lambda$ is
$(v,w)\mapsto\psi_\R(iv,w)+i\psi_\R(v,w)$.

This shows $\eta_\Z$ determines all the data of
the isomorphism class of $(A,\lambda,\eta_\Z)$,
except for the structure of complex vector space on $\R^{2d}\simeq\Lie(A)$.
This structure of complex vector space is equivalent to the data of
an $\R$-linear endomorphism $J$ of $\R^{2d}$ such that $J^2=-1$ (the
endomorphism $J$ corresponds to multiplication by $i$). We also
need the $\R$-bilinear map $\R^{2d}\times\R^{2d}\ra\C$ defined by
$(v,w)\mapsto\psi_\R(J(v),w)+i\psi_\R(v,w)$ to be
a positive definite Hermitian form on $\R^{2d}$, which is equivalent to the
following conditions:
\begin{itemize}
\item[(a)]
$\psi_\R(J(v),J(w))=\psi_\R(v,w)$ for all $v,w\in\R^{2d}$;
\item[(b)] the $\R$-bilinear form $(v,w)\mapsto\psi_\R(J(v),w)$ on $\R^{2d}$
(which is symmetric by (a)) is positive definite.
\end{itemize}

Conversely, if we have a complex structure $J$ on $\R^{2d}$ satisfying
(a) and (b), then we get a positive definite Hermitian form $H$ on
$\R^{2d}$ whose imaginary part takes
integer values on the lattice $\Z^{2d}$, so the complex torus
$A=\R^{2d}/\Z^{2d}$ has a polarization $\lambda$ induced by $H$, 
hence is an abelian variety (for example
by the Kodaira embedding theorem),
and we get an element $(A,\lambda,\id_{\Z^{2d}})$ of $\til{\Mcal}_d$.

To sum up, we have defined a bijection from $\til{\Mcal}_d$ to the set $X'$ of
endomorphisms $J$ of $\R^{2d}$ such that $J^2=-1$ and that 
$J$ satisfies conditions
(a) and (b). 

Now observe that, if $W$ is a $\R$-vector space, then the data
of an endomorphism $J$ of $W$ such that $J^2=-1$ (i.e. of the structure
of a $\C$-vector space on $W$) is equivalent to the data of a $\C$-linear
endomorphism $J_\C$ of $W\otimes_\R\C$ such that 
$W\otimes_\R\C=\Ker(J_\C-i\cdot\id_W)\oplus\Ker(J_\C+i\cdot\id_W)$ and
$\Ker(J_\C+i\cdot\id_W)=
\overline{\Ker(J_\C-i\cdot\id_W)}$, where $v\mapsto\overline{v}$
is the involution of
$W\otimes_\R\C$ induced by complex conjugation on $\C$. This is equivalent
to giving a $\C$-vector subspace $E$ of $W\otimes_\R\C$ such that
$W\otimes_\R\C=E\oplus\overline{E}$, i.e. a $d$-dimensional complex
subspace $E$ of $W\otimes_\R\C$ such that $E\cap\overline{E}=\{0\}$
(take $E=\Ker(J_\C-i\cdot\id_W)$).
\footnote{In fancier terms, we are saying that putting a structure of
complex vector space on $W$ is the same as putting a pure Hodge structure
of type $\{(-1,0),(0,-1)\}$ on it (or of type $\{(1,0),(0,1)\}$, depending
on your normalization). When $W=\R^{2d}$ and the complex structure comes
from an element $(A,\lambda,\eta_\Z)$ of $\til{\Mcal}_d$, then this
Hodge structure is the one coming from the Hodge
structure on $\H_1(A,\R)$ via the isomorphism
$\H_1(A,\R)\flnom{\eta_\Z\otimes\R}\Z^{2d}\otimes_\Z\R\iso\R^{2d}$.}

We apply this to $W=\R^{2d}$. Let $J$ be a complex structure on $\R^{2d}$, and
let $E$ be the corresponding $\C$-vector subspace of $W\otimes_\R\C=\C^{2d}$. 
Then condition
(a) on $J$ is equivalent to the fact that:
\begin{itemize}
\item[(a')] $\psi_\C(v,w)=0$ for all $v,w\in E$
\end{itemize}
(i.e. to the fact that $E$ is a Lagrangian subspace\footnote{A maximal
isotropic subspace.} of $\C^{2d}$),
and condition (b) on $J$ is equivalent to the fact that
\begin{itemize}
\item[(b')] 
$-i\psi_\C(v,\overline{v})\in\R_{>0}$ for all $v\in E\setminus\{0\}$. 
\end{itemize}
Note that these two conditions on a $\C$-vector subspace
$E$ of $\C^{2d}$ imply that $E\cap\overline{E}=\{0\}$. So we get a bijection
from $X'$ to the set of Lagrangian subspaces $E$ of $\C^{2d}$ 
satisfying (b').

If we represent Lagrangian subspaces of $\C^{2d}$ by their bases, seen
as complex matrices of size $2d\times d$, then the action of $\Sp_{2d}(\R)$
is just left multiplication.
For example, if we see $\R^{2d}$ as a complex vector space via
the isomorphism $\R^{2d}\iso\C^d$ sending
$(x_1,\ldots,x_d,y_1,\ldots,y_d)$ to $(x_1+i y_1,\ldots,x_d+iy_d)$, then
the corresponding endomorphism of $\R^{2d}$ is
$\begin{pmatrix} 0 & -I_d \\ I_d & 0\end{pmatrix}\in X'$, and the
corresponding subspace of $\C^{2d}$
is the one with basis $\begin{pmatrix} i I_d \\ I_d\end{pmatrix}$.

More generally,
if $Y\in\hgoth_d^+$, the subspace of $\C^{2d}$ with basis
$\begin{pmatrix}  Y \\ I_d\end{pmatrix}$ is a Lagrangian subspace satisfying 
condition (b'), and every such Lagrangian subspace is of that form. So we
get bijections $\til{\Mcal}_d\simeq X'\simeq\hgoth_d^+$,
and we can check that the second bijection is
$\Sp_{2d}(\R)$-equivariant. Unraveling the definitions, we see that
$Y\in\hgoth_d^+$ corresponds to the element $(A_Y,\lambda_Y,\eta_{\Z,Y})$
of $\til{\Mcal}_d$ such that
$A_Y=\C^d/(\Z^d+Y\Z^d)$, $\lambda_Y$ is the principal
polarization given by the Hermitian
form with matrix $\im(Y)^{-1}$ on $\C^d$, and $\eta_{\Z,Y}:\Z^d+
Y\Z^d\iso\Z^{2d}$ is the isomorphism sending $a\in\Z^d$ to
$(a,0)\in\Z^d\times\Z^d=\Z^{2d}$ and $Ya\in Y\Z^d$ to $(0,a)\in
\Z^d\times\Z^d=\Z^{2d}$.

Now we want an interpretation of the quotients $\Gamma\sous\hgoth_d^+$,
for $\Gamma$ an arithmetic subgroup of $\Sp_{2d}(\Q)$. We will do this
for the groups $\Gamma(n)=\Ker(\Sp_{2d}(\Z)\ra\Sp_{2d}(\Z/n\Z))$, where $n$
is a positive integer ($\Gamma(n)$ is called the \emph{principal congruence
subgroup} at level $n$). An arithmetic subgroup of $\Sp_{2d}(\Q)$ that contains
some $\Gamma(n)$ is called a \emph{congruence subgroup}. For $d=1$, there exist
arithmetic subgroups that are not congruence subgroups (see the introduction
of Raghunathan's paper \cite{Rag} for a counterexample, attributed to
Fricke-Klein);
for $d\geq 2$, it is known that every arithmetic subgroup is a congruence subgroup
(this was proved independently by Bass-Lazard-Serre in \cite{BLS} and
Mennicke in \cite{Men}).

We will need the notion of a level structure; we give the general definition
here.

\begin{definition}\label{def_level_structure}
Let $S$ be a scheme,
$(A,\lambda)$ be a principally polarized abelian scheme of relative
dimension $d$ over $S$, and $n$ be a positive integer. A
\emph{level $n$ structure} on $(A,\lambda)$ is a couple
$(\eta,\varphi)$, where 
$\eta:A[n]\iso\underline{\Z/n\Z}^{2d}_S$ and
$\varphi:\underline{\Z/n\Z}_S\iso\mu_{n,S}$ 
are isomorphisms of group schemes such that
$\varphi\circ\psi_{\Z/n\Z}\circ(\eta\times\eta)$ is the Weil pairing associated
to $\lambda$ on $A[n]$.

\end{definition}

\begin{remark}
A level $n$ structure on $(A,\lambda)$ can only exist if $n$ is invertible
on $S$ and $\mu_{n,S}$ is a constant group scheme.

\end{remark}

Note that isomorphisms $\varphi:\underline{\Z/n\Z}_S\iso\mu_{n,S}$
correspond to sections $\zeta\in\mu_n(S)$ generating $\mu_{n,S}$
(i.e. to primitive $n$th roots of $1$ over $S$), by sending
$\varphi$ to $\zeta=\varphi(1)$. So we will also see level
structures as couples $(\eta,\zeta)$, with $\zeta\in\mu_n(S)$ primitive.

Let $\zeta_n=e^{-2i\pi/n}\in\mu_n(\C)$. If $Y\in\hgoth_d^+$,
then $\frac{1}{n}\eta_{\Z,Y}$ defines an isomorphism of groups
$\eta_Y:A_Y[n]\iso(\Z/n\Z)^{2d}$, and it follows from formula
\eqref{eq_Weil_pairing} that $(\eta,\zeta_n)$ is a level
$n$ structure on $(A_Y,\lambda_Y)$.

Using the fact that $\Sp_{2d}(\Z)\ra\Sp_{2d}(\Z/n\Z)$ is surjective
for every $n\in\Nat$, which follows from strong approximation for
$\Sp_{2d}$,
\footnote{See the papers \cite{Plat1} and \cite{Plat2} of Platonov.}
we finally get:

\begin{prop}\label{prop_complex_unif}
Let $n$ be a positive integer.
The map $Y\mapsto (A_Y,\lambda_Y,\eta_Y)$ induces a bijection from
$\Gamma(n)\sous\hgoth_d^+$ to the set of isomorphism classes of
triples $(A,\lambda,\eta)$, where $(A,\lambda)$ is a principally
polarized complex abelian variety of dimension $d$ and
$\eta:A[n]\iso(\Z/n\Z)^{2d}$ is an isomorphism of groups such that
$(\eta,\zeta_n)$ is a level $n$ structure on $(A,\lambda)$.

\end{prop}

So we have reinterpreted $\Gamma(n)\sous\hgoth_d^+$ as the set of
$\C$-points of a certain moduli problem. Considering this moduli problem
over more general bases, we have an
``obvious'' way to make $\Gamma(n)\sous\hgoth_d^+$
into an algebraic variety, which we now review.

\subsubsection{The connected Siegel modular variety}

Let $\Of_n=\Z[1/n][T]/(T^n-1)$. If $S$ is a scheme over $\Of_n$,
we denote by $\varphi_0:\underline{\Z/n\Z}_S
\iso\mu_{n,S}$ the isomorphism sending $1$ to the class of $T$.

\begin{definition}\label{def_M'dn}
Let $\Mcal'_{d,n}$ be the contravariant
functor from the category of $\Of_n$-schemes to the
category of sets sending $S$ to the set of isomorphisms classes of
triples $(A,\lambda,\eta)$, 
where $(A,\lambda)$ is a principally polarized abelian scheme of relative
dimension $d$ over $S$ and $\eta:A[n]\iso\underline{\Z/n\Z}^{2d}_S$ is an
isomorphism of group schemes such that $(\eta,\varphi_0)$ is a level
$n$ structure on $(A,\lambda)$.

An isomorphism from $(A,\lambda,\eta)$ to $(A',\lambda',\eta')$
is an isomorphism of abelian varieties $u:A\iso A'$ such that
$\lambda'\circ u=u^\vee\circ\lambda$ and $\eta'=\eta\circ(u,u)$, and,
for $f:T\ra S$ a morphism
of $\Of_n$-schemes, the map
$\Mcal'_{d,n}(S)\ra\Mcal'_{d,n}(T)$ is base change along $f$.

\end{definition}

\begin{theorem}[Mumford, cf. \cite{CF}]\label{thm_Mdn'}
Suppose that $n\geq 3$.
Then the functor $\Mcal'_{d,n}$ 
is representable by a smooth quasi-projective $\Of_n$-scheme purely
of dimension $d(d+1)/2$ and with connected geometric fibers, which we
still denote by $\Mcal'_{d,n}$ and call the
\emph{connected Siegel modular variety of level $n$}.

\end{theorem}

\begin{remark}\label{rmk_small_level}
If $n\in\{1,2\}$, then triples $(A,\lambda,\eta)$ as in 
Definition~\ref{def_M'dn}
may have automorphisms, so we should see $\Mcal'_{d,n}$ as a stack.
It will then be representable
by a smooth Deligne-Mumford stack over $\Of_n$ that is a finite \'etale
quotient of the scheme $\Mcal'_{d,3n}$.

\end{remark}

We can now strengthen Proposition~\ref{prop_complex_unif} to the following
result.

\begin{prop}\label{prop_complex_unif2}
Let $n\geq 3$ be an integer.
Then the map $Y\mapsto (A_Y,\lambda_Y,\eta_Y)$ induces an isomorphism
of complex algebraic varieties from
$\Gamma(n)\sous\hgoth_d^+$ to $\Mcal'_{d,n}(\C)$.

\end{prop}

Let $\alpha$ be the map of the proposition. 
To show that $\alpha$ is regular, we use the fact that the triples
$(A_Y,\lambda_Y,\eta_Y)$ are obtained as the fibers of an abelian scheme
$\Acal\ra\Gamma(n)\sous\hgoth_d^+$ given together with a polarization $\lambda$
and a level $n$ structure $\eta$.
The abelian scheme $\Acal$ is a locally symmetric space for the non-reductive
group $\Q^{2d}\rtimes\Sp_{2d}$, where $\Sp_{2d}$ acts on $\Q^{2d}$ via its standard
representation, and is an example of a \emph{mixed Shimura variety} (see for example
Pink's thesis \cite{Pink-these} for the theory of mixed Shimura varieties and their
compactifications), and $\lambda$ and $\eta$ can also be given by explicit formulas
interpolating those for $\lambda_Y$ and $\eta_Y$.
Once we have this fact, the map $\alpha$ corresponds to the
$\Gamma(n)\sous\hgoth_d^+$-point of $\Mcal'_{d,n}$ defined by $(\Acal,\lambda,\eta)$,
hence it is a morphism of algebraic varieties.

Let us explain why $\alpha^{-1}$ is
holomorphic (for more details, see the proof of Proposition~1.1 in
Milne's~\cite{Mi-SVM}); by Borel's uniqueness theorem (see Theorem~\ref{thm_BB}),
this will imply that it is an isomorphism of algebraic varieties.
We don't know anything explicit about the
algebraic variety structure on $\Mcal:=\Mcal'_{d,n}(\C)$, 
but we do know that we have a universal abelian scheme $f:\Acal\ra\Mcal$
of relative dimension $d$, with a principal polarization $\lambda$ and a level
$n$ structure $\eta$. By definition, the map $\alpha^{-1}$ sends $x\in\Mcal$ to
the element of $\Gamma(n)\sous\hgoth_d^+$ corresponding to the triple 
$(\Acal_x,\lambda_x,\eta_x)$, so we need to understand how to construct this.
The exponential map gives a surjective morphism $\exp:\Lie(\Acal/\Mcal)\ra\Acal$ of
abelian groups schemes over $\Mcal$ (where the source is the relative Lie
algebra of $\Acal\ra\Mcal$), whose kernel is the relative homology of $\Acal\ra
\Mcal$ with coefficients in $\Z$. So we
get an exact sequence of abelian sheaves over $\Mcal$:
\[0\ra \Lf:=(R^1f_*\underline{\Z})^\vee\ra\Lie(\Acal/\Mcal)\flnom{\exp}\Acal\ra 0\]
Let $x\in\Mcal$. The level structure $\eta_x$ defines an isomorphism $(\Z/n\Z)^{2d}\iso\Acal_x[n]$,
that is, a basis of $\Acal_x[n]$; we can see the entries $e_1,\ldots,e_{2d}$
of this basis as elements
of $\frac{1}{n}\Lf_x/\Lf_x\subset\Lie(\Acal_x/\Mcal_x)/\Lf_x\iso\Acal_x$, and lift them
to a basis $\til{e}_1,\ldots,\til{e}_{2d}$
of $\frac{1}{n}\Lf_x$, which gives an isomorphism $\Z^{2d}\iso\Lf_x$, hence
an isomorphism $\R^{2d}\iso\Lie(\Acal_x/\Mcal_x)$. Now we go
back to the discussion in Subsection~\ref{sssection_cav} and we see that
the element $Y_x$ of $\hgoth_d^+$ whose class is equal to $\alpha^{-1}(x)$ is
basically a basis of the Lagrangian subspace corresponding to the complex
structure on $\R^{2d}$ coming from the isomorphism $\R^{2d}\iso\Lie(\Acal_x/\Mcal_x)$ that
we just defined. So we need to show that we can find a neighborhood $U$ of $x$ such that
$\til{e}_1,\ldots,\til{e}_{2d}$ extend to holomorphic sections of $\frac{1}{n}\Lf_{\mid U}\ra U$
that define a level $n$ structure on $\Acal_y$ at each point $y$ of $U$. The extensions
exist because the morphism $\frac{1}{n}\Lf\ra\Mcal$ is \'etale, and the
condition on the extensions is open because it only involves finite groupe schemes over
$\Mcal$.

In particular, we showed that $\Gamma(n)\sous\hgoth_d^+$ is the set
of complex points of an algebraic variety defined over the
number field $\Q(\zeta_n)$. Unfortunately, this number field depends
on the level $n$. The issue is that we need a fixed primitive $n$th
root of $1$ in order to define the moduli problem $\Mcal'_{d,n}$, so
we need to be over a basis where such a primitive $n$th root exists.
To fix this problem, we will allow the primitive $n$th root of $1$ to
vary.

\vspace{.3cm}

\subsubsection{The Siegel modular variety}
\label{sssection_SMV}

\begin{definition}\label{def_Mdn}
Let $n$ be a positive integer.
The \emph{Siegel modular variety} $\Mcal_{d,n}$ is the contravariant functor
from the category of $\Z[1/n]$-schemes to the category of sets sending
a scheme $S$ to the set of isomorphism classes of triples $(A,\lambda,\eta,
\varphi)$,
where $(A,\lambda)$ is a principally polarized abelian scheme of relative
dimension $d$ over $S$ and $(\eta,\varphi)$ 
is a level $n$ structure on $(A,\lambda)$.

An isomorphism from $(A,\lambda,\eta,\varphi)$ to $(A',\lambda',\eta',\varphi')$
is an isomorphism of abelian varieties $u:A\iso A'$ such that
$\lambda'\circ u=u^\vee\circ\lambda$ and $\eta'=\eta\circ(u,u)$.

\end{definition}

The group $\GSp_{2d}(\widehat{\Z})$ acts on $\Mcal_{d,n}$: 
if $g\in\GSp_{2d}(\widehat{\Z})$ and $(A,\lambda,\eta,\varphi)\in\Mcal_{d,n}(S)$,
then
\[g\cdot(A,\lambda,\eta,\varphi)=(A,\lambda,g\circ\eta,c(g)^{-1}\varphi).\]
The kernel of this action
is the group $K(n)=\Ker(\GSp_{2d}(\widehat{\Z})\ra\GSp_{2d}(\Z/n\Z))$.
If $n$ divides $m$, then we have a morphism
$\Mcal_{d,m}\ra\Mcal_{d,n}$ that forgets part of the level $m$ structure;
this morphism is representable finite \'etale, and in fact it is a
torsor under the finite group $K(n)/K(m)$.

We have the following variant of Theorem~\ref{thm_Mdn'}.

\begin{theorem}[Mumford, cf. \cite{CF}]\label{thm_Mdn}
Suppose that $n\geq 3$.
Then the functor $\Mcal_{d,n}$ 
is representable by a smooth quasi-projective $\Z[1/n]$-scheme purely
of dimension $d(d+1)/2$, which we
still denote by $\Mcal_{d,n}$ and call the
\emph{Siegel modular variety of level $n$}.

\end{theorem}

\begin{remark}\label{rmk_MdK}
The scheme $\Mcal_{d,n}$ is the Shimura variety for
$\GSp_{2d}$ with level $K(n)=\Ker(\GSp_{2d}(\widehat{\Z})\ra
\GSp_{2d}(\Z/n\Z))$, 
or rather its integral model.
If $K$ is an open compact subgroup of $\GSp_{2d}(\Af)$ that is small enough,
\footnote{For example, $K\subset K(N)$ with $N\geq 3$.}
then we can also define the Shimura variety $\Mcal_{d,K}$ with level
$K$, and even its integral model: 
choose $n$ such that $K(n)\subset K$. Then $K(n)$ is a normal subgroup
of $K$, so the group $K/K(n)\subset\GSp_{2d}(\Z/n\Z)$ acts on $\Mcal_{d,n}$,
and we set $\Mcal_{d,K}=\Mcal_{d,n}/(K/K(n))$. It is easy to check that
$\Mcal_{d,K,\Q}$ does not depend on the choice of $n$; however, with our
definition, the localization of $\Z$ over which $\Mcal_{d,K}$ is defined
depends on the choice of $n$, so this is not ideal.

In fact, for $K$ an open compact subgroup of $\GSp_{2d}(\Af)$,
we have a direct definition of a level $K$ structure on a principally
polarized abelian scheme (see Section~5 of Kottwitz's paper~\cite{K-NP}). For
$K$ small enough, the scheme $\Mcal_{d,K,\Q}$ is the moduli space of principally
polarized abelian schemes with level $K$ structure. In general,
this moduli space is representable by a Deligne-Mumford stack.
We can also define the moduli space over a localization of $\Z$, but
the primes that we invert depend on $K$; see the discussion in
Subsubsection~\ref{sssection_can_integral_models}.

\end{remark}

Let us explain the relationship between
$\Mcal_{d,n}$ and $\Mcal'_{d,n}$. We define a map 
$s:(\Z/n\Z)^\times\ra\GSp_{2d}(\Z/n\Z)$ by
$s(\alpha)=\begin{pmatrix}0 & \alpha I_d \\ I_d & 0\end{pmatrix}$;
note that $s$ is a section of the multiplier $c:\GSp_{2d}(\Z/n\Z)\ra
(\Z/n\Z)^\times$, and that it is not a morphism of groups.

\begin{prop}
The morphism $\Mcal'_{d,n}\times(\Z/n\Z)^\times\ra\Mcal_{d,n,\Of_n}$ sending
$((A,\lambda,\eta),\alpha)$ to $(A,\lambda,s(\alpha)\circ\eta,
\varphi_0\circ\alpha)$ (where we see $\alpha$ as an automorphism of
$\underline{\Z/n\Z}_S$ for any scheme $S$) is an isomorphism.

\end{prop}

As a corollary, we get a description of the complex points of $\Mcal_{d,n}$.
Let $\hgoth_d=\hgoth_d^+\cup(-\hgoth_d^+)$ be the set of symmetric
matrices $Y\in M_d(\C)$ such that $\im(Y)$ is
positive definite or negative definite. The action of $\Sp_{2d}(\R)$ on
$\hgoth_d$
extends to a transitive action of $\GSp_{2d}(\R)$, 
given by the same formula. The stabilizer of $iI_d\in\hgoth_d$ in
$\GSp_{2d}(\R)$ is $\R_{>0}K_\infty$, where $K_\infty$ is as before the
stabilizer of $iI_d$ in $\Sp_{2d}(\R)$,
so $\hgoth_d\simeq\GSp_{2d}(\R)/
\R_{>0}K_\infty$ as real analytic manifolds.

\begin{cor}
We have an isomorphism of complex algebraic varieties
\[\Mcal_{d,n}(\C)\simeq\GSp_{2d}(\Q)\sous(\hgoth_d\times\GSp_{2d}(\Af)/
K(n)\]
extending the isomorphism of Proposition \ref{prop_complex_unif2},
where $K(n)=\Ker(\GSp_{2d}(\widehat{\Z})\ra\GSp_{2d}(\Z/n\Z))$ and
$\GSp_{2g}(\Q)$ acts diagonally on $\hgoth_d\times\GSp_{2d}(\Af)$.

\end{cor}

This follows from the fact that
\[\GSp_{2d}(\Q)\sous(\hgoth_d\times\GSp_{2d}(\Af)/K(n)\simeq
\GSp_{2d}(\Q)^+\sous(\hgoth_d^+\times\GSp_{2d}(\Af)/K(n)),\]
where $\GSp_{2d}(\Q)^+=\{g\in\GSp_{2d}(\Q)\mid c(g)>0\}$, and from
strong approximation for $\Sp_{2d}$,
\footnote{See the papers \cite{Plat1} and \cite{Plat2} of Platonov.} 
which implies that $c$ induces a bijection
\[\GSp_{2d}(\Q)_+\sous\GSp_{2d}(\Af)/K(n)\iso\Q_{>0}\sous\Af^\times/c(K(n))
\simeq\widehat{\Z}^\times/(1+n\widehat{\Z})\simeq(\Z/n\Z)^\times.\]
For every $i\in(\Z/n\Z)^\times$, we choose $x_i\in\GSp_{2d}(\Af)$
lifting $i$ and we set
\[\Gamma(n)_i=\GSp_{2d}(\Q)^+\cap x_i K(n) x_i^{-1}=\Sp_{2d}(\Q)\cap
x_i K(n) x_i^{-1}.\]
Then the $\Gamma(n)_i$ are arithmetic subgroups of $\Sp_{2d}(\Q)$, 
and we have
\[\GSp_{2d}(\Q)\sous(\hgoth_d\times\GSp_{2d}(\Af)/K(n)\simeq
\coprod_{i\in(\Z/n\Z)^\times}\Gamma(n)_i\sous\hgoth_d^+\]
as complex manifolds. 

In fact, for $i\in(\Z/n\Z)^\times$, we can take $x_i=\begin{pmatrix} 0 & a_i
I_d \\ I_d & 0\end{pmatrix}$ with $a_i\in\widehat{\Z}^\times$ lifting $i$.
In particular, we have $x_i\in\GSp_{2d}(\widehat{\Z})$; as
$K(n)$ is a normal subgroup of $\GSp_{2d}(\widehat{\Z})$, we get
$x_i K(n) x_i^{-1}=K(n)$, hence $\Gamma(n)_i=\Gamma(n)$, and finally
\[\GSp_{2d}(\Q)\sous(\hgoth_d\times\GSp_{2d}(\Af)/K(n)\simeq
\coprod_{i\in(\Z/n\Z)^\times}\Gamma(n)\sous\hgoth_d^+.\]

\begin{remark}
If $K$ is a small enough open compact subgroup of $\GSp_{2d}(\Af)$, 
then we get an isomorphism of
complex manifolds:
\[\Mcal_{d,K}(\C)\simeq\GSp_{2d}(\Q)\sous(\hgoth_d\times\GSp_{2d}(\Af)/K).\]

\end{remark}

\subsubsection{Hecke correspondences}

We can also descend the Hecke correspondences of
Subsection~\ref{subsection_lss}
to morphisms of schemes.

We proceed as in
Section~3 of Laumon's paper~\cite{Laumon}. Let
$g\in\GSp_{2d}(\Af)$, and let $K,K'$ be small enough open compact
subgroups of $\GSp_{2d}(\Af)$ such that $K'\subset K\cap gK g^{-1}$.
We want to define finite \'etale morphisms $T_1,T_g:\Mcal_{d,K'}\ra
\Mcal_{d,K}$, and the Hecke correspondence associated to $(g,K,K')$ is the
couple $(T_1,T_g)$.

Choose
$n\geq 3$ such that $K(n)\subset K'$;
then $\Mcal_{d,K'}=\Mcal_{d,n}/(K'/K(n))$
and $\Mcal_{d,K}=\Mcal_{d,n}/(K/K(n))$.
The morphism $T_1$ should just forget part of the level structure: 
as $K'/K(n)\subset K/K(n)$, we have an obvious choice for
$T_1:\Mcal_{d,K'}\ra\Mcal_{d,K}$: we take the morphism induced by
the identity of $\Mcal_{d,n}$.

To define $T_g$, we first consider the following special case:
if $g\in M_{2d}(\widehat{\Z})\cap\GSp_{2d}(\Af)$, let
$x=(A,\lambda,\eta,\varphi)\in\Mcal_{d,n}(S)$. Let
$u$ be the endomorphism of $(\Z/n\Z)^{2d}$ with matrix $g$.
Then $T_g$ sends the class of $x$ in $\Mcal_{d,K'}(S)$ to
the class of $(A',\lambda',\eta',\varphi)\in\Mcal_{d,n}(S)$,
where $A'$ is the quotient of $A$ by the finite flat subgroup scheme
$\eta^{-1}(\underline{\Ker u})$ of $A[n]$ and $\lambda',\eta'$ are the morphisms
deduced from $\lambda,g\circ\eta$.

Note that, if $g=a\cdot I_{2d}$ with $a\in\widehat{\Z}\cap\Af^\times$ and
$K'=K$, then the morphism $T_g:\Mcal_{d,K}\ra\Mcal_{d,K}$ is an isomorphism.

Finally, for a general $g\in\GSp_{2d}(\Af)$, we write $g=g_0^{-1}g_1$ with
$g_0\in(\widehat{\Z}\cap\Af^\times)I_{2d}$ and $g_1\in M_{2d}(\widehat{\Z})\cap
\GSp_{2d}(\Af)$, and we set $T_g=T_{g_0}^{-1}\circ T_{g_1}$.

\begin{remark}
If we use instead the general
definition of a level $K$ structure from Section~5 of \cite{K-NP},
then it becomes much easier to define the Hecke correspondences on the
moduli problems; 
see Section~6 of loc. cit. 

\end{remark}

\begin{remark}
We have two ways to think of $\Mcal_{d,n}(\C)$: as an adelic double quotient
or as finite disjoint union of spaces $\Gamma(n)\sous\hgoth_d^+$, which
are locally symmetric spaces associated to the semi-simple group
$\Sp_{2d}$. The first description is more convenient to see the
action of adelic Hecke operators, and the second description is a bit more
concrete and has simpler combinatorics.
Note also that the complex manifold $\Gamma(n)\sous
\hgoth_d^+$ is isomorphic to the set of $\C$-points of the algebraic variety 
$\Mcal'_{d,n}$, but this algebraic variety is defined over the field
$\Q[T]/(T^n-1)$, which depends on $n$. On the other hand, the adelic double
quotient $\GSp_{2d}(\Q)\sous(\hgoth_d\times\GSp_{2d}(\Af)/K(n)$ is
isomorphic to the set of $\C$-points of the algebraic variety
$\Mcal_{d,n}$, which is defined over $\Q$ (and even $\Z[1/n]$). 
So if we want to consider
Shimura varieties as a projective system of algebraic varieties
over a number field, then it makes sense to use the adelic double quotients,
because they are all defined over the same field.

\end{remark}

\subsubsection{Background on abelian schemes}
\label{ref_sss_background}

Let $S$ be a scheme. We denote by $\Sch/S$ the category of $S$-schemes.

\begin{definition}
An \emph{abelian scheme} over $S$ is an
$S$-group scheme $A\ra S$ which is smooth and proper with geometrically
connected fibers. If $S$ is the spectrum of a field $k$, an
abelian scheme over $S$ is also called an \emph{abelian variety} over
$k$.

A \emph{morphism of abelian schemes} over $S$ is a morphism of $S$-group
schemes between abelian schemes over $S$.

\end{definition}

\begin{prop}
Let $A$ be an abelian scheme over $S$. Then
the $S$-group scheme $A\ra S$ is commutative.

\end{prop}

\begin{definition}\label{def_Picard_functor}
Let $A$ be an abelian scheme over $S$, and let $e:S\ra A$ be its zero
section. We consider the following two functors from $(\Sch/S)^\op$ to the
category of sets:
\begin{itemize}
\item[(a)] The functor $\Pic_{A/S,e}$ sending an $S$-scheme $T\ra S$ to the
set of isomorphism classes of couples $(\Lf,\varphi)$, where
$\Lf$ is an invertible sheaf on $A\times_S T$
and $\varphi:\Of_T\iso e_T^*\Lf$ is an isomorphism,
with $e_T=e\times_S T:T\ra A\times_S T$. 
An isomorphism from $(\Lf,\varphi)$ to $(\Lf',\varphi')$ is an isomorphism
of $\Of_T$-modules $\alpha:\Lf\iso\Lf'$ such that $(e_T^*\alpha)\circ\varphi=
\varphi'$;
\item[(b)] The subfunctor $\Pic_{A/S,e}^0$ of $\Pic_{A/S,e}$ sending
an $S$-scheme $T\ra S$ to the set of isomorphism classes of
couples $(\Lf,\varphi)$ as in (a) such that, for every point
$t$ of $T$, every smooth projective curve over the residue field
$\kappa(t)$ of $t$ and every morphism of $\kappa(t)$-schemes
$f:C\ra A\times_S t$, the line bundle $f^*(\Lf_{\mid A\times_S t})$ is
of degree $0$ on $C$.

\end{itemize}
\end{definition}

\begin{remark}\label{rmk_Pic}
\begin{itemize}
\item[(1)] The functor $\Pic_{A/S,e}$ can be made into a functor into the
category of abelian groups: if $T$ is an $S$-scheme and
$(\Lf,\varphi)$, $(\Lf',\varphi')$ represent elements of $\Pic_{A/S,e}(T)$,
their product is represented by $(\Lf\otimes_{\Of_T}\Lf',\varphi\otimes
\varphi')$, where $\varphi\otimes\varphi'$ is the isomorphism
\[\xymatrix@C=10pt{\Of_T\ar@{=}[r]&\Of_T\otimes_{\Of_T}\Of_T\ar[r]^-\sim_-
{\varphi\otimes\varphi'}&
(e_T^*\Lf)\otimes_{\Of_T}(e_T^*\Lf')
\ar@{=}[r] &e_T^*(\Lf\otimes_{\Of_T}\Lf').}\]
Moreover, for every $S$-scheme $T$, the set $\Pic_{A/S,e}^0(T)$ is a subgroup
of $\Pic_{A/S,e}(T)$.

\item[(2)] If $X\ra S$ is a scheme over $S$, then the relative Picard
functor $\Pic_{X/S}$ on $\Sch/S$ is the fppf sheafification of the functor
$T\mapsto\Pic(X\times_S T)$, where, for $Y$ a scheme, we denote by $\Pic(Y)$
the set of isomorphism classes of line bundles on $Y$ (that is an abelian
group for the tensor product); see \cite[\href{https://stacks.math.columbia.edu/tag/0D25}{Situation 0D25}]{stacks-project}.
We can also define a subfonctor $\Pic_{X/S}^0$ of $\Pic_{X/S}$ as in
Definition \ref{def_Picard_functor}.
By  \cite[\href{https://stacks.math.columbia.edu/tag/0D28}{Lemma 0D28}]{stacks-project}, if $A$ is an abelian scheme over $S$, then there is an isomorphism
of functors in abelian groups $\Pic_{A/S}\simeq\Pic_{A/S,e}$, inducing an
isomorphism $\Pic_{A/s}^0\simeq\Pic_{A/S,e}^0$.

\item[(3)] We can upgrade $A\mapsto\Pic_{A/S,e}$ and $A\mapsto\Pic_{A/S,e}^0$
to contravariant functors in $A$: if $f:A\ra B$ is a morphism of abelian
schemes over $S$, then it induces a natural transformation
$f^*:\Pic_{B/S,e}\ra\Pic_{A/S,e}$ sending $(\Lf,\varphi)$ to
$(f^*(\Lf),f^*(\varphi))$, and $f^*$ sends $\Pic_{B/S,e}^0$ to
$\Pic_{A/S,e}^0$.

\end{itemize}
\end{remark}

\begin{theorem}
Let $A$ be an abelian scheme over $S$. Then $\Pic_{A/S,e}^0$ is representable
by an abelian scheme over $S$.

\end{theorem}

\begin{proof}
We know that $\Pic_{A/S,e}^0$ is representable by an algebraic space over $S$ by
a result of M. Artin (see \cite{Artin} or  \cite[\href{https://stacks.math.columbia.edu/tag/0D2C}{Proposition 0D2C}]{stacks-project}).
We can check on the moduli problem that this algebraic space is
proper and smooth, and its fibers over points of $S$ are abelian
varieties by the classical theory of the dual abelian variety (see
sections II.8 and III.13 of Mumford's book \cite{Mumford}).
It remains to prove that the algebraic space representing
$\Pic_{A/S,e}^0$ is a scheme; this is due to Raynaud, and a proof
is given in Theorem 1.9 of \cite{CF}.

\end{proof}

\begin{definition}\label{def_dual_av}
Let $A$ be an abelian scheme over $S$. The abelian scheme over $S$
representing $\Pic_{A/S,e}^0$ is called the \emph{dual abelian scheme}
of $A$ and denoted by $A^\vee$. In particular, we get a couple
$(\Pf_A,\varphi_A)$ representing the element of
$\Pic_{A/S,e}^0(A^\vee)$ corresponding to $\id_{A^\vee}$, with
$\Pf_A$ a line bundle on $A\times_S A^\vee$, called the 
\emph{Poincar\'e line bundle}.

\end{definition}

If $f:A\ra B$ is a morphism of abelian schemes over $S$, we denote
by $f^\vee:B^\vee\ra A^\vee$ the morphism corresponding to the natural
transformation $f^*:\Pic_{B/S,e}^0\ra\Pic_{A/S,e}^0$ of Remark
\ref{rmk_Pic}.

\begin{remark}
Let $e:S\ra A^\vee$ be the unit section. Then the pullback of
$\Pf_A$ by $A\times_S e:A=A\times_S S\ra A\times_S A^\vee$ is the
line bundle on $A$ corresponding to the element $e$ of
$A^\vee(S)=\Pic_{A/S,e}^0(S)$; in other words, it is isomorphic
to the trivial line bundle $\Of_A$. So $\Pf_A$ defines an element
of $\Pic_{A^\vee/S,e}^0(A)$, that is, a morphism of $S$-schemes
$A\ra A^{\vee\vee}$, called the \emph{biduality morphism}.
The \emph{biduality theorem} says that the biduality morphism is
an isomorphism. For $S$ the spectrum of a field, this is proved in
Section III.13 of \cite{Mumford}, and the general case reduces to this
by looking at the fibers over points of $S$.

\end{remark}

Let $A$ be an abelian scheme over $S$ and
let $\Lf$ be a line bundle on $A$. We denote by $\mu,p_1,p_2,\varepsilon:
A\times_S A\ra A$ the addition morphism, the first projection, the second
projection and the zero morphism respectively. Then the line
bundle $(\mu^*\Lf)\otimes(p_1^*\Lf^{\otimes -1})\otimes(p_2^*\Lf^{\otimes -1})\otimes
(\varepsilon^*\Lf)$ on $A\times_S A$ is trivial when restricted
to $S\times_S A$ via the zero section of $A$, hence it defines an
element of $\Pic_{A/S,e}^0(A)$, corresponding to a morphism of
$S$-schemes $\lambda(\Lf):A\ra A^\vee$, which is a morphism
of abelian varieties by the theorem of the cube
(see for example Section III.10 of \cite{Mumford}).

\begin{definition}\label{def_polarization}
Let $A$ be an abelian scheme over $S$. A \emph{polarization} on $A$
is an isogeny $\lambda:A\ra A^\vee$ (i.e. a finite faithfully flat morphism
of abelian schemes)
such that, for
every algebraically closed field $k$ and every morphism
$\Spec k\ra S$, the morphism $\lambda\times_S\Spec k:A\times_S\Spec k\ra
A^\vee\times_S\Spec k=(A\times_S\Spec k)^\vee$ is of the form
$\lambda(\Lf)$, for $\Lf$ an ample line bundle on $A\times_S\Spec k$.
We say that a polarization is \emph{principal} if it is an isomorphism.

A \emph{principally polarized abelian scheme} over $S$ is a pair
$(A,\lambda)$, where $A$ is an abelian scheme over $S$ and
$\lambda$ is a principal polarization on $A$.

\end{definition}

\begin{remark}\label{rmk_def_Weil_pairing}
Let $\lambda$ be a polarization on $A$, and let $n$ be a positive
integer.
Then, composing
$\lambda:A[n]\ra A^\vee[n]$ with the canonical Cartier pairing
$A[n]\times A^\vee[n]\ra\mu_{n,S}$,
\footnote{See page 183 of Mumford's book~\cite{Mumford} for the case
of abelian varieties over a field. For the case of general abelian schemes, see
the beginning of Section~1 of Oda's paper~\cite{Oda}, applied to the isogeny
given by multiplication by $n$.}
we get a pairing $A[n]\times A[n]\ra\mu_{n,S}$, called the
\emph{Weil pairing} associated to $\lambda$. If $\lambda$ is principal.
this is a perfect pairing.

\end{remark}

\subsection{Shimura varieties over $\C$}

Remember the upshot of Subsection~\ref{subsection_Siegel_modular_variety}:
if we want algebraic varieties that are all defined over the same
number field, and Hecke correspondences that are also defined on this
number field,
it is better to work with adelic double quotients
for a reductive group such as $\GSp_{2d}$ rather than with locally symmetric
spaces for a semi-simple group such as $\Sp_{2d}$. This 
(and Theorem~\ref{thm_HSD}) motivates the
definition of Shimura data, due to Deligne in \cite{Deligne-SV}.

\subsubsection{The Serre torus}
\label{sssection_Se}

Let $\Se$ be $\C^\times$ seen as an algebraic group over $\R$; this is
called the \emph{Serre torus}. In other
words, the group $\Se$ is the Weil restriction of scalars from $\C$ to
$\R$ of $\GL_1$, so that $\Se(R)=(R\otimes_\R\C)^\times$ for every $\R$-algebra
$R$. We denote by $w$ the injective morphism $\GL_{1,\R}\ra\Se$ corresponding to
the inclusion $\R^\times\subset\C^\times$.

We have $\Se(\C)=(\C\otimes_\R\C)^\times\iso\C^\times\times\C^\times$, where
the isomorphism sends $a\otimes 1+b\otimes i$ to $(a+ib,a-ib)$. 
So the abelian group
$\Hom(\Se_\C,\GL_{1,\C})$ of characters of $\Se$ is free of rank $2$ and
generated by the characters $z$ and $\overline{z}$ corresponding to the
two projections of $\C^\times\times\C^\times$ on $\C^\times$. We denote by
$r:\GL_{1,\C}\ra\Se_\C$ the injective morphism corresponding to the injection
of the first factor in $\C^\times\times\C^\times$.
\label{def_r}

If $V$ is a real vector space and $\rho:\Se\ra\GL(V)$ is a morphism
of algebraic groups (i.e. a representation of $\Se$ on $V$), then
we have $V_\C:=V\otimes_\R\C=\bigoplus_{p,q\in\Z}V^{p,q}$, where $V^{p,q}$ is the
subspace of $V_\C$ on which $\Se_\C$ acts by the character $z^{-p}
\overline{z}^{-q}$;
moreover, as $\rho$ is defined over $\R$, we have $\overline{V^{p,q}}=V^{q,p}$
for all $p,q\in\Z$. Let $m\in\Z$.
We say that $\rho$ is \emph{of weight $m$} if
$\rho\circ w:\GL_1\ra\GL(V)$ is equal to $x\mapsto x^{-m}\id_V$.

\begin{remark}
If $\rho:\Se\ra\GL(V)$ is of weight $m$,
we have $V^{p,q}=0$ unless $p+q=m$, so the decomposition $V_\C=\bigoplus_{p,q}
V^{p,q}$ is a pure Hodge structure of weight $m$ on $V$. In fact, representations
of weight $m$ of $\Se$ on $V$ are in bijection with pure Hodge structures of
weight $m$ on $V$.

\end{remark}

\subsubsection{Shimura data}

\begin{definition}\label{def_SD}
A \emph{Shimura datum} is a couple $(G,h)$, where $G$ is a connected
reductive algebraic group over $\Q$ and $h:\Se\ra G_\R$ is a morphism
of real algebraic groups such that:
\begin{itemize}
\item[(a)] The image of $h\circ w:\GL_{1,\R}\ra G_\R$ is central;
\item[(b)] If $\ggoth=\Lie(G_\R)$ and $\ggoth_\C=\bigoplus_{p,q\in\Z}
\ggoth^{p,q}$ is the decomposition induced by the representation
$\Ad\circ h:\Se\ra\GL(\ggoth)$, then we have $\ggoth^{p,q}=0$
unless $(p,q)\in\{(-1,1),(0,0),(1,-1)\}$;
\item[(c)] Conjugation by $h(i)$ induces a Cartan involution of
$G_\der(\R)$ (see Theorem~\ref{thm_HSD});
\item[(d)] $G_\ad$ has no normal subgroup (defined over
$\Q$) whose group of $\R$-points
is compact.
\footnote{Note that $G_\ad(\R)$ could still have compact
normal algebraic subgroups, as long as they are not defined
over $\Q$.}

\end{itemize}
\end{definition}

Note that condition (b) implies condition (a), because it implies that
the image of $h\circ w$ centralizes $\Lie(G_\R)$.

Let $(G,h)$ be a Shimura datum. We denote by $K_\infty$ the centralizer
of $h$ in $G(\R)$ and by $X$ the set of $G(\R)$-conjugates of $h$.
Then $K_\infty$ contains the center of $G(\R)$, and $K_\infty\cap G_\der(\R)^0$
is equal to the centralizer of $h(i)$ in $G_\der(\R)^0$, hence
is a maximal compact subgroup of $G_\der(\R)^0$ by condition (c). We
have $X\simeq G(\R)/K_\infty$, and Theorem~\ref{thm_HSD} implies that
there is a $G(\R)$-invariant complex structure on $X$ such that the
connected components of $X$ are Hermitian symmetric domains.

\begin{remark}
There are many variants in the definition of a Shimura datum. The definition
we gave here is equivalent to Definition~5.5 of Milne's notes~\cite{Milne-SV},
but it is possible to impose additional conditions to make some
statements simpler; see for example the discussion of additional
axioms on page~63 of~\cite{Milne-SV}. On the other hand, some authors
use a slightly more general definition, where $X$ is allowed to be
finite over the set of morphisms $\Se\ra G_\R$. This is for example
the case of Pink in~\cite{Pink-these} (see Definition~2.1)
and~\cite{Pink-l-adic} (see~(3.1)); in the first of these
references, Pink also generalizes the definition to include
groups that are not necessarily reductive.

\end{remark}

\begin{example}\label{ex_SD_GSp}
Take $G=\GSp_{2d}$. Up to conjugation, there exists a unique
morphism $h:\Se\ra\GSp_{2d}$ satisfying conditions (a)-(c) of
Definition~\ref{def_SD} and such that $(h\circ w)(x)=xI_{2d}$ for every
$x\in\R^\times$. An element of that class is given by
\[h(a+ib)=\begin{pmatrix}aI_d & -b I_d \\ b I_d & a I_d\end{pmatrix}.\]
For this $h$, we have $K_\infty=\GSp_{2d}(\R)\cap\GO(2d)$, and we
can check that the map $\GU(d)\ra K_\infty$ sending
$X+iY\in\GU(d)$ (with $X,Y\in M_d(\R)$) to $\begin{pmatrix} X & Y  \\
- Y &  X \end{pmatrix}$ is an isomorphism of Lie groups. So
$K_\infty=\R_{>0}K'_\infty$, where $K'_\infty=\Sp_{2d}(\R)\cap\Orth(d)$ is
the maximal compact subgroup of $\Sp_{2d}(\R)$ that was called $K_\infty$
in Subsections~\ref{subsection_lss} and~\ref{subsection_Siegel_modular_variety}.
This implies that $X\simeq\hgoth_d$.

The couple $(\GSp_{2d},h)$ is called a \emph{Siegel Shimura datum}.

\end{example}

Let $K$ be an open compact subgroup of $G(\Af)$. We set
\[M_K(G,h)(\C)=G(\Q)\sous X\times G(\Af)/K,\]
where the group $K$ acts by right translations on the factor
$G(\Af)$, and the group $G(\Q)$ acts by left translations on both factors
simultaneously. 
This is the \emph{Shimura variety at level $K$} associated to the Shimura
datum $(G,h)$.

As in Subsection~\ref{subsection_lss}, if
$(x_i)_{i\in I}$ is a system of representatives of the
finite quotient $G(\Q)\sous G(\Af)/K$, and if $\Gamma_i=G(\Q)\cap
x_i Kx_i^{-1}$ for every $i\in I$, then the $\Gamma_i$ are arithmetic subgroups
of $G(\Q)$, and we have
\[M_K(G,h)(\C)=\coprod_{i\in I}\Gamma_i\sous X.\]
Hence it follows from Theorem~\ref{thm_BB} that $M_K(G,h)(\C)$ is
the set of complex points of a quasi-projective algebraic variety over $\C$,
smooth if $K$ is small enough.

Again as in Subsection~\ref{subsection_lss}, we have Hecke correspondences
between the varieties $M_K(G,h)(\C)$, which are finite maps, hence morphisms of
algebraic varieties. This defines an action of $G(\Af)$ on the projective
system $(M_K(G,h)(\C))_{K\subset G(\Af)}$, or on its limit
\[M(G,h)(\C)=\varprojlim_K M_K(G,h)(\C).\]
Under the additional assumption that the $\R$-split and $\Q$-split
maximal subtori of the center of $G$ are equal (which is sometimes required
in the definition of a Shimura datum), the limit is equal to
$G(\Q)\sous X\times G(\Af)$
(see Proposition~4.19 of~\cite{Milne-SV}).
Moreover, for every open compact subgroup
$K$ of $G(\Af)$,
we have $M_K(G,h)(\C)=M(G,h)(\C)/K$.

\vspace{.3cm}

\subsubsection{Morphisms of Shimura varieties}
\label{sssection_morphisms_SV}

Let $(G_1,h_1)$, $(G_2,h_2)$ be Shimura data, and let $u:G_1\ra
G_2$ be a morphism of algebraic groups such that $u\circ h_1$
and $h_2$ are conjugated under $G_2(\R)$; we say that $u$ is
a \emph{morphism of Shimura data}.
Then $u$ induces a morphism
of complex manifolds $X_1\ra X_2$. Hence, for all $K_1\subset G_1(\Af)$,
$K_2\subset G_2(\Af)$ open compact subgroups such that $u(K_1)\subset K_2$,
we get a morphism of quasi-projective varieties $u(K_1,K_2):M_{K_1}(G_1,K_1)(\C)
\ra M_{K_2}(G_2,h_2)(\C)$. We can also think of the collection of all
$u(K_1,K_2)$ as a morphism of
$\C$-schemes $u:M(G_1,h_1)(\C)\ra M(G_2,h_2)(\C)$.

\begin{prop}[Deligne, see Proposition~1.15 of~\cite{Deligne-SV}]
If $G_1$ is an algebraic subgroup of $G_2$ and $u$ is the inclusion, then, 
for every open compact subgroup
$K_1$ of $G_1(\Af)$, there exists an open compact subgroup
$K_2\supset K_1$ of $G_2(\Af)$ such that $u(K_1,K_2):M_{K_1}(G_1,K_1)(\C)
\ra M_{K_2}(G_2,h_2)(\C)$ is a closed immersion.

\end{prop}

\subsubsection{Connected components}

For $(G,h)$ equal to the Shimura datum of Example~\ref{ex_SD_GSp} and
$K=\Ker(\GSp_{2d}(\widehat{\Z})\ra\GSp_{2d}(\Z/n\Z))$, we have seen that the
multiplier $c:\GSp_{2d}\ra\GL_1$ induces a bijection
\[\pi_0(M_K(G,h)(\C))\iso\Q_{>0}\sous\Af^\times/c(K)=
\Q^\times\sous\Ade^\times/c(K_\infty K).\]

In fact, it follows from real approximation and the Hasse principle that
this works for many Shimura varieties:

\begin{theorem}[Deligne, see \cite{Deligne-SV} 2.7]\label{thm_CC_SV}
Let $\nu:G\ra T:=G/G_\der$ be the quotient morphism, and
suppose that $G_\der$ is simply connected.
Then, for every open compact subgroup $K$ of $G(\Af)$,  
the map $\nu$ induces a bijection
\[\pi_0(M_K(G,h)(\C))\iso T(\Q)\sous T(\Ade)/\nu(K_\infty\times K).\]

\end{theorem}

\subsection{Canonical models}
\label{subsection_canonical_models}

In the situation of Example~\ref{ex_SD_GSp}, we have seen that the
algebraic varieties $M_K(G,h)(\C)$ and all the Hecke correspondences are
defined over $\Q$. We would like to generalize this kind of result to
other Shimura varieties.

\subsubsection{Model of a Shimura variety}
First we need to say what we mean by a model. 

\begin{definition}\label{def_model}
Let $(G,h)$ be a Shimura datum,
and $F$ be a subfield of $\C$. A \emph{model} of the projective
system $(M_K(G,h)(\C))_K$ over $F$ is the data:
\begin{itemize}
\item for every open compact subgroup $K$ of $G(\Af)$, of a quasi-projective
variety $M_K$ over $F$ and an isomorphism
$\iota_K:M_K\otimes_F\C\iso M_K(G,h)(\C)$;
\item for every $g\in G(\Af)$ and all open compact subgroups
$K,K'$ of $G(\Af)$ such that $gK'g^{-1}\subset K$, of a morphism of
$F$-varieties $T_{g,K,K'}:M_{K'}\ra M_K$,

\end{itemize}
such that:
\begin{itemize}
\item[(i)] For all $g,K,K'$ as above, the morphism $\iota_K\circ T_{g,K,K',\C}
\circ\iota_{K'}^{-1}:M_{K'}(G,h)(\C)\ra M_K(G,h)(\C)$ sends the class
of $(x,h)$ in $M_{K'}(G,h)(\C)$ to the class of $(x,hg)$ in $M_K(G,h)(\C)$;

\item[(ii)] If $K$ is an open compact subgroup of $G(\Af)$ and
$g\in K$, then $T_{g,K,K}=\id_{M_K}$;

\item[(iii)] If $K,K',K''$ are open compact subgroups of $G(\Af)$ and
$g,h\in G(\Af)$ are such that $gK'g^{-1}\subset K$ and $hK''h^{-1}\subset K'$,
then $T_{g,K,K'}\circ T_{h,K',K''}=T_{gh,K,K''}$;

\item[(iv)] If $K,K'$ are open compact subgroups of $G(\Af)$ such that
$K'$ is a normal subgroup of $K$, then the morphisms $T_{g,K',K'}$ for
$g\in K$ define an action of $K/K'$ on $M_{K'}$ (this follows from
(ii) and (iii)), and $T_{1,K',K}:M_{K'}\ra M_K$ induces an isomorphism
$M_{K'}/(K/K')\ra M_K$.

\end{itemize}

\end{definition}

\begin{remark}\label{remark_CN}
We phrased Definition~\ref{def_model} in this way to stay close to the
definition of a canonical model given by Deligne in D\'efinition~3.1
of~\cite{Deligne-SV}. In Deligne's words, a family of $E$-schemes $(M_K)$ with
a family of morphisms $(T_{g,K,K'})$ satisfying conditions (ii), (iii) and (iv)
is a \emph{scheme over $E$ with a continuous action of $G(\Af)$}, and the
family of isomorphisms $(\iota_K)$ satisfying condition (i) is \emph{an isomorphism
of $(M_K)\otimes_F\C$ with $M(G,h)(\C)$ compatible with the $G(\Af)$-actions}.
However, we did not not phrase the conditions in the most optimal way possible;
for example, condition (i) implies conditions (ii) and (iii).

\end{remark}

If we have a model $(M_K)_K$ of $(M_K(G,h)(\C))_K$ over $F$, we write
$M=\varprojlim_K M_K$ (where the transition morphisms are given by the
$T_{1,K',K}$). This is an $F$-scheme with an action of $G(\Af)$, and we have
a $G(\Af)$-equivariant isomorphism $M\otimes_F\C\iso M(G,h)(\C)$.

In particular
we get an action of $\Gal(\overline{F}/F)$ on $\pi_0(M(G,h)(\C))\iso
\pi_0(M\otimes_F\overline{F})$, which must commute with the action
of $G(\Af)$. Under the hypothesis of Theorem~\ref{thm_CC_SV},
we have $\pi_0(M\otimes_F\overline{F})\iso\pi_0(T(\Ade)/T(\Q))/\pi_0(K_\infty)$
with $T=G/G_\der$,
and $G(\Af)$ acts transitively on this set of connected components
(Proposition~2.2 of \cite{Deligne-SV}). So every element of $\Gal(\overline{F}/F)$
acts by translation by an element of $\pi_0(T(\Ade)/T(\Q))/\pi_0(K_\infty)$,
and the action of $\Gal(\overline{F}/F)$ comes from a morphism of groups
$\Gal(\overline{F}/F)\ra\pi_0(T(\Ade)/T(\Q))/\pi_0(K_\infty)$, that necessarily
factors through the maximal abelian quotient $\Gal(\overline{F}/F)^\ab$.

Suppose that $F$ is a number field. Then global class field theory
\footnote{Normalized so that local uniformizers correspond to geometric
Frobenius elements.}
gives an isomorphism
\[\Gal(\overline{F}/F)^\ab\iso\pi_0(F^\times\sous\Ade_F^\times)\]
where $\Ade_F$ is the ring of adeles of $F$, so the action of $\Gal(
\overline{F}/F)$ on $\pi_0(M\otimes_F\overline{F})$ comes from a morphism
of groups
\[\lambda_M:\pi_0(F^\times\sous\Ade_F^\times)\ra\pi_0(T(\Q)\sous T(\Ade))/
\pi_0(K_\infty),\]
called the \emph{reciprocity law} of the model.

\subsubsection{The case of tori}
\label{sssection_SV_tori}

We consider the case where $G=T$ is a torus. Let
$h:\Se\ra T_\R$ be any morphism of real algebraic groups. Then $h$ trivially
satisfies the conditions of Definition~\ref{def_SD}, so we get a Shimura
datum $(T,h)$, and $X=T(\R)/\Cent_{T(\R)}(h)$ is a singleton.
For every open compact subgroup $K$ of $T(\Af)$,
\[M_K(T,h)(\C)=T(\Q)\sous T(\Af)/K\]
is a finite group. Let $\overline{T(\Q)}$ be the closure of $T(\Q)$ in
$T(\Af)$; then
\[M(T,h)(\C)=\overline{T(\Q)}\sous T(\Af),\]
which is a profinite group isomorphic to $\pi_0(T(\Q)\sous T(\Af))$ by the
obvious map.
Giving a model of the Shimura variety of $(T,h)$ over
a subfield $F$ of $\C$ is the same as giving an action of
$\Gal(\overline{F}/F)$ over $\pi_0(T(\Q)\sous T(\Af))$ (commuting with the action
of $T(\Af)$ by translations), i.e. a morphism of groups
$\Gal(\overline{F}/F)\ra \pi_0(T(\Q)\sous T(\Af))$.
If $F$ is a number field, this is equivalent to giving 
a morphism of groups
\[\lambda_M:\pi_0(F^\times\sous\Ade_F^\times)\ra\pi_0(T(\Q)\sous T(\Ade)).\]

It is natural to construct such a morphism from a morphism of algebraic
groups $F^\times\ra T$, where $F^\times$ is seen as an algebraic group over
$\Q$ (so that, for example, we have $F^\times(\Ade)=(\Ade\otimes_\Q F)^\times=
\Ade_F^\times$).
We already have a morphism $h:\Se\ra T_\R$, which gives a morphism
of complex algebraic groups $h_\C:\Se_\C\ra T_\C$. 
Remember that $\Se_\C\simeq\GL_{1,\C}\times\GL_{1,\C}$, and that
we denoted by $r:\GL_{1,\C}\ra\Se_\C$ the embedding of the first factor
(see page~\pageref{def_r}). We get a morphism $h_\C\circ r:\GL_{1,\C}\ra
T_\C$. As $T$ is an algebraic group over $\Q$, this morphism is defined
over a finite extension of $\Q$ in $\C$, and we call this extension
$F$. We get a morphism of $F$-algebraic groups $\GL_{1,F}\ra T_F$, hence
a morphism of $\Q$-algebraic groups $F^\times\ra\Res_{F/\Q}T_F$, where
$\Res_{F/Q}T_F$ is the algebraic group that sends a $\Q$-algebra $R$ to
$T(R\otimes_\Q F)$. Composing this with the norm $N_{F/\Q}:\Res_{F/\Q}T_F\ra
T$, we finally get a morphism $r(h):F^\times\ra T$, called the
\emph{reciprocity morphism for $(T,h)$}.
We take
$\lambda_M$ to be induced by $r(h)$.
\footnote{This is the opposite of Deligne's convention in~\cite{Deligne-SV} (3.9.1);
I have chosen to follow Pink's convention in~\cite{Pink-these}
Chapter~11.}

So if $(G,h)$ is a Shimura datum with $G$ a torus, we get a canonically
defined model of the associated Shimura variety over the field of
definition of $h_\C\circ r$.

\subsubsection{The reflex field}
\label{sssection_rf}

Let $(G,h)$ be a Shimura datum satisfying the hypothesis of
Theorem~\ref{thm_CC_SV}, and let $\nu:G\ra T:=G/G_\der$ be the quotient
morphism. For every open compact subgroup $K$ of $G(\Af)$,  
we have a bijection, induced by $\nu$:
\[\pi_0(M_K(G,h)(\C))\iso T(\Q)\sous T(\Ade)/\nu(K_\infty\times K).\]

Suppose that we have a model $(M_K)_K$ of the Shimura
variety of $(G,h)$ over a number field $F\subset\C$. We would expect
the Shimura variety of $(T,\nu\circ h)$ to also have a model over
$F$, and the isomorphism above to be $\Gal(\overline{F}/F)$-equivariant,
where the action on the right hand side is given by the 
morphism $r(\nu\circ h):F^\times\ra T$
constructed in~\ref{sssection_SV_tori}. 

Remember that $r:\GL_{1,\C}\ra\Se_\C\simeq\GL_{1,\C}\times\GL_{1,\C}$ 
is the embedding of the first factor
(see page~\pageref{def_r}). By the previous paragraph, we would
expect $F$ to contain the field of definition of $\nu_\C\circ h_\C\circ r$.
In fact it would make sense to take $F$ to be the field of definition
of $h_\C\circ r:\GL_{1,\C}\ra G_\C$, except that $h$ is only significant up
to conjugation. This motivates the following definition.

\begin{definition}\label{def_reflex_field}
Let $(G,h)$ be a Shimura datum. The \emph{reflex field} $F(G,h)$
of $(G,h)$ is the
field of definition of the conjugacy class of $h_\C\circ r:\GL_{1,\C}\ra G_\C$.

\end{definition}

Let $F=F(G,h)$. Then $F$ is a finite extension of
$\Q$ in $\C$, and, for every morphism $\rho$ of $G$ into a commutative algebraic
group, the morphism $\rho_\C\circ h_\C\circ r$ is defined over $F$.
Note that $h_\C\circ r$ itself is not necessarily defined over $F$.

\begin{example}\label{ex_GU}
Let $E=\Q[\sqrt{-d}]$ 
be an imaginary quadratic extension of $\Q$, let $p\geq q\geq 1$
be integers, and set $n=p+q$. Let $J=\begin{pmatrix}
I_p & 0 \\ 0 & -I_q \end{pmatrix}
\in\GL_n(\Z)$.
For every commutative ring $R$, we denote by
$x\mapsto\overline{x}$ the involution of $R\otimes_\Z\Of_E$ induced by
the nontrivial element of $\Gal(E/\Q)$, and, for every $Y\in M_n(R\otimes_\Q 
E)$, we write $Y^*={}^t\overline{Y}$.

The general unitary group $\GU(p,q)$ is the $\Z$-group scheme defined by
\[\GU(p,q)(R)=\{g\in\GL_n(R\otimes_\Z \Of_E)\mid\exists c(g)\in R^\times,
g^* J g=c(g) J\}\]
for every commutative ring $R$. Then $\GU(p,q)_\Q$ is a connected reductive
algebraic group, and we have a morphism of group schemes $c:\GU(p,q)\ra
\GL_1$, whose kernel is the unitary group $\U(p,q)$.

Let $h:\Se\ra\GU(p,q)_\R$ be the morphism defined by
\[h(z)=\begin{pmatrix} z I_p & 0 \\ 0 & \overline{z} I_q\end{pmatrix}
\in\GU(p,q)(\R).\]
Then $(G,h)$ is a Shimura datum, and $K_\infty$ is the set of matrices
$\begin{pmatrix} g_1 & 0 \\ 0 & g_2\end{pmatrix}$ such that
$g_1\in\GL_p(\C)$, $g_2\in\GL_q(\C)$ and there exists $c\in\R$ with
$g_1^*g_1=cI_p$ and $g_2^*g_2=cI_q$.

We use the isomorphism $\C\otimes_\Q E\iso\C\times\C$ sending
$x\otimes 1+y\otimes\sqrt{-d}$ to $(x+\sqrt{-d}y,x-\sqrt{-d}y)$ to
identify $\GU(p,q)(\C)$ to a subgroup of $\GL_n(\C)\times\GL_n(\C)$;
note that the involution $g\mapsto\overline{g}$ of $\GU(p,q)(\C)$ corresponds
to switching the two factors.
With this convention,
we have $h_\C\circ r(z)=\left(\begin{pmatrix} z I_p & 0 \\ 0 & I_q\end{pmatrix},
\begin{pmatrix} I_p & 0 \\ 0 & z I_q\end{pmatrix}\right)$
for every $z\in\C^\times$. 

It is easy to see that $h_\C\circ r$ is defined over $E$ but not over $\Q$.
On the other hand, the reflex field of $(G,h)$ is $E$ if $p>q$ and $\Q$
if $p=q$.

\end{example}

\begin{example}\label{ex_rf_GSp}
Let $(G,h)$ be the Shimura datum of Example~\ref{ex_SD_GSp}
(so that $G=\GSp_{2d}$). For every $z\in\C^\times$, we have
\[h_\C\circ r(z)=\begin{pmatrix}
\frac{1}{2}(z+1) I_d & -\frac{1}{2i}(z-1)I_d \\
\frac{1}{2i}(z-1)I_d & \frac{1}{2}(z+1) I_d
\end{pmatrix}=P\begin{pmatrix} z I_d & 0 \\ 0 & I_d\end{pmatrix}P^{-1},\]
where
\[P=\begin{pmatrix}\frac{1}{\sqrt{2}} I_d & \frac{i}{\sqrt{2}} I_d \\
\frac{i}{\sqrt{2}} I_d & \frac{1}{\sqrt{2}} I_d \end{pmatrix}.\]
So the reflex field of $(G,h)$ is $\Q$.

\end{example}

\subsubsection{Canonical models}

We are now ready to define canonical models.

\begin{definition}\label{def_canonical_model}
Let $(G,h)$ be a Shimura datum and let $F=F(G,h)$. A \emph{canonical
model} of $M(G,h)(\C)$ is a model $(M_K)_K$ over $F$ such that, for every
torus $u:H\subset G$ and every $h':\Se\ra H_\R$ such $u\circ h'$ and
$h$ are $G(\R)$-conjugated (i.e. such that $u$ induces a morphism of
Shimura data from $(H,h')$ to $(G,h)$), the morphism
\[u:M(H,h')(\C)\ra M(G,h)(\C)\]
is defined over the compositum $F\cdot F(H,h')\subset\C$, where we use
as model of $M(H,h')(\C)$ over $F(H,h')$ the one defined
in~\ref{sssection_SV_tori}.

\end{definition}

\begin{examples}\label{ex_can_model}
\begin{itemize}
\item[(1)]
If $G$ is a torus, then the model of~\ref{sssection_SV_tori} is
a canonical model of $M(G,h)(\C)$.

\item[(2)] If $(G,h)$ is the Shimura datum of Example~\ref{ex_SD_GSp}
(so that $G=\GSp_{2d}$), then the schemes $(\Mcal_{d,K,\Q})_{K\subset\GSp_{2d}(\Af)}$
of~\ref{sssection_SMV} form a canonical model of
$M(G,h)(\C)$. This is not obvious but follows from the main theorem
of complex multiplication; see Section~4 of \cite{Deligne-SV}.

\end{itemize}
\end{examples}

At the time of Deligne's paper \cite{Deligne-SV}, it was not
known whether all Shimura varieties have canonical models (spoiler: this is
now known to be true, see Theorem~\ref{thm_ECM}), 
but it was possible to prove their uniqueness.
If $u:(H,h')\ra(G,h)$ is a morphism of Shimura data as in
Definition~\ref{def_canonical_model} (so that $H$ is a subtorus of $G$),
the points of the image of $M(H,h')(\C)$ in $M(G,h)(\C)$ are called
\emph{special points}. The fact that canonical models are unique if they
exist relies on the following two points:
\begin{itemize}
\item[(i)] The set of all $G(\Af)$-translates of
special points are dense in $M(G,h)(\C)$ (see Proposition~5.2
of~\cite{Deligne-SV}; in fact, we just need to take $G(\Af)$-translates of
one special point);
\item[(ii)] For every finite extension $F'\subset\C$ of $F(G,h)$, there
exists $u:(H,h')\ra(G,h)$ as above such that $F(H,h')$ and $F'$ are
linearly disjoint over $F(G,h)$ (Th\'eor\`eme~5.1 of~\cite{Deligne-SV}).

\end{itemize}

From this, we can deduce:

\begin{theorem}[Corollaire~5.4 of~\cite{Deligne-SV}]\label{thm_can_model_unique}
Let $u:(G_1,h_1)\ra(G_2,h_2)$ be a morphism of Shimura data, and suppose that the
Shimura varieties of $(G_1,h_1)$ and $(G_2,h_2)$ have canonical models.
Then the morphism $M(G_1,h_1)(\C)\ra M(G_2,h_2)(\C)$ of Shimura varieties
corresponding to $u$
is defined over any common extension $F$ of $F(G_1,h_1)$ and $F(G_2,h_2)$ in
$\C$.

\end{theorem}

\begin{cor}
Let $(G,h)$ be a Shimura datum. Then a canonical model of $M(G,h)(\C)$ is
unique up to unique isomorphism if it exists.

\end{cor}

\begin{cor}
Let $(G,h)$ be a Shimura datum satisfying the hypothesis of
Theorem~\ref{thm_CC_SV}, let $\nu:G\ra T:=G/G_\der$ be the quotient morphism,
and let $F=F(G,h)$. Then the action of
$\Gal(\overline{F}/F)$ on $\pi_0(M(G,h)(\C))\simeq\pi_0(T(\Q)\sous T(\Ade))
/\pi_0(K_\infty)$ is given by the reciprocity morphism
for $(T,\nu\circ h)$.

\end{cor}

\begin{proof}
Apply Theorem~\ref{thm_can_model_unique} to $\nu:(G,h)\ra(T,\nu\circ h)$.

\end{proof}

Using the same techniques as for Theorem~\ref{thm_can_model_unique}, we also
get the following very useful result.

\begin{prop}[Corollaire~5.7 of~\cite{Deligne-SV}]\label{prop_construction_cm}
Let $u:(G_1,h_1)\ra (G_2,h_2)$ be a morphism of Shimura
data such that the underlying morphism of algebraic groups is a closed
immersion and that $F(G_1,h_1)\subset F(G_2,h_2)$. If
$M(G_2,h_2)(\C)$ has a canonical model, then so does $M(G_1,h_1)(\C)$.

\end{prop}

\section{Lecture 2: arithmetic Shimura varieties}

In lecture 1, we have defined Shimura varieties over $\C$ and
introduced the notion of canonical model of a Shimura variety. We also
discussed in some detail the example of the Siegel modular varieties, 
which have canonical models coming from their modular interpretation.
In this lecture, we first want to present different types of Shimura
varieties, with each type being contained in the next one:
\begin{itemize}
\item \textbf{The Siegel modular variety}: it is a moduli space of
principally polarized abelian schemes with some level structure;

\item \textbf{PEL type Shimura varieties}: they have an interpretation
as moduli spaces of polarized abelian schemes with multiplication by
the ring of integers $\Of$ of some number field and some level structure
(here ``P'' means ``polarization'', ``E'' means ``endomorphisms'' in
reference to the action of $\Of$ and ``L'' means ``level structure'');

\item\textbf{Hodge type Shimura varieties}: they come from Shimura
data $(G,h)$ that have an injective morphism into a Siegel Shimura datum
(Example~\ref{ex_SD_GSp});

\item\textbf{Abelian type Shimura varieties}: their Shimura datum is
``isogenous'' to a Hodge type Shimura datum (in a way to be made precise later);

\item\textbf{General Shimura varieties}: all Shimura varieties.

\end{itemize}

The further we go down in the list, the less is known about the geometry
of the Shimura variety (and the higher the price for what we know), because
many of the techniques we have rely on the interpretation of the Shimura
varieties as moduli spaces of abelian schemes, and this is only
really available for PEL type Shimura varieties (there is a modular
interpretation for
Hodge type Shimura varieties, but it is harder to use).

For example, PEL type Shimura varieties naturally come with an \emph{integral
model} defined over a localization of the ring of integers of their
reflex field, but it took a lot of effort
to construct integral models for Hodge type and abelian type Shimura
varieties, and to formulate their properties
(see~\ref{sssection_can_integral_models} and Theorems~\ref{Kisin-IM1}
and~\ref{Kisin-IM2}); as far as we know, nothing
is known for general Shimura varieties.

We are not claiming that the classification above is the only measure of the
complexity of a Shimura variety, whatever that means. For example, as we
will discuss in lecture 3, if one wants to study the cohomology of Shimura
varieties and their zeta functions, then the simplest case is not the case of
Siegel modular varieties, but rather of compact PEL type Shimura varieties
whose group has no endoscopy and a simply connected derived subgroup;
we will introduce some examples of these, known as
\emph{Kottwitz's simple Shimura varieties}, in~\ref{ssection_SSV}.

\subsection{PEL type Shimura varieties}

These Shimura varieties were introduced by Kottwitz in \cite{K-NP}, but
we will follow the presentation of Lan for the moduli problems
(cf. \cite{Lan-PEL} 1.4.1), which is
closer to our definition of the Siegel moduli problem. The
equivalence between the two definitions is proved in \cite{Lan-PEL} 1.4.3.

\vspace{.3cm}

\subsubsection{PEL data}

\begin{definition}[see \cite{Lan} 5.1 or Definition~1.2.13 of~\cite{Lan-PEL}]
\label{def_PELD}
An \emph{(integral) PEL datum} is a quintuple $(\Of,*,\Lambda,\lan\cdot,\cdot\ran,h)$, where:
\begin{itemize}
\item[(1)] $\Of$ is an order in a finite-dimensional semisimple $\Q$-algebra
$B$ (that is, $\Of$ is a subring of $B$ that is a free $\Z$-module and
spans the $\Q$-vector space $B$);

\item[(2)] $*$ is a \emph{positive involution} of $\Of$, i.e. an
anti-automorphism of rings of order $2$ such that, for every
$x\in\Of\setminus\{0\}$, we have $\Tr_{(B\otimes_\Q\R)/\R}(xx^*)>0$;

\item[(3)] $\Lambda$ is an $\Of$-module that is finitely generated and free
as a $\Z$-module;

\item[(4)] $\lan\cdot,\cdot\ran:\Lambda\times \Lambda\ra\Z$ is a $\Z$-bilinear
alternating map such that, for all $x,y\in \Lambda$ and $b\in\Of$, we have
\[\lan bx,y\ran=\lan x,b^* y\ran;\]

\item[(5)] $h:\C\ra\End_{B\otimes_\Q\R}(\Lambda\otimes_\Z\R)$ is an $\R$-algebra
morphism such that:
\begin{itemize}
\item[(a)] For $z\in\C$ and $x,y\in \Lambda\otimes_\Z\R$, we have
\[\lan h(z)(x),y\ran=\lan x,h(\overline{z})(y)\ran;\]

\item[(b)] The $\R$-bilinear pairing $\lan \cdot,h(i)(\cdot)\ran$,
which is symmetric by (a), is also positive definite.

\end{itemize}
\end{itemize}

\end{definition}

Let $(\Of,*,\Lambda,\lan\cdot,\cdot\ran,h)$ be a PEL datum. We define a group
scheme $G$ over $\Z$ by
\[G(R)=\{g\in\End_{\Of\otimes_\Z R}(\Lambda\otimes_\Z R)
\mid\exists c(g)\in R^\times,\ 
\lan g(\cdot),g(\cdot)\ran=c\lan\cdot,\cdot\ran\}\]
for every commutative ring $R$.
We also get a morphism of group schemes $c:G\ra\GL_1$.
The morphism of $\R$-algebras $h:\C\ra\End_{B\otimes_\Q\R}(\Lambda\otimes_\Z\R)$
induces a morphism of $\R$-algebraic groups $h:\Se\ra G_\R$.

\begin{prop}\label{prop_PELD_to_PELSD}
The couple $(G_\Q^0,h)$ satisfies conditions (a)-(c) in the definition of
a Shimura datum (Definition~\ref{def_SD}).

\end{prop}

Let us explain why this proposition is true. Condition (a) of
Definition~\ref{def_SD} follows
from the fact that $h:\C^\times\ra\End_{B\otimes_\Q\R}(\Lambda\otimes_\Z\R)$ is
a morphism of $\R$-algebras, and condition (c) of loc. cit. follows
from condition (5)(b) of Definition~\ref{def_PELD}. 
We prove condition (b), on the decomposition of $\Lie(G_\C)$ into
eigenspaces for the action of $\Ad\circ h_\C$.
Let $V=\Lambda\otimes_\Z\R$,
a finite-dimensional $\R$-vector space. Then $h:\C\ra\End_\R(V)$ is a morphism
of $\R$-algebras, hence defines
a structure of $\C$-vector space on $V$. In particular, we have
$V\otimes_\R\C=V_1\oplus V_2$, where $\C$ acts on $V_1$ (resp.
$V_2$) via multiplication by $z\mapsto z$ (resp. $z\mapsto\overline{z}$). Condition (b)
follows from this, because $\End(V)\simeq V^*\otimes_\R V$ and
$\ggoth:=\Lie(G_\R)\subset\End_\R(V)$.

We keep the same notation. Then $A:=V/\Lambda$ is a
complex torus. We would like this torus to be an abelian variety, so
we need a polarization on it, that is, a positive definite
Hermitian form $H$ on $V$ such that $\im(H)$ takes integer values
on $\Lambda$. But we already have an alternating form on $\Lambda$, so
we know from~\ref{sssection_cav} how to proceed:
define $H$ by
\[H(v,w)=\lan h(i)(v),w\ran+i\lan v,w\ran.\]
The fact that $\im(H)$ takes integral values on $\Lambda$ is clear,
the pairing $H$ is Hermitian by condition (5)(a) and positive definite
by condition (5)(b). So the torus $A$ is an abelian variety, with
dual abelian variety $A^\vee=V/\Lambda^\vee$, where
$\Lambda^\vee=\{v\in V\mid\forall w\in\Lambda, \lan v,w\ran\in\Z\}$
is the dual lattice of $\Lambda$. The polarization $\lambda:A\ra A^\vee$
defined by $H$ is then just the map induced by $\Lambda\subset\Lambda^\vee$.
We have a morphism of algebraic groups $\Se\flnom{h} G_\R\subset\GL(V)$, and the
decomposition $V\otimes_\R\C=\bigoplus_{p,q\in\Z} V^{p,q}$ induced by this morphism
(see~\ref{sssection_Se}) is the same as the Hodge structure
coming from the isomorphism $V\simeq\H_1(A,\R)$. 
Note also that the action of $\Of$ on $L$ defines a morphism of rings
$\iota:\Of\ra\End_\C(A)$ satisfying the \emph{Rosati condition}, which
means that, for every $b\in\Of$, we have
\[\lambda\circ i(b^*)=i(b)^\vee\circ\lambda.\]

\begin{remark}[\cite{K-NP}, sections 5 and 7]
Let $(\Of,*,\Lambda,\lan\cdot,\cdot\ran,h)$ be a PEL datum, and let
$G$ be the associated group scheme. If the $\Q$-algebra $B=\Of\otimes_\Z\Q$
is simple, then all the simple factors of $G_\der(\C)$ are of the same
type, which is $A$, $C$ or $D$. If we are in type $A$ or $C$, then
$G_\Q$ is connected and reductive. In type $D$, the group
$G_\Q$ is reductive and
has $2^{[F_0:\Q]}$ connected components, where $F_0$ is the field of fixed
points of $*$ in the center of $\Of\otimes_\Z\Q$, so it is never connected.

\end{remark}

We now give some examples covering all the types.

\begin{example}\label{ex_PEL_C}
Take $\Of=\Z$, $*=\id_\Of$, 
$\Lambda=\Z^{2d}$,
$\lan\cdot,\cdot\ran$ the perfect symplectic pairing 
with matrix $\begin{pmatrix} 0 & I_d \\ -I_d & 0\end{pmatrix}$ in the
canonical basis of $\Z^{2d}$, and
$h:\C\ra M_{2d}(\R)$ defined by
\[h(a+ib)=\begin{pmatrix}aI_d & -b I_d \\ b I_d & a I_d\end{pmatrix}.\]
Then we get a PEL datum, and the couple $(G_\Q,h)$ is the Siegel Shimura
datum of Example~\ref{ex_SD_GSp}.

\end{example}

\begin{example}\label{ex_PEL_A}
Let $B=E\subset\C$ be an imaginary quadratic extension of $\Q$,
$\Of$ be an order in $E$ (for example the ring of integers $\Of_E$), 
$*$ be the restriction to $\Of$ of complex conjugation, 
\footnote{If $(1,a)$ is a $\Z$-basis of $\Of_E$, then every
order $\Of$ of $E$ is contained in $\Of_E$ and
of the form $\Z\oplus fa\Z$, where
$f=[\Of_E:\Of]$ is the \emph{conductor} of $\Of$. In particular,
the order $\Of$ is stable by the nontrivial element of $\Gal(E/\Q)$,
which is the restriction of complex conjugation to $E$.}
$\Lambda=
\Of^{p+q}$ with $p\geq q\geq 0$. Choose $\varepsilon\in\Of$ such
that $-i\varepsilon\in\R_{>0}$, let $H$ be the Hermitian pairing on
$\Lambda$ with matrix $\begin{pmatrix} I_p & 0 \\ 0 & -I_q\end{pmatrix}$
and $\lan\cdot,\cdot\ran$ be the
alternating pairing $\Tr_{\Of/\Z}(\varepsilon H)$ on $\Lambda$. Finally, define
$h:\C\ra\End_\C(\Lambda\otimes_\Z\R)=M_{p+q}(\C)$ by
$h(z)=\begin{pmatrix} z I_p & 0 \\ 0 & \overline{z} I_q\end{pmatrix}$.
Then we get a PEL datum, and the couple $(G_\Q,h)$ is the Shimura
datum of Example~\ref{ex_GU}.

If $q=0$, then $(G,h)$ does not satisfy condition (d) of
Definition~\ref{def_SD}.

\end{example}

\begin{example}\label{ex_PEL_Kottwitz}
At the other extreme of type A, we have the PEL data considered by Kottwitz in~\cite{K-SSV}
(see the beginning of Section~1 of that paper). Let $F_0$ be a totally real number field,
$F$ be a totally imaginary quadratic extension of $F_0$ and $B$ be a division algebra 
with center $F$; we suppose given a positive involution $*$ on $B$ whose restriction
to $F$ is the nontrivial element of $\Gal(F/F_0)$.

We take $\Lambda\otimes_\Z\Q=B$ with its action by left translations.
The alternating form $\lan .,.\ran$ on $\Lambda\otimes_\Z\Q$
is given by $\lan x,y\ran=\Tr_{B/\Q}(xb y^*)$, where $b$ is an element of
$B^\times$ such that $b^*=-b$.
We have $\End_B(\Lambda\otimes_\Z\Q)=B^\op$, so $G_\Q=\{x \in (B^\op)^\times|xx^\ddag\in\Q^\times\}$,
where $(.)^\ddag$ is the involution of $B^\op$ define by $x^\ddag=bx^*b^{-1}$ (which is not
positive in general).

\textbf{Warning}: What Kottwitz calls $D$ and $*$ in~\cite{K-SSV} is what we call $B^\op$ and
$\ddag$ here, so in particular the involution $*$ of~\cite{K-SSV} is not assumed to be positive
(and the PEL datum is not defined there).
We changed the notation to be consistent with the rest of the text and with
Kottwitz's other paper~\cite{K-NP}.

We have $F_0\otimes_\Q\R\simeq\R^r$, $F\otimes_\Q\R\simeq\C^r$ and $B\otimes_\Q\R\simeq M_n(\C)^r$,
where $n=\sqrt{\dim_\Q B}$; the involution $*$ is conjugate on each
factor to the positive involution $A\mapsto{}^t\overline{A}$. So, if
$G_0=\{x\in (B^\op)^\times|xx^\ddag=1\}$, then $G_{0,\R}$ is isomorphic to a product
$\prod_{i=1}^r\U(p_i,q_i)$ of unitary groups, and $G_\R$ is isomorphic to
the corresponding
general group (i.e. the subgroup of $\prod_{i=1}\GU(p_i,q_i)$ where all the
multipliers are equal). We define $h:\C^\times\ra G_\R$ by sending
$z\in\C^\times$ to element of $G(\R)$ corresponding to the family of matrices
$\left(\begin{pmatrix} z I_{p_i} & 0 \\ 0 & \overline{z} I_{q_i}\end{pmatrix}\right)_{1\leq i\leq r}$.
This extend to a morphism of $\R$-algebras $h:\C\ra B^\op\otimes_\Q\R$ sending $0$ to $0$;
as $\C$ is commutative, we can also see $h$ as a morphism of $\R$-algebras
$\C\ra B\otimes_\Q\R$.

We have defined what is known as a rational PEL datum. To get an integral PEL
datum, we must suppose that there exists an order $\Of$ in $B$ such that
$\Of^*=\Of$ and $\lan.,.\ran$ takes integer values on $\Of$.

\end{example}

\begin{example}\label{ex_PEL_D}
Let $B$ be a quaternion algebra over $\Q$ such that $B\otimes_\Q\R\simeq\Ham$,
and let $\Of$ be an order in $B$ that is stable by the involution
of $\Ham$ defined by $(x+iy+jz+kt)^*=x-iy-jz-kt$. Let $\Lambda=\Of^{2n}$,
let $\lan\cdot,\cdot\ran$ 
be $\Tr_{\Of/\Z}\circ H$, where $H$ is the skew-Hermitian
pairing on $\Lambda$ with matrix $\begin{pmatrix} 0 & I_n \\ -I_n & 0
\end{pmatrix}$. Define $h:\C\ra\End_{B\otimes_\Q\R}(\Lambda\otimes_\Z\R)=
M_n(\Ham)$ by $h(a+ib)=\begin{pmatrix} a I_n & -b I_n \\ bI_n & aI_n
\end{pmatrix}$. 

The group $(G^0)_{\R,\der}$ is often denoted by $\SO^*_{2n}$; it is
a quasi-split outer form of the split orthogonal
group $\SO_{2n}$, hence is of type $D_n$. But note that the algebraic
group $G_\Q$ is not connected.

\end{example}

\begin{definition}
We say that a Shimura datum is of \emph{PEL type} if it is of the form
$(G^0_\Q,h)$, where $(G,h)$ comes from a PEL datum.
The corresponding Shimura varieties are called \emph{PEL type Shimura
varieties}.

\end{definition}

\subsubsection{PEL moduli problems}

Just as in the case of the Siegel modular variety, PEL type Shimura
varieties are the solution of a moduli problem,
\footnote{Well, almost. See Proposition~\ref{prop_PEL_can_model}
below for a more precise statement.}
known as a \emph{PEL moduli problem}. We now discuss these.

Let $(\Of,*,\Lambda,\lan\cdot,\cdot\ran,h)$ be a PEL datum, and let
$(G,h)$ be defined as before.
If $V=\Lambda\otimes_\Z\R$, then we saw that the morphism $h:\C\ra
\End_{\Of\otimes_\Z\R}(V)$ defines a decomposition
$V\otimes_\R\C=V^{-1,0}\oplus V^{0,-1}$.
The \emph{reflex field} $F$ of the PEL datum is the field of definition of
the isomorphism class of $V^{-1,0}$ as an $\Of\otimes_\Z\C$-module, that is,
the subfield of $\C$ generated by the elements
$\Tr(b,V^{1,0})$, for $b\in\Of$. It
is also equal to the field of the definition of the conjugacy class of
$h_\C\circ r:\GL_{1,\C}\ra G_\C$, 
that is, the reflex field of $(G^0_\Q,h)$ when this couple is a Shimura
datum.

\begin{definition}\label{def_p_good}
We say that a prime number $p$ is \emph{good} for the PEL datum if:
\begin{itemize}
\item $p$ is unramified in $\Of$ (i.e. it does not divide the discriminant
of $\Of/\Z$);
\item $p$ does not divide $[\Lambda^\vee:\Lambda]$, where
$\Lambda^\vee$ is as before the lattice 
$\{v\in V\mid\forall w\in\Lambda, \lan v,w\ran\in\Z\}$;
\item $p\not=2$ if the PEL datum has a factor of type $D$ (i.e. if
$\Of\otimes_\Z\R$ has a simple factor isomorphic to an algebra
$M_n(\Ham)$ with its canonical positive involution).

\end{itemize}

\end{definition}

If $p$ is not good we say that it is bad.
Note that, if $p$ is good, then $F$ and $G_\Q$ are unramified at $p$;
in fact, the group $G(\Z_p)$ is then a hyperspecial maximal compact
subgroup of $G(\Q_p)$.

Let $T$ be a set of good primes (finite or infinite), and let
$\Of_{F,T}$ be the localization $\Of_F[\frac{1}{p},\ p\not\in T]$.

\begin{definition}\label{def_PEL_moduli}
Let $n$ be a positive integer that is prime to all the elements of $T$. Then
the PEL moduli problem at level $n$ defined by the fixed PEL datum is
the contravariant
functor $\Mcal_n$ from the category of $\Of_{F,T}$-schemes to the
category of sets sending an $\Of_{F,T}$-scheme $S$ to the set of
isomorphism classes of quadruples $(A,\lambda,\iota,(\eta,\varphi))$, where:
\begin{itemize}
\item $A$ is an abelian scheme over $S$;
\item $\lambda:A\ra A^\vee$ is a polarization whose degree is prime to
all the elements of $T$;
\item $\iota:\Of\ra\End_S(A)$ is a morphism of rings satisfying
the Rosati condition: for every $b\in\Of$, we have
$\lambda\circ\iota(b^*)=\iota(b)^\vee\circ\lambda$;
\item $(\eta,\varphi)$ is a level $n$ structure on $A$, i.e.,
$\eta:A[n]\iso(\Lambda/n\Lambda)_S$ is an $\Of$-equivariant isomorphism
of group schemes and $\varphi:\underline{\Z/n\Z}_S\iso\mu_{n,S}$ is an
isomorphism of group schemes such that $\varphi\circ\lan\cdot,\cdot\ran
\circ\eta$ is the Weil pairing defined by $\lambda$ on $A[n]_S$, and
moreover $(\eta,\varphi)$ are liftable to level $m$ structures for every
prime-to-$T$ multiple $m$ of $n$, in the sense of \cite{Lan-PEL}
Definition~1.3.6.2.

\end{itemize}
We furthermore require that this quadruple satisfy the following
\emph{determinant condition} (see Definition~1.3.4.1 of~\cite{Lan-PEL}): 
let $\alpha_1,\ldots,\alpha_t$ be a basis
of the $\Z$-module $\Of$, and let $X_1,\ldots,X_t$ be indeterminates.
Then $\det(\iota(\alpha_1) X_1+\ldots+\iota(\alpha_t) X_t,\Lie(A))$ 
is a polynomial in
$\Of_S[X_1,\ldots,X_t]$, and the condition says that this polynomial
is equal to the image by the map $\Of_{F,T}\ra\Of_S$ of the
polynomial $\det(\alpha_1X_1+\ldots+\alpha_tX_t,V^{1,0})\in
\Of_{F,T}[X_1,\ldots,X_t]$.
\footnote{It is not totally obvious that the second polynomial, which is
a priori in $F[X_1,\ldots,X_t]$, has its coefficients in $\Of_{F,T}$.
But it is also not too hard to check. See for example pages 389-390 of
\cite{K-NP}.}

\end{definition}

\begin{remark}
\begin{itemize}
\item[(1)] We did not specify the relative dimension of the abelian
scheme $A$ in the moduli problem, because it is necessarily equal to
$\dim_\C(V^{1,0})$
by the determinant condition.

\item[(2)] For $n\in\{1,2\}$, 
the objects of the moduli problem $\Mcal_n$ can have
nontrivial automorphisms, so it would make more sense to see
$\Mcal_n$ as a functor with values in groupoids, i.e. as a stack. See
Remark~\ref{rmk_small_level}.

\item[(3)] As in Remark~\ref{rmk_MdK}, it is also possible to define the moduli
problem $\Mcal_K$
for more general levels $K$, i.e. for open compact subgroups
$K$ of $G(\Af)$, though there is a condition on $K$ corresponding to
the condition that $n$ be prime to $T$. Let $\Af^T$ be the ring of
prime-to-$T$ adeles of $\Q$, i.e. the restricted product of the
$\Q_p$ for $p\not\in T$, and $\Ade_T$ be the restricted product of the
$\Q_p$ for $p\in T$. We have $\Af=\Ade_T\times\Af^T$, and the condition on
$K$ is that $K=K_TK^T$, where $K^T\subset G(\Af^T)$ and $K_T=\prod_{p\in T}K_p$
with $K_p=G(\Z_p)$ for every $p\in T$.

\end{itemize}
\end{remark}

The following theorem is a consequence of Mumford's representability
result (Theorem~\ref{thm_Mdn}),
as explained in Section~5 of Kottwitz's paper~\cite{K-NP}. See
also Corollaries~1.4.1.12 and~7.2.3.10 of Lan's book~\cite{Lan-PEL}
for more details.

\begin{theorem}
If $n\geq 3$, then the functor $\Mcal_n$ is representable by a smooth
quasi-projective scheme over $\Of_{F,T}$.

\end{theorem}

In fact, by Theorem~1.4.1.11 and Corollary~7.2.3.10 of \cite{Lan-PEL}, the
functor $\Mcal_K$ is representable by
a smooth quasi-projective scheme for $K$ small enough.

\begin{remark}
We also have an action of the Hecke operators defined by elements
$g\in G(\Af^T)$ (i.e. the Hecke operators that are trivial at primes of $T$)
on the tower $(\Mcal_K)$; as for the Siegel moduli problem, the element
$g$ acts on the level structure.
See Remark~1.4.3.11 of \cite{Lan-PEL}
(and the comparison result of Proposition~1.4.3.4 of loc. cit.).

\end{remark}

\subsubsection{PEL moduli problems and canonical models}

Consider the PEL datum of Example~\ref{ex_PEL_C}. The associated
couple $(G_\Q,h)$ is the Siegel Shimura datum of Example~\ref{ex_SD_GSp}, and
we have seen in Example~\ref{ex_can_model}(2) that the corresponding
moduli problem, with $T=\varnothing$, defines a canonical model of the
Shimura variety of $(G_\Q,h)$. 

We would like something like this to be true for general PEL data,
but one obstruction is that the moduli problem only depends on the
completions of the symplectic $\Of$-module $(\Lambda,\lan .,,\ran)$ at the places of $\Q$;
this is not obvious in Definition~\ref{def_PEL_moduli}, but
it becomes so
if we use the moduli problem of Definition~1.4.2.1 of~\cite{Lan-PEL}
or of Section~5 of~\cite{K-NP} (the equivalence of the two moduli
problems is proved in \cite{Lan-PEL} 1.4.3). 
Suppose that we have a finitely generated $\Of$-module $\Lambda'$
that is free as a $\Z$-module,
with a symplectic form $\lan.,.\ran'$ satisfying condition (4) of
Definition~\ref{def_PELD} and
such that, for every place
$v$ of $\Q$, $\Lambda\otimes_\Z\Q_v$ and $\Lambda'\otimes_\Z\Q_v$ are
isomorphic as symplectic $B\otimes_\Q{\Q_v}$-modules. Then we can see $h$ as a morphism
$\C\ra\End_{B\otimes_\Q\R}(\Lambda'\otimes_\Z\R)$, and it will satisfy condition
(5) of Definition~\ref{def_PELD}, so we get a PEL datum $(\Of,*,\Lambda',\lan.,.\ran',h)$
that defines the same moduli problem as the original PEL datum. In particular,
we will get a Shimura datum $({G'}^0_\Q,h)$, and we should also see the canonical
model for the Shimura variety of $({G'}^0_\Q,h)$ in the PEL moduli problem. A priori,
there is no reason for $G$ and $G'$ to be isomorphic, and they won't be
in general; we only know that
$G_{\Q_v}\simeq G'_{\Q_v}$ for every place $v$ of $\Q$.

In fact, we have the following result.

\begin{prop}\label{prop_PEL_can_model}
Suppose that the semisimple $\Q$-algebra $\Of\otimes_\Z\Q$ has no
simple factor of type $D$ (see Definition~\ref{def_p_good}). Then:
\begin{itemize}
\item[(i)] (Remark~1.4.4.4 of~\cite{Lan-PEL}.)
Let $T\subset T'$ be two sets of good prime numbers,
let ${K}^{T'}$ be an open compact subgroup of $G(\Af^{T'})$, and let
$K=K^{T'}\prod_{p\in T'}G(\Z_p)$. Let $\Mcal_K$ (resp. $\Mcal'_K$) be the
moduli problem over $\Spec(\Of_{F,T})$ (resp. $\Spec(\Of_{F,T'})$)
from Definition~\ref{def_PEL_moduli}, where we take the set of good
primes to be $T$ (resp. $T'$). Then the forgetful functor
$\Mcal_K\ra\Mcal'_K\times_{\Spec(\Of_{F,T'})}\Spec(\Of_{F,T})$ is an
isomorphism.

\item[(ii)] (Sections~7-8 of~\cite{K-NP}.) 
Let $\Ker^1(\Q,G)$ be the kernel of the diagonal map
$\H^1(\Q,G)\ra\prod_{v}\H^1(\Q_v,G)$, where we take the product over all
places $v$ of $\Q$; this is a finite set.
Suppose that $(G_{\Q},h)$ is a Shimura datum, 
\footnote{As we excluded case $D$, the group $G_\Q$ is connected. So, by
Proposition~\ref{prop_PELD_to_PELSD}, the only obstacle to $(G_\Q,h)$ being a Shimura
datum is condition (d) of Definition~\ref{def_SD}.}
and use $T=\varnothing$
to define the moduli problem of Definition~\ref{def_PEL_moduli}, so
that it is a moduli problem over $\Spec(F)$.
\footnote{By (i), this moduli problem is the generic fiber of the
moduli problems defined by nonempty sets $T$.}
Then the projective system $(\Mcal_K)_{K\subset G(\Af)}$ is a disjoint union
indexed by $i\in\Ker^1(\Q,G)$
of projective systems $(\Mcal_K^{(i)})_{K\subset G(\Af)}$,
and each $(\Mcal_K^{(i)})_{K\subset G_i(\Af)}$ is a
canonical model of the Shimura variety of $(G_\Q,h)$.

\end{itemize}
\end{prop}

\begin{remark}
\begin{itemize}
\item[(1)] 
For every $i\in\Ker^1(\Q,G)$, let $G_i$ be the corresponding
inner form of $G_\Q$;\footnote{See~\ref{sssection_inner_U} for the definition
of inner forms.} 
we have $G_{i,v}\simeq G_v$ for every
place $v$ of $\Q$, and we denote by $h_i$
the morphism $\Se\flnom{h} G_\R\simeq G_{i,\R}$.
If $(G_\Q,h)$ is a Shimura datum, then all $(G_i,h_i)$ are, and
the projective system $(\Mcal_K^{(i)})_{K\subset G(\Af)}$ 
in (ii) of the proposition is actually a canonical model of the Shimura
variety of $(G_i,h_i)$ (note that
$G_i(\Af)\simeq G(\Af)$).
But, as noted by Kottwitz at the end of Section~8 of~\cite{K-NP}, under
our hypothesis that $\Of\otimes_\Z\Q$ has no factor of type $D$, all
the groups $G_i$ are isomorphic to $G_\Q$ and the $h_i$ correspond to
$h$. (In fact Kottwitz calculates $\Ker^1(\Q,G)$ in Section~7 of~\cite{K-NP},
and finds that in our situation it is either trivial or isomorphic
to $\Ker^1(\Q,Z(G))$ by the canonical map, where $Z(G)$ is the
center of $G$.)

\item[(2)] As mentioned in (1), there is some information about $\Ker^1(\Q,G)$ in
Section~7 of~\cite{K-NP}, and ways to calculate it.
For example, if $\Of\otimes_\Z\Q$ is simple
of type $C$, we have $\Ker^1(\Q,G)=\{1\}$, which explains why the
moduli problem of Definition~\ref{def_Mdn} gives a canonical model
of the Siegel Shimura variety and not of a finite disjoint union of
copies of it. If $\Of\otimes_\Z\Q$ is simple of type $A$, then
$\Ker^1(\Q,G)$ is automatically trivial in ``half'' of the cases, 
and it is always isomorphic to $\Ker^1(\Q,Z(G))$; for example, it is
trivial for the Shimura datum of Example~\ref{ex_GU}.

\item[(3)] Kottwitz does not say much about the case where
$\Of\otimes_\Z\Q$ has simple factors of type $D$. The situation is
complicated for many reasons: the group $G_\Q$ is not connected,
$\Ker(\Q,G)$ is not trivial, and point (i) of
Proposition~\ref{prop_PEL_can_model} is not true in general, so we must
also be careful about the choice of $T$. In any case, it is still
true that the Shimura varieties $M_K(G^0_\Q,h)(\C)$ are open and closed
subschemes of $\Mcal_{K,\C}$ (see \cite{Lan-Comp} 2.5), and we might even get
canonical models out of this, but we will not pursue this here because there
are other ways to construct canonical models for these Shimura data 
(see~\ref{sssection_Hodge_type}).

\item[(4)] If we want $(G_\Q,h)$ to be a PEL Shimura datum, then this puts
pretty strict conditions on the center of $G_\Q$. For example, if
$F$ is a nontrivial totally real extension of $\Q$, then the group
$\Res_{F/\Q}\GL_2$ (defined by $\Res_{F/\Q}\GL_2(R)=\GL_2(R\otimes_\Q F)$) is
part of a Shimura datum, but this Shimura datum cannot be PEL; more generally,
we have the issue with the group $\Res_{F/\Q}\GSp_{2d}$. This is somewhat
annoying, as sometimes we really do want to consider the Shimura varieties
for these precise groups (see for example Nekovar and Scholl's
\cite{NS}). Fortunately, these Shimura data are of abelian type
(see~\ref{sssection_abelian_type}), so their Shimura
varieties are still understood reasonably well.

\end{itemize}
\end{remark}

\subsubsection{Canonical integral models}
\label{sssection_can_integral_models}

PEL moduli problems don't just give canonical models of Shimura varieties,
they also give models over various localizations of the ring of integers
of the reflex field, at least when $\Ker^1(\Q,G)$ is trivial.
Here are some things that we can learn from this example:
\begin{itemize}
\item[(1)] The ring of integers over which we can expect to have a ``good''
integral model depends on the level $K$. More precisely, to have a
good integral model defined over $\Of_{F,(p)}$, we need $K$ to be of
the form $K^pK_p$, where $K^p\subset G(\Af^p)$ and $K_p$ is a hyperspecial
maximal compact subgroup of $G(\Q_p)$. Likewise, the Hecke correspondences
that will extend to finite \'etale morphisms between integral models
over $\Of_{F,(p)}$ are the ones that are trivial at $p$, i.e. defined
by elements of $G(\Af^p)$.

\item[(2)] We need a notion of what a ``good'' integral model is.
If the Shimura varieties for $(G,h)$ are compact, then we can just ask
for the integral model to be projective smooth over the localization of
$\Of_F$ that we are using. But this does
not suffice in the noncompact case. 

\end{itemize}

\begin{remark}
In these notes, we will only talk about integral models at places
of good reduction, so we want these models to be smooth over
the base ring and this is why we impose the conditions of (1). 
There is of course a (much more difficult and very interesting) theory for
places of bad reduction.

\end{remark}

To solve problem (2), Milne suggested only looking at models with a certain
extension property. Let $(G,h)$ be a Shimura datum, let $F=F(G,h)$, let
$p$ be a prime number at which $G$ is unramified, and let $K_p\subset G(\Q_p)$
be a hyperspecial maximal compact subgroup.
\footnote{The condition on $p$ means that $G$ extends to a reductive
group scheme $\Gf$ over $\Z_p$, and then we can take $K_p=\Gf(\Z_p)$.}
For every level $K$, write $M_K=M_K(G,h)$.
We want to define the notion of 
a canonical integral model $\Mcal_{K_p}$ over $\Of_{F,(p)}$
of the projective system $(M_{K_pK^p})_{K^p\subset G(\Af^p)}$, or of its
limit $M_{K_p}$.
The idea, first suggested by Milne in \cite{Milne-LR} (see also Moonen's paper
\cite{Moo}), is to require that, for every $S$ in a class of
``admissible test schemes'' over $\Of_{F,(p)}$, any morphism
$S\otimes_{\Of_{F,(p)}}F\ra M_{K_p}$ should extend to a morphism
$S\ra\Mcal_{K_p}$. The problem is to decide what class of admissible test
schemes one should use. We will follow Kisin's presentation
in \cite{Kisin-IM}.

\begin{definition}[See \cite{Kisin-IM} 2.3.7.]
\label{def_integral_model}
A \emph{canonical integral model} of the projective system
$(M_{K_pK^p})_{K^p\subset G(\Af^p)}$ 
(or of its limit $M_{K_p}$) over
$\Of_{F,(p)}$ is a projective system $(\Mcal_{K_pK^p})_{K^p\subset G(\Af^p)}$ of
smooth $\Of_{F,(p)}$-schemes with finite \'etale transition maps, given with
finite \'etale morphisms $T_{g,K^p,{K'}^p}:
\Mcal_{K_p{K'}^p}\ra\Mcal_{K_pK^p}$ for all $g\in G(\Af^p)$ and
$K^p,{K'}^p$ open compact subgroups of $G(\Af^p)$ such that ${K'}^p\subset
K^p\cap gK^pg^{-1}$, and with an isomorphism of projective systems
$\iota:(\Mcal_{K_pK^p})\otimes_{\Of_{F,(p)}}F\iso(M_{K_pK^p})$,
such that:
\begin{itemize}
\item[(a)] The morphisms $T_{g,K^p,{K'}^p}$ satisfy the analogues of conditions
(ii), (iii) and (iv) of Definition~\ref{def_model}, and they correspond
to the morphisms $T_{g,K_pK^p,K_p{K'}^p}$ between canonical models by the
isomorphism $\iota$ (in other words, $\iota$ is $G(\Af^p)$-equivariant);
\item[(b)] The scheme $\Mcal_{K_p}:=\varprojlim_{K^p}\Mcal_{K_pK^p}$ 
satisfies the following extension property: 
if $S$ is a regular formally smooth $\Of_{F,(p)}$-scheme, then
any morphism $S\otimes_{\Of_{F,(p)}}F\ra M_{K_p}$ extends to a morphism
$S\ra\Mcal_{K_p}$.

\end{itemize}
\end{definition}

As in the definition of canonical models, the properties listed in this definition 
are redundant, see Remark~\ref{remark_CN}.

In particular, by applying the extension property with $S=\Mcal_{K_p}$, we
see that integral canonical models are unique up to unique isomorphism.
Now the problem is existence. For PEL type Shimura varieties of type
$A$ or $C$ satisfying the condition that $\Ker^1(\Q,G)$ is trivial, the PEL
moduli problem will give a canonical integral model, though that is not
trivial: the main ingredient is the N\'eron-Ogg-Shafarevich criterion
on the good reduction of abelian varieties. For more details, see Theorem~2.10
of~\cite{Milne-LR}, and also Corollary~3.8 and Example~3.12(ii) of ~\cite{Moo}.

But what about other
Shimura varieties ?

\subsection{Hodge type and abelian type Shimura varieties}

\subsubsection{Hodge type Shimura varieties}
\label{sssection_Hodge_type}

For every $d\geq 1$, we denote by $(\GSp_{2d},h_d)$ the Siegel Shimura
datum of Example~\ref{ex_SD_GSp}.

The following condition was introduced in Deligne's paper
\cite{De-TC} (Section~2.3) and named in Milne's paper
\cite{Milne-CM} (at the end of Section~3).

\begin{definition}
A Shimura datum $(G,h)$ is \emph{of Hodge type} if there exists an integer
$d\geq 1$ and a morphism
of Shimura data $u:(G,h)\ra(\GSp_{2d},h_d)$ such that the underlying
morphism of algebraic groups $G\ra\GSp_{2d}$ is injective. In this case, we also
say that the corresponding Shimura varieties are of Hodge type.

\end{definition}

\begin{example}
\begin{itemize}
\item[(1)] Every Shimura datum of PEL type is of Hodge type, pretty much
by definition (or by Proposition~2.3.2 of~\cite{De-TC}): 
if $(\Of,*,\Lambda,\lan\cdot,\cdot\ran,h)$ is a PEL datum
and $G$ is the corresponding group scheme, then $G$ embeds into the group
scheme $H$ defined by
\[H(R)=\{g\in\End_R(\Lambda\otimes_\Z R)\mid\exists c(g)\in R^\times,\ 
\lan g(\cdot),g(\cdot)\ran=c(g)\lan\cdot,\cdot\ran\}\]
for every commutative ring $R$. As the alternating pairing $\lan\cdot,\cdot\ran$
is nondegenerate by condition~(5)(b) of Definition~\ref{def_PELD},
and as all nondegenerate alternating pairings on $\Lambda\otimes_\Z\Q$ are
equivalent, we have $H_\Q\simeq\GSp_{2d,\Q}$ for
$2d=\dim_\Q(\Lambda\otimes_\Z\Q)$.
Let $h'$ be the composition of $h:\Se\ra G_\R$ and of
the embedding $G_\R\ra H_\R=\GSp_{2d,\R}$. We have seen in the discussion
after Proposition~\ref{prop_PELD_to_PELSD} that $h:\C\ra\End_\R(\Lambda
\otimes_\Z\R)$ induces a Hodge structure of type $\{(-1,0),(0,-1)\}$ on
$\Lambda\otimes_\Z\C$,
\footnote{Which means that, in the decomposition
$\Lambda\otimes_\Z\C=\bigoplus_{p,q\in\Z}V^{p,q}$ induced by $h$, we have
$V^{p,q}=0$ unless $(p,q)\in\{(-1,0),(0,-1)\}$.}
so $h'$ satisfies condition (b) of Definition~\ref{def_SD}. It also
satisfies condition (c) of loc. cit. because the $\R$-bilinear pairing
$\lan\cdot,h(i)(\cdot)\ran$ on $\Lambda\otimes_\Z\R$ is symmetric
definite positive. Finally, for every $a\in\R^\times$, the element
$h(a)$ of $\End_\R(\Lambda\otimes_\Z\R)$ is $a\cdot\id$ (because $h$ is
a morphism of $\R$-algebras), so $h'(a)=a I_{2d}\in\GSp_{2d}(\R)$. 
This implies that $h'$ and $h_d$ are conjugated by $\GSp_{2d}(\R)$.

\item[(2)] The list of groups $G$ that have Shimura data of Hodge type
is given (at least in theory) in Section~2.3 of \cite{De-TC}. For
example, the group $G$ can be of type $B$, while that is not possible
for PEL type Shimura data.

\end{itemize}
\end{example}

The following result is due to Deligne; it follows from
Corollaire~5.7 of \cite{Deligne-SV}, which was already cited as
Proposition~\ref{prop_construction_cm}.

\begin{prop}
Every Shimura variety of Hodge type admits a canonical model.

\end{prop}

\begin{remark}
There is a general philosophy that Shimura varieties should be moduli
spaces of motives (the conditions that we put on a Shimura
datum $(G,h)$ are basically there to force the $G(\R)$-conjugacy class
of $h$ to be a parameter space for Hodge structures); see Section~3 of
\cite{Milne-CM} for more precise hopes.

For $(G,h)$ of Hodge type, 
we are a bit closer to that hope: Milne has proved that
$M_K(G,h)(\C)$ is a moduli space of abelian varieties with
Hodge cycles of a certain type and level structure (see
Theorem~3.11 of~\cite{Milne-CM}). As Hodge cycles on complex
abelian varieties are absolute by a theorem of Deligne
(see Theorem~2.11 of Chapter~I of~\cite{DMOS}),
we can also see the action of $\Aut(\C/F)$ (where $F=F(G,h)$) on
this modular interpretation.
So far, this has not allowed people to give a moduli interpretation of
the integral models of the Shimura variety of $(G,h)$, but
it does help with the construction of integral models, that we now discuss.

\end{remark}

Let $(G,h)$ be a Shimura datum of Hodge type, and let $u:G\ra\GSp_{2d}$ be
an injective morphism inducing a morphism of Shimura data $(G,h)\ra
(\GSp_{2d},h_d)$. Let $F=F(G,h)$ be the reflex field of $(G,h)$.
We fix a prime number $p$ such that $G_{\Q_p}$ extends to
a reductive group scheme $\Gf$ over $\Z_p$. 
To simplify the presentation,
we will assume that the embedding $G_{\Q_p}\ra\GSp_{2d,\Q_p}$ extends
to an embedding $\Gf\ra\GSp_{2d,\Z_p}$, though that is not necessary.
We set $K_p=\Gf(\Z_p)$ and $K'_p=\GSp_{2d}(\Z_p)$.

For each sufficiently small open compact
subgroup $K^p$ of $G(\Af^p)$, we fix an open compact subgroup
${K'}^p$ of $\GSp_{2d}(\Af^p)$ such that, setting $K=K_pK^p$ and
$K'=K'_p{K'}^p$, the morphism $u$ defines a closed immersion
$M_K(G,h)\ra M_{K'}(\GSp_{2d},h_d)_F$ (this is possible by Proposition~1.15
of~\cite{Deligne-SV}). 

Let $\Mcal_{K'}$ be the model of $M_{K'}(\GSp_{2d},h_d)$ over
$\Of_{F,(p)}$ given by base change from its canonical integral model over
$\Z_{(p)}$.
We denote by $\Mcal_K$ the normalization of the closure of the image
of $M_K(G,h)$ in $\Mcal_{K'}\supset M_{K'}(\GSp_{2d},h_d)_F$.
Kisin proved the following result.

\begin{theorem}[Theorem~2.3.8 of~\cite{Kisin-IM}]
Suppose that $p>2$.
Then $\varprojlim_{{K}^p}\Mcal_{K_pK^p}$ is a canonical integral
model of the Shimura variety of $(G,h)$.

\label{Kisin-IM1}
\end{theorem}

In particular, the schemes $\Mcal_K$ do not depend on the choice of
${K'}^p$ or on the embedding $G\ra\GSp_{2d}$.

\begin{remark}
In fact, Theorem~2.3.8 of~\cite{Kisin-IM} is more general, 
and gives a construction
of a canonical integral model without the assumption on the embedding
$G_{\Q_p}\ra\GSp_{2d,\Q_p}$. It even allows the case $p=2$ under some conditions.

\end{remark}

\begin{remark}
We have given the original construction of Kisin, but Yujie Xu recently proved
in~\cite{Xu_Y}
that the normalization step is unnecessary, so $\Mcal_K$ is just a closed
subscheme of $\Mcal_{K'}$.

\end{remark}

\subsubsection{Abelian type Shimura varieties}
\label{sssection_abelian_type}

\begin{definition}
Let $(G,h)$ be a Shimura datum. We say that $(G,h)$ is
\emph{of abelian type} if there exists a Shimura datum of Hodge type
$(G_1,h_1)$ and a central isogeny $G_{1,\der}\ra G_\der$ that induces an
isomorphism of Shimura data $(G_{1,\ad},h_{1,\ad})\iso(G_\ad,h_\ad)$, where
$h_{1,ad}$ (resp. $h_\ad$) is the composition of $h_1$ (resp.
$h$) and of the quotient morphism $G_1\ra G_{1,\ad}$ (resp. $G\ra G_\ad$).
In this case, we also say that the corresponding Shimura varieties are
of abelian type.

\end{definition}

Another way to formulate the definition is to say that a Shimura variety
$M_K(G,h)(\C)$ is of abelian type if all its connected components are
finite quotients of connected components of Shimura varieties of
Hodge type (see Section~9 of~\cite{De-TC}). Deligne has classified all
connected Shimura varieties of abelian type in 2.3 of~\cite{De-TC}.
The very rough upshot is that all Shimura data $(G,h)$ with $G$ of type
$A$, $B$ and $C$ are of abelian type; if $G$ is of type $D$, it's
complicated, and if $G$ is of type $E_6$ or $E_7$, then the Shimura
datum is never of abelian type. See page 61 of \cite{Lan} for more details.

\begin{theorem}[Deligne, see Corollaire~2.7.21 of~\cite{De-TC}]
Let $(G,h)$ be a Shimura datum of abelian type.
Then the Shimura variety of $(G,h)$ admits a canonical model.

\end{theorem}

In fact, Deligne reduces the construction of a canonical model of
$M_K(G,h)(\C)$ to that of canonical models of its connected components
(over finite extensions of $F(G,h)$). See for example Corollaire~2.7.18
of~\cite{De-TC}.

\begin{remark}
Shimura varieties of abelian type are \emph{not} moduli spaces
of abelian varieties in general. However, Milne proved in~\cite{Milne-SVM}
that they are moduli spaces
of motives if $h\circ w:\GL_{1,\R}\ra G_\R$ is defined over $\Q$, and this
leads to a more direct 
proof of existence of their canonical models even without that
condition. 

\end{remark}

\begin{theorem}[Milne, see Theorem~3.31 of~\cite{Milne-SVM}; see
also Brylinski's paper \cite{Br}]
Let $(G,h)$ be a Shimura datum of abelian type
such that $h\circ w:\GL_{1,\R}\ra G_\R$ is defined over $\Q$. Then
each $M_K(G,h)$ is a moduli space of abelian motives (over the reflex field
of $(G,h)$).

\end{theorem}

By reducing to the case of Shimura varieties of Hodge type, Kisin
was able to prove the existence of canonical integral models of
Shimura varieties of abelian type for $p>2$. The case $p=2$ was then
settled by Kim and Madapusi.

\begin{theorem}[Corollary~3.4.14 of~\cite{Kisin-IM} and Theorem~1 of~\cite{Kim-Mada}]
Let $(G,h)$ be a Shimura datum of abelian type, let $p$ be a prime number
such that $G_{\Q_p}$ extends to a reductive group scheme $\Gf$ over
$\Z_p$, and let $K_p=G(\Z_p)$. Then $M_{K_p}:=\varprojlim_{K^p\subset G(\Af^p)}
M_{K_p K^p}(G,h)$ admits a canonical integral model over $\Of_{F,(p)}$, where
$F=F(G,h)$.

\label{Kisin-IM2}
\end{theorem}

\subsection{General Shimura varieties}

Let $(G,h)$ be a Shimura datum that is not of abelian type. Then we
know very little, but we do know that canonical models exist.
This result was first announced by Milne
(based on earlier results of Kazhdan and Borovoi), but there was a gap
in the proof, which was fixed by Moonen.

\begin{theorem}[Deligne, Borovoi, Milne-Shih, Milne, Moonen, 
cf. \cite{Moo} Section~2.]
\label{thm_ECM}
Let $(G,h)$ be a Shimura datum. Then the Shimura variety of $(G,h)$ admits
a canonical model.

\end{theorem}

In fact, Section~2 of~\cite{Moo} contains a good summary of the
different construction methods of canonical models. The proof in
the general case does not proceed by reduction to the case of Siegel
modular varieties (unlike the previous proofs in the abelian type case),
but uses results of Borovoi (see \cite{Bor1} and \cite{Bor2}), 
Deligne, Milne-Shih (Chapter~V of the book~\cite{DMOS})
and Milne (\cite{Milne-CMI}) on a conjecture
of Langlands, that says that a conjugate of a Shimura variety over $\C$
by an automorphism of $\C$ is still a Shimura variety (Langlands's conjecture
is much more precise than this, see for example Theorem~2.14 of~\cite{Moo}).

\subsection{Kottwitz's simple Shimura varieties}
\label{ssection_SSV}

After considering more and more complicated Shimura varieties
in the previous subsections, we will now introduce a very simple
PEL family, that has been studied by Kottwitz in~\cite{K-SSV}; their
PEL datum already appeared in Example~\ref{ex_PEL_Kottwitz}.
These
Shimura varieties are simple for several reasons:
\begin{itemize}
\item they are compact;
\item they are PEL of type $A$, hence moduli spaces of abelian schemes with
extra structures;
\item their reductive group ``has no endoscopy'' 
(see Proposition~\ref{prop_no_end} for a precise statement).

\end{itemize}

\subsubsection{Inner forms of unitary groups}
\label{sssection_inner_U}

We fix a totally real extension $F_0$ of $\Q$ and
a totally imaginary quadratic extension $F$ of $F_0$.
Such an extension $F$ of $\Q$ is called a \emph{CM extension}.
We denote by $z\mapsto\overline{z}$ the nontrivial element of
$\Gal(F/F_0)$.

Let $n$ be a positive integer. The quasi-split unitary group
$\U^*(n)$ over $F_0$ is defined to be the unitary group of the
Hermitian $F_0$-space $F^n$, with the form
\[((x_1,\ldots,x_n),(y_1,\ldots,y_n))\mapsto\sum_{i=1}^n x_i\overline{y_{n+1-i}}.\]
In other words, for every commutative $F_0$-algebra $R$, we have 
\[\U^*(n)=\{g\in\GL_n(F\otimes_{F_0}R)\mid {}^t\overline{g}Jg=J\},\]
where $J$ is the $n\times n$ antidiagonal matrix with all nonzero coefficients
equal to $1$.

We want to describe all inner forms of $\U^*(n)$. Remember that, if $G$ and
$H$ are algebraic groups over a field $k$, we say that they are
\emph{inner forms} of each other if there exists an isomorphism
$\varphi:G_{\overline{k}}\iso H_{\overline{k}}$ such that, for every
$\sigma\in\Gal(\overline{k}/k)$, the automorphism
$\varphi^{-1}\circ{}^\sigma\varphi$ of $G(\overline{k})$ is inner
(${}^\sigma\varphi$ is the isomorphism $\sigma\varphi\sigma^{-1}:
G_{\overline{k}}\iso H_{\overline{k}}$, where $\sigma$ acts on $G_{\overline{k}}$ and $H_{\overline{k}}$
via its action on $\overline{k}$).
Inner forms of $G$ are in bijection with elements of
$\H^1(k,G_\ad):=\H^1(\Gal(\overline{k}/k),G_\ad)$, where $G_\ad=G/Z(G)$
and $Z(G)$ is the center of $G$.
See Section~III.1 of Serre's book~\cite{Serre-CG} and Section~3 of
Springer's Corvallis notes~\cite{Springer-C} for more about (inner) forms.

Here, observing that $\U^*(n)_\ad$ is the group of automorphisms of
the couple formed by the central simple algebra $M_n(F)$ over $F$ and
the involution $g\mapsto J{}^t\overline{g}J$, 
we see that inner forms of $\U^*(n)$
are all of the form $\U(B,*)$, where $B$ is a central simple algebra 
of dimension $n^2$ over
$F$ and $*$ is an involution on $B$ extending the involution
$z\mapsto\overline{z}$ on $F$, and $\U(B,*)$ is the $F_0$-group defined
by
\[\U(B,*)(R)=\{g\in B\otimes_{F_0}R\mid gg^*=1\},\]
for every commutative $F_0$-algebra $R$. 
As the notation indicates, the group $\U(B,*)$ is a unitary group.
More precisely, we can write $B=M_{n/m}(D)$ where $m$ is a divisor of $n$ and
$D$ is a division algebra of dimension $m^2$ over $F$; if we choose
a simple $B$-module $V$, then we can take $D=\End_B(V)^\op$ and see
$V$ as a $D^\op$-module. By Theorem~4.1 of the lecture notes~\cite{Tignol} of Tignol,
there exists an involution on $D^\op$ that restricts to
$z\mapsto\overline{z}$ on $F$ and a nondegenerate bilinear form $\lan.,.\ran$
on the $D^\op$-module $V$, Hermitian with respect to that involution, such
that $*$ is equivalent to the involution $\ddag$ on $B\simeq\End_{D^\op}(V)$
corresponding to the Hermitian form $\lan.,.\ran$, i.e. defined by the
condition $\lan bx,y\ran=\lan x,b^\ddag y\ran$ for $b\in B$ and
$x,y\in V$.

On the other hand, if $G$ is an inner form of $\U^*(n)$, then
$G_{F_{0,v}}$ is an inner form of $\U^*(n)_{F_{0,v}}$ for every place $v$ of
$F_0$, and $G_{F_{0,v}}$ is isomorphic to $\U^*(n)_{F_{0,v}}$ itself for all but
finitely many $v$.
Inner forms over local fields are easier to classify (because
the absolute Galois groups of local fields are simpler); but then we must
be able to decide when a family of inner forms of the $\U^*(n)_{F_{0,v}}$ comes
from a ``global'' inner form of $\U^*(n)$ (defined over $F_0$). Using
Galois cohomology calculations and the calculation of the Brauer groups
of local and global fields, we get the following two propositions
(see for example Section~2
of Clozel's paper \cite{Clo-RG}, and see Gille and Szamuely's
book~\cite{GS} for more information on central simple algebras and Brauer
groups):

\begin{prop}\label{prop_IFU1}
Let $v$ be a place of $F_0$.
\begin{enumerate}
\item Suppose that $v$ is finite and does not split in $F$. If $n$ is odd,
then the only inner form of $\U^*(n)_{F_{0,v}}$ is $\U^*(n)_{F_{0,v}}$ itself
(up to isomorphism).
If $n$ is even, then there are two isomorphism classes of inner
forms of $\U^*(n)_{F_{0,v}}$.

\item Suppose that $v$ splits in $F$ (in particular, $v$ is finite), and
let $w$ be a place of $F$ above $v$. Then
$\U^*(n)_{F_{0,v}}\simeq\GL_{n,F_{0,v}}\simeq\GL_{n,F_w}$, and its inner forms
are (up to isomorphism)
the groups $\GL_m(D)$, for $m$ dividing $n$ and $D$ a central division
algebra over $F_w$ of dimension $(n/m)^2$.

\item Suppose that $v$ is infinite, hence a real place of $F_0$.
Then the inner forms of $\U^*(n)_{F_{0,v}}$ are (up to isomorphism)
the real unitary groups
$\U_{p,q}$ of signature $(p,q)$, for $p+q=n$, and we have
$\U_{p,q}\simeq\U_{r,s}$ if and only if $(r,s)=(p,q)$ or $(r,s)=(q,p)$.

\end{enumerate}
\end{prop}

Suppose that $n$ is even, let $v$ be a place of $F_0$ and
let $G$ be an inner form of $\U^*(n)_{F_{0,v}}$.
\begin{itemize}
\item If $v$ is finite and does not split in $F$, set
$\epsilon(G)=1$ if $G\simeq\U^*(n)_{F_{0,v}}$ and $\epsilon(G)=-1$ otherwise;

\item If $v$ is finite and splits in $F$, set
$\epsilon(G)=(-1)^m$ if $G\simeq\GL_m(D)$
with $m$ dividing $n$ and $D$ a central division
algebra over $F_w$ of dimension $(n/m)^2$
(note that $m$ only depends on $v$, not on $w$);

\item If $v$ is infinite, set
$\epsilon(G)=(-1)^{n/2-p}$ if $G\simeq\U_{p,q}$.

\end{itemize}

\begin{prop}\label{prop_IFU2}
For every place $v$ of $F_0$, let $G_v$ be
an inner form of $\U^*(n)_{F_{0,v}}$.
Suppose that $G_v\simeq
\U^*(n)_{F_{0,v}}$ for all but finitely many $v$.
\begin{enumerate}
\item If $n$ is odd, there exists an inner form $G$ of $\U^*(n)$
such that $G_{F_{0,v}}\simeq G_v$ for every $v$.

\item If $n$ is even, there exists an inner form $G$ of $\U^*(n)$
such that $G_{F_{0,v}}\simeq G_v$ for every $v$ if and only if
$\prod_v\epsilon(G_v)=1$.

\end{enumerate}
\end{prop}

\begin{remark}\label{rmk_IFU}
Let $G$ be an inner form of $\U^*(n)$. We know that $G\simeq\U(B,*)$, with
$B$ a central simple algebra over $F$ and $*$ 
an involution on $B$ extending the involution
$z\mapsto\overline{z}$ on $F$. In the next subsubsection, it will be of
interest to us to know when $B$ is a division algebra.
Let $v$ be a place of $F_0$, and let $G_v=G_{F_{0,v}}$.
\begin{itemize}
\item If $v$ does not split in $F$, let $w$ be a place of $F$ above $v$.
By the discussion at the beginning of~\ref{sssection_inner_U} (i.e.
Theorem~4.1 of~\cite{Tignol}, which applies to any field), there
exists a division algebra $D_v$ over $F_w$ with an involution
extending the nontrivial element of $\Gal(F_w/F_{0,v})$ and a
Hermitian module $(V, \lan.,.\ran)$ over $D_v$ such that $G_v$ is the
unitary group of $(V,\lan.,.\ran)$. As $F_w$ is a local field,
we must have $D_v=F_w$, and then $V\simeq F_w^n$ and the Hermitian
form $\lan.,.\ran$ is characterized  by its determinant (up to
isomorphism).

\item If $v$ splits in $F$, then we have
$G_v\simeq M_{m_v}(D_v)$, with $m_v$ dividing $n$
and $D_v$ a central division algebra of dimension $(n/m_v)^2$ over $F_{0,v}$.
\end{itemize}

We can now deduce from the classification of central simple algebras
over $F$ that $B$ is a division algebra if
the gcd of the family $(m_v)$ is equal to $1$, because it forces the order
of the class of $B$ in the Brauer group of $\Q$ to be equal to $n$.
(See Corollary~6.5.4 of~\cite{GS} for the classification
of central simple algebras over a number field.) 
The simplest way to make sure that this condition is satisfied is to
take one of the $m_v$ equal to $1$, for example to take $G_v$ of the form
$D_v^\times$ for a place $v$
of $F_0$ split in $F$, where $D_v$ is a central division
algebra over $F_{0,v}$.

\end{remark}

\subsubsection{Simple Shimura varieties}
\label{subsubsection_SSV}

\begin{definition}[See \S1 of~\cite{K-SSV}.]
\label{def_KSS}
A \emph{Kottwitz simple Shimura variety} is a Shimura variety defined
by the Shimura datum $(G_\Q,h)$ associated to a PEL datum
$(\Of,*,\Lambda,\lan\cdot,\cdot\ran,h)$
as in Example~\ref{ex_PEL_Kottwitz}.

\end{definition}

In particular:
\begin{itemize}
\item[(a)] $D:=\Of\otimes_\Z\Q$ is a division algebra with center a CM extension
$F$ of $\Q$;
\item[(b)] $*$ extends the nontrivial automorphism of $F/F_0$, where
$F_0$ is the maximal totally real subextension of $F$;
\item[(c)] $\Lambda=\Of$;
\item[(d)] $\lan x,y\ran=\Tr_{\Of/\Z}(xby^*)$ for all $x,y\in\Lambda$, where $b\in D^\times$ is such that
$b^*=-b$.

\end{itemize}

We have $\End_\Of(\Lambda)=\Of^\op$ (acting by right multiplication), so,
for every commutative ring $R$,
\[G(R)=\{g\in(\Of^\op\otimes_\Z R)^\times\mid gg^*\in R^\times\}.\]
We have a morphism $c:G\ra\GL_1$ sending $g\in G(R)$ to $gg^*\in\GL_1(R)$, and, if
we denote its kernel by $G_0$, then $G_{0,\R}$ is a product of unitary groups
of the form $\U(p,q)$. We gave the formula for $h$ in Example~\ref{ex_PEL_Kottwitz}: if
$z\in\C$ has absolute value $1$, then the projection of $h(z)$ to the $\U(p,q)$ factor
of $G_0(\R)$ is $\begin{pmatrix}z I_p & 0 \\ 0 &\overline{z} I_q\end{pmatrix}$.
(There is a choice here, as the couple $(p,q)$ is only determined up to
order.)

\quash{
Denote by $\Res_{F_0/\Q}$ the functor of Weil restriction of scalars
from $F_0$ to $\Q$; in particular, if $H$ is an algebraic group over
$F_0$, then $\Res_{F_0/\Q}H$ is the algebraic group over $\Q$ defined
by $\Res_{F_0/\Q}H(R)=H(R\otimes_\Q F_0)$ for every commutative $\Q$-algebra
$R$. Then
\[G_{0,\Q}=\Res_{F_0/\Q}\U,\]
with $\U$ the algebraic group over $F_0$ defined by
\[U(R)=\{g\in(D^\op\otimes_{F_0}R)^\times\mid gg^*=1\}\]
for every $F_0$-algebra $R$.

The dimension of $D$ over $F$ is of the form $n^2$ with $n\in\Nat$, and
$\U$ is an inner form
of the quasi-split unitary group $\U^*(n)$ over $F_0$ defined
by the extension $F/F_0$. By Propositions~\ref{prop_IFU1}
and~\ref{prop_IFU2} and Remark~\ref{rmk_IFU}, we have a description
of such inner forms.
In particuler, we know that
\[G_{0,\R}\simeq\prod_{\tau\in\Phi}\U_{p_\tau,q_\tau},\]
where $\Phi$ is the set of real places of $F_0$ and, for every
$\tau\in\Phi$, $p_\tau$ and $q_\tau$ are nonnegative integers such
that $p_\tau+q_\tau=n$. Definition~\ref{def_KSS} says nothing
about $h:\C^\times\ra G(\R)$, but there is a unique choice that makes
$(G_\Q,h)$ a Shimura datum:
if $z\in\C^\times$ is of norm $1$, we take
\[h(z)=\left(\begin{pmatrix} z I_{p_\tau} & 0 \\ 0 & \overline{z}I_{q_\tau}
\end{pmatrix}
\right)_{\tau\in\Phi}\in\prod_{\tau\in\Phi}\U_{p_\tau,q_\tau}\simeq
G_0(\R).\]
Then there is a unique extension $h$ to a morphism $h:\C^\times\ra G(\R)$
such that $h(a)=a\cdot\id$ if $a\in\R^\times$.
}

Note that we can choose the signatures of $G_{0,\R}$ arbitrarily, by manipulating
what happens at finite places. On the other hand, the Shimura varieties
of $(G_\Q,h)$ are always compact, because of the following lemma.
Here we use the fact that the Shimura varieties defined by a Shimura datum
$(G,h)$ are compact if and only if $G_\der$ has $\Q$-rank $0$, which was
recalled in 1.1 (see for example~\cite[Theorem~5.10]{Ji}).

\begin{lemma}
The group $G_\der$ is of $\Q$-rank $0$.

\end{lemma}

\begin{proof}
Let $N:D^\times\ra F^\times$ be the reduced norm, seen as a morphism of
algebraic groups over $F$.
We have $G_\der=\Res_{F_0/\Q}\SU$, where $\SU=\{g\in\U\mid N(g)=1\}$,
\footnote{We can do the calculation after base changing to $\C$, where it
reduces to the fact that $\GL_{n,\der}=\SL_n$.}
so $G_{\der,F}\simeq(\Ker N)^{[F_0:\Q]}$.
Let $T$ be a maximal torus of $G_\der$. Then $T_F$ is a maximal torus of
$(\Ker N)^{[F_0:\Q]}$, hence a product of maximal tori of $\Ker N$. 
Let $T'$ be a maximal torus of $\Ker N$.
As $D$ is a division algebra, there exists a degree $n$ 
extension $F'$ of $F$ such that
$T'=\{x\in {F'}^\times\mid N_{F'/F}(x)=1\}$, and so the
maximal split subtorus of $T'$ is trivial. This implies that the maximal split
subtorus of $T$ is trivial.

\end{proof}

We now discuss endoscopy.

\vspace{.3cm}

\subsubsection{Endoscopy}

In its simplest form, endoscopy is the following phenomenon: let
$G$ be an algebraic group over a field $k$. Then the $G(\overline{k})$-
conjugacy classes in $G(k)$ can be larger than the $G(k)$-conjugacy classes.
Some vocabulary: $G(\overline{k})$-conjugacy classes in $G(k)$ are
often called \emph{stable conjugacy classes}, and $G(\overline{k})$-conjugate
elements are called \emph{stably conjugate}. Actually we are cheating here,
and our definition is only correct if the semisimple elements of $G$ have
connected centralizers; the correct definition is in Section~3 of
Kottwitz's paper~\cite{K-RCC}. Note that semisimple elements of $G$ have
connected centralizers if $G_\der$ is simply connected: this is a result
of Steinberg (see Corollary~8.5 of~\cite{Steinberg})
for $G$ semisimple, and Kottwitz explains in Section~3 of~\cite{K-RCC}
how to extend it to the case of a reductive group.

\begin{examples}\label{ex_conj_class}
\begin{itemize}
\item[(1)] Suppose that $G=\GL_n$. It is then a classical exercise that
any two elements of $G(k)$ that are $G(\overline{k})$-conjugate are
actually $G(k)$-conjugate. We say that $\GL_n$ has no endoscopy.

\item[(2)] We can generalize (1) to inner forms of $\GL_n$, i.e.
algebraic groups of the form $B^\times$, where $B$ is a central simple
algebra over $k$.

\item[(3)] Take $k=\R$ and $G=\SL_2$. As $\GL_2(\C)=\C^\times\cdot\SL_2(\C)$,
elements of $G(\R)$ are $G(\C)$-conjugate if and only if they are
$\GL_2(\C)$-conjugate. For example, the matrices
$\begin{pmatrix} 0 & 1 \\ -1 & 0\end{pmatrix}$ and
$\begin{pmatrix} 0 & -1 \\ 1 & 0\end{pmatrix}$
are in the same $\SL_2(\C)$-conjugacy class, but we can check by a
direct calculation that they are not in the same $\SL_2(\R)$-conjugacy class.
So $\SL_2$ has endoscopy, and we can generalize that example
to $\SL_n$ for $n\geq 2$.

\item[(4)] If $k=\R$ and $G=\U_n:=\{g\in\GL_n(\C)\mid g^*g=I_n\}$, then
again it is a classical exercise to check that $\U_n(\C)$-conjugacy classes
in $\U_n(\R)$ coincide with $\U_n(\R)$-conjugacy classes.

\end{itemize}
\end{examples}

\begin{remark}
We care about stable conjugacy classes because, in the Langlands
philosophy, groups are related via their $L$-groups (see
\ref{ssection_Langlands_group} and \ref{sssection_Langlands_corr}).
As inner forms have the same $L$-groups, this means that we should
be able to move information between inner forms; but we cannot compare
conjugacy classes in two inner forms, only stable conjugacy classes.
More generally, if $G$ and $H$ are algebraic groups over $k$
and there is a morphism ${}^LH\ra{}^LG$ between their $L$-groups,
then we can use this to transport stable conjugacy classes of regular semisimple
elements from $H$ to $G$.

\end{remark}

It quickly becomes tiring to calculate conjugacy classes by hand, so
we need more efficient methods to check for endoscopy.
Let
$\gamma\in G(k)$, and let $G_\gamma\subset G$ be the centralizer of
$\gamma$.

Let $\delta\in G(k)$ be stably conjugate to $\gamma$, and let
$g\in G(\overline{k})$ such that $\delta=g\gamma g^{-1}$. 
For every $\sigma\in\Gal(\overline{k}/k)$, we have
\[g\gamma g^{-1}=\delta=\sigma(\delta)=\sigma(g)\sigma(\gamma)\sigma(g)^{-1}=
\sigma(g)\gamma\sigma(g)^{-1},\]
hence $g^{-1}\sigma(g)\in G_\gamma(\overline{k})$.
So we get a $1$-cocyle $c:\Gal(\overline{k}/k)\ra G_\gamma(\overline{k})$,
$\sigma\mapsto g^{-1}\sigma(g)$, and we can check that the image 
$\inv(\gamma,\delta)$ of
this $1$-cocyle in $\H^1(k,G_\gamma)$ does not depend on the choice
of $g$. Moreover, the image of $\inv(\gamma,\delta)$ in
$\H^1(k,G)$ is trivial.

\begin{prop}(See Section~3 of~\cite{K-RCC}.)
Suppose that $\gamma$ is semisimple and that $G_\gamma$ is connected.
Then the map $\delta\mapsto\inv(\gamma,\delta)$ gives a bijection
from the set of stable conjugacy classes in the $G(k)$-conjugacy
class of $\gamma$ to $\Ker(\H^1(k,G_\gamma)\ra\H^1(k,G))$.

\end{prop}

This is particularly useful when $G_\der$ is simply connected and
$\gamma$ is regular and semisimple, as the centralizer of $\gamma$ is
then a maximal torus (because it is connected, as explained before
Example~\ref{ex_conj_class}).

\begin{example}\label{ex_end_2}
\begin{itemize}
\item[(1)] 
\footnote{I am very grateful to
Mikhail Borovoi for patiently explaining how to deduce this
example from Example~4.4 of his paper~\cite{Borovoi}. Any remaining mistakes
are all mine.}
Let $p\geq q\geq 0$ be integers.
Take $k=\R$ and $G=\U_{p,q}$, the unitary group of the Hermitian form
on $\C^{p+q}$ with matrix $\begin{pmatrix} I_p & 0 \\ 0 & -I_q\end{pmatrix}$,
where $I_r$ denotes as usual the identity matrix in $\GL_r(\Z)$.
One can show that any maximal torus $T$ of $G$ is isomorphic to
$\Se^r\times\U_1^s$, with $2r+s=p+q$. (Remember that $\Se=\Res_{\C/\R}\GL_{1,\C}$,
i.e. it is $\C^\times$ regarded as an algebraic group over $\R$.)
By Shapiro's lemma,
we have
\[\H^1(\R,\Se)=\H^1(\C,\GL_{1,\C})=0.\]
Using the exact sequence
\[1\ra\U_1\ra\Se\flnom{\mathrm{Nm}}\GL_{1,\R}\ra 1\]
where $\mathrm{Nm}:\C^\times\ra\R^\times$ is the map $z\mapsto z\overline{z}$,
we see that $\H^1(\R,\U_1)=\R^\times/
\mathrm{Nm}(\C^\times)=\Z/2\Z$.
So
\[\H^1(\R,T)=(\Z/2\Z)^s.\]

Let $T_0$ be the group of diagonal matrices in $G$. Then $T_0$ is a maximal
torus of $G$ and $T_0\simeq\U(1)^{p+q}$, 
so we get that $\H^1(\R,T_0)=(\Z/2\Z)^{p+q}$.
On the other hand, by the main result (Theorem~3.1) of Borovoi's
paper \cite{Borovoi}, we have
\[\H^1(\R,G)=\H^1(\R,T_0)/W_0(\R),\]
and $W_0$ is the Weyl group scheme $W_{T_0}=N_G(T_0)/Z_G(T_0)$. 
It is also proved in Section~3 of the same paper that
$W_0(\R)=W_0(\C)$, so $W_0(\R)$ is isomorphic to the symmetric
group $\Sgoth_{p+q}$. However, as permutation matrices are not
always inside $\U_{p,q}$, the group $\Sgoth_{p+q}$ does not act on $T_0$ by
permuting the entries.

We calculate the quotient $\H^1(\R,T_0)/\Sgoth_{p+q}$. Note that
$\H^1(\R,T_0)\simeq T_0(\R)_2$, where $T_0(\R)_2=\{\diag(z_1,\ldots,z_{p+q})\in
T_0(\R)\mid\forall i,\ z_i\in\{\pm 1\}\}$, so we need to understand the
action of $\Sgoth_{p+q}$ on $T_0(\R)_2$. 
We have $T_0(\R)_2=\{c_\Xi,\ \Xi\subset\{1,2,\ldots,p+q\}\}$, where, for
every subset $\Xi$ of $\{1,2,\ldots,p+q\}$, we write $c_\Xi$ for the
matrix $\diag(z_1,\ldots,z_{p+q})$ such that $z_i=-1$ for $i\in\Xi$ and
$z_i=1$ for $i\not\in\Xi$.
We regard $\Sgoth_p\times\Sgoth_q$ as 
a subgroup of $\Sgoth_{p+q}$ by identifying $\Sgoth_q$ to the permutation
group on $\{p+1,p+2,\ldots,p+q\}$. Then $\Sgoth_{p+q}$ is generated by
$\Sgoth_p\times\Sgoth_q$ and by the transposition $\sigma:=(1,p+1)$.
Also, elements of $\Sgoth_p\times\Sgoth_q$ act on $T_0(\R)_2$ by permuting
the coordinates, because the corresponding permutation matrices are
in $\U_{p,q}$. So two elements $c_\Xi$, $c_{\Xi'}$ of $T_0(\R)_2$
are in the same $\Sgoth_p\times\Sgoth_q$-orbit if and only if
$\card(\Xi\cap\{1,\ldots,p\})=\card(\Xi'\cap\{1,\ldots,p\})$ and
$\card(\Xi\cap\{p+1,\ldots,p+q\})=\card(\Xi'\cap\{p+1,\ldots,p+q\})$.

It remains to understand the action of $\sigma$. It will be represented
by an element of $N_{T_0}(\C)$, so we need a convenient representation
of $G(\C)$. For this we consider the $\R$-algebra embedding
$\Mcal:\C\ra M_2(\R)$, $a+ib\mapsto\begin{pmatrix}a & b \\ -b & a\end{pmatrix}$.
Applying this to every matrix entry, we get an embedding of real algebraic
groups
$G\subset\GL_{2(p+q),\R}$, hence also $G_\C\subset\GL_{2(p+q),\C}$.

Suppose that $(p,q)=(1,1)$. Then the matrix 
$n=i\begin{pmatrix} 0 & I_2\\ I_2 & 0
\end{pmatrix}\in\GL_4(\C)$ is in $N_{T_0}(\C)\subset G(\C)$ and
represents $\sigma$, so the action of $\sigma$ on $T_0(\R)_2$ is given
by $c\mapsto n^{-1}c\overline{n}$. It is then easy to see that $\sigma$ fixes
$c_{\{1\}}$ and $c_{\{2\}}$, and exchanges $c_\varnothing=1$ and $c_{\{1,2\}}=-1$.
So $\H^1(\R,G)$ has three elements.

The general case is similar: we write an explicit representative
$n\in N_{T_0}(\C)\subset\GL_{2(p+q)}(\C)$ of $\sigma$ by putting the
$I_2$ blocks in the correct positions. For every $\Xi\subset
\{1,2,\ldots,p+q\}$, write $p_\Xi=\card(\Xi\cap\{1,\ldots,p\})$
and $q_\Xi=\card(\Xi\cap\{p+1,\ldots,p+q\})$.
By a calculation similar to that of the
previous paragraph, we get that $c_\Xi$ and $c_{\Xi'}$ are in the same
$\Sgoth_{p+q}$-orbit if and only $\Xi,\Xi'$ satisfy the following three
conditions:
\[\left\{\begin{array}{l}
p_\Xi=p_{\Xi'}\\
q_\Xi=q_{\Xi'}\\
p_\Xi-q_\Xi=p_{\Xi'}-q_{\Xi'}.
\end{array}\right.\]
So a set of representations for the $\Sgoth_{p+q}$-orbits in
$T_0(\R)_2$ is given by $\{c_\varnothing\}\cup\{c^1_r,\ 1\leq r\leq p\}\cup
\{c^2_s,\ 1\leq s\leq q\}$, where $c^1_r=c_{\{1,2,\ldots,r\}}$ and
$c^2_s=c_{\{p+1,p+2,\ldots,p+s\}}$.

We deduce that $\H^1(\R,G)$ is of cardinality $p+q+1$, and that
\[\Ker(\H^1(\R,T_0)\ra\H^1(\R,G))=\{c_\Xi\mid p_\Xi=q_\Xi\}.\]

Alternatively, to calculate $\H^1(\R,G)$,
one can use the fact that $\H^1(\R, G)$ for $G=U_{p,q}$
classifies isomorphism classes of non-degenerate Hermitian forms on $\C^{p+q}$
(see Serre \cite{Serre-LF}, Section~X.2, Proposition~4),
and one can classify the isomorphism classes of such forms 
using the Hermitian version of Sylvester's law of inertia.

\item[(2)] We take $k$ to be a finite extension of $\Q_p$,
$E$ to be an unramified quadratic extension of $k$,
and $G=\U^*(n)_{E/k}$ to be the quasi-split unitary group defined by
that extension, i.e. the unitary group of the
Hermitian $k$-space $E^n$, with the form
\[((x_1,\ldots,x_n),(y_1,\ldots,y_n))\mapsto\sum_{i=1}^n x_i\overline{y_{n+1-i}},\]
where $x\mapsto\overline{x}$ is the nontrivial element of
$\Gal(E/k)$.
Any maximal torus $T$ of
$G$ is isomorphic to $(\Res_{E/k}\GL_1)^r\times\U(1)_{E/k}^s$, where
$\U(1)_{E/k}$ is the subgroup of norm one elements in $\Res_{E/k}\GL_1$
and $2r+s=n$.
As before, we have $\H^1(k,\Res_{E/k}\GL_1)=0$ and
$\H^1(k,\U(1)_{E/k})=k^\times/N_{E/k}(E^\times)=\Z/2\Z$. 

But now we have $\H^1(k,G)=\Z/2\Z$, and
the map $\H^1(k,T)\ra\H^1(k,G)$ is the sum map.
Indeed, the derived group
$G_\der=\SU^*(n)_{E/k}$ of $G$ is simply connected, so
$\H^1(k,G_\der)=0$ by a theorem of Kneser (see \cite{Kne}), and then
we use the fact that $G=G_\der\rtimes\U(1)_{E/k}$.

So we see that, for some maximal tori $T$, the set
$\Ker(\H^1(k,T)\ra\H^1(k,G))$ is not a singleton, hence there are
stable conjugacy classes containing more than one conjugacy class.
As in (1), the more anisotropic factors the torus $T$ has, the bigger
its $\H^1$.

\end{itemize}
\end{example}

The examples show that, if $G$ is the general unitary
group of a Kottwitz simple Shimura variety as in Definition~\ref{def_KSS},
then it is not reasonable to expect that $G$ will behave as $\GL_n$ and
have absolutely no endoscopy over any field. In fact, what Kottwitz actually
proved about these groups is the following:

\begin{prop}[Lemma~2 of~\cite{K-SSV} and Theorem~6.6 of~\cite{K-EST}]
\label{prop_no_end}
Let $\overline{\Ade}$ be the restricted product of the $\oQ_v$, for
$v$ a place of $\Q$. Let $\gamma$ be a semisimple element of
$G(\Q)$, and let $\delta$ be an element of $G(\Ade)$ that is
$G(\overline{\Ade})$-conjugate to $\gamma$. Then there exists an
element of $G(\Q)$ that is $G(\Ade)$-conjugate to $\delta$.
\footnote{Note however that this new
element of $G(\Q)$ may not be $G(\Q)$-conjugate to $\gamma$.}

\end{prop}

This implies that in the trace formula for $G$
(see~\ref{sssection_Arthur_TF}), we will be able to group the orbital
integrals on the geometric side by stable conjugacy class and obtain an
expression that is easier to transfer between groups. See Section~4
of~\cite{K-SSV}.

The proof of Proposition~\ref{prop_no_end} is a more complicated
global version of the Galois cohomology calculations of
Example~\ref{ex_end_2}. First Kottwitz reduces to the case of
$G_\der=\Res_{F_0/\Q}\U$, where $\U$ is the unitary group (over $F_0$)
given by $\U(R)=\{g\in(D\otimes_{F_0}R)^\times\mid gg^*=1\}$.
Then he reduces to a similar result for $\U$, now involving $F_0$-rational
points and the ring of adeles of $F_0$. The main point is that, if
$T$ is a maximal torus of $\U$, then $T_F$ is a maximal torus of $D^\times$,
and maximal tori of $D^\times$ are all of the form $\Res_{K/F}\GL_1$, for $K$
a degree $n$ extension of $F$. This allows us to control the Galois
cohomology of $T$.

\section{Lecture 3: the cohomology of Shimura varieties}

In this lecture, we discuss techniques to calculate the cohomology
of Shimura varieties, concentrating on the compact non-endoscopic
case (see~\ref{ssection_KC} for a short survey of the general case).

We fix a connected reductive group $G$ over $\Q$; let $Z(G)$ be the center
of $G$, and let $S_G$ be the maximal $\Q$-split torus in $Z(G)$.
Except in Subsection~\ref{ssection_KC},
we will often assume that $G_\der$ is of $\Q$-rank $0$, which means that $G$
has no proper parabolic subgroup.

\subsection{Matsushima's formula}

\subsubsection{Discrete automorphic representations}
\label{sssection_dar}

Let $A_G=S_G(\R)^0$.
If $G_\der$ is of $\Q$-rank $0$, then
the quotient $G(\Q)\sous G(\Ade)/A_G$ is compact; this follows from
the theorem on page 461 of the paper \cite{MosTam} by Mostow and
Tamagawa and from the fact that, with the notation of that paper, we have
$G(\Ade)=A_G\cdot G(\Ade)^1$ (which is easy to prove). 
Without this assumption, the quotient $G(\Q)\sous G(\Ade)/A_G$ is still of
finite volume.
\footnote{This is due to Borel in this generality once we know about the
relation to locally symmetric spaces (cf~\ref{subsection_lss});
see~\cite{Borel-GA} 8.4 and 13.2.}
We denote by
$L^2_G$ the space of complex $L^2$ functions on 
$G(\Q)\sous G(\Ade)/A_G$. The group $G(\Ade)$ acts on $L^2_G$ by
right translation on the argument, and this defines a continuous unitary
representation of $G(\Ade)$.

If $G(\Q)\sous G(\Ade)/A_G$ is compact, the representation $L^2_G$
decomposes as a Hilbertian sum of irreducible
representations $\pi$ of $G(\Ade)$, with finite multiplicities
$m(\pi)$.
\footnote{This is because nice enough functions on $G(\Ade)$ act on
$L^2_G$ by compact operators, see for example Proposition~2.2.1
and Corollary~2.2.2 of Booher's notes~\cite{Boo}.}
This is not true in general, but we can still consider the
part $L^2_{G,\disc}$ of $L^2_G$ that decomposes discretely, and we still
denote by $m(\pi)$ the multiplicity of an irreducible representation
$\pi$ of $G(\Ade)$ in $L^2_{G,\disc}$.

We denote by $\Pi_G$ the set of equivalences classes of irreducible
representations $\pi$ of $G(\Ade)$ such that $m(\pi)\not=0$.
Elements of $\Pi_G$ are called \emph{discrete automorphic representations}
of $G(\Ade)$ (or just of $G$).
If $G_\der$ is of rank $0$, then the concepts of
automorphic representations, discrete
automorphic representations and cuspidal automorphic representations
coincide.

If $\pi\in\Pi_G$, then, as $G(\Ade)=G(\Af)\times G(\R)$, 
we can write $\pi=\pi_f\otimes\pi_\infty$, where
$\pi_f$ (resp. $\pi_\infty$) is an irreducible representation of
$G(\Af)$ (resp. $G(\R)$).

\subsubsection{The theorem}

To state Matsushima's formula, we first need some definitions.
Let $\ggoth$ be the Lie algebra of $G_\R$, let $K'_\infty$ be a maximal
compact subgroup of $G(\R)$, and set $K_\infty=A_G\cdot K'_\infty$.
If $\pi\in\Pi_G$, then
we denote by $\H^*(\ggoth,K_\infty;\pi_\infty)$ 
the \emph{$(\ggoth,K_\infty)$-cohomology} of $\pi_\infty$, i.e. the
cohomology of the complex $C^q(\ggoth,K_\infty;\pi_\infty)=\Hom_{K_\infty}(
\wedge^q(\ggoth/\kgoth),(\pi_\infty)^\infty)$, where $\kgoth=\Lie(K_\infty)$
and $(\pi_\infty)^\infty$ is the space of smooth vectors in $\pi_\infty$ 
(which is stable by $K_\infty$ because $\pi_\infty$ is
$K'_\infty$-finite and $A_G$ acts trivially).
\footnote{See Section~I.5 of Borel and Wallach's book \cite{BoWa}
for more about $(\ggoth,K)$-cohomology. The main point for us is that
this is something that can in theory be calculated.}

The following theorem for connected components of complex Shimura varieties
is Corollary~VII.3.4 of \cite{BoWa}, and the adelic reformulation can be
found in Section~2 of Arthur's paper \cite{A-L2}.

\begin{theorem}[Matsushima's formula]\label{thm_Matsushima}
Let $(G,h)$ be a Shimura datum with $G_\der$ of $\Q$-rank $0$. Then we
have a $G(\Af)$-equivariant isomorphism of graded $\C$-vector spaces
\[\varinjlim_K\H^*(M_K(G,h)(\C),\C)\simeq\bigoplus_{\pi\in\Pi_G}
\pi_f\otimes\H^*(\ggoth,K_\infty;\pi_\infty)^{m(\pi)},\]
where $G(\Af)$ acts on the factors $\pi_f$ on the right hand side.

\end{theorem}

Here $\H^*(M_K(G,h)(\C),\C)$ is Betti cohomology with coefficients in $\C$.

The theorem is equivalent to the following corollary: Let $K$ be an open
compact subgroup of $G(\Af)$, and
remember from
page~\pageref{def_Hecke_correspondences} of Lecture~1 that the
\emph{Hecke algebra $\Hf_{G,K}$ at level
$K$} is the space of bi-$K$-invariant functions
from $G(\Af)$ to $\Q$ with compact support, with the convolution product
as multiplication.\footnote{Here we fixed any Haar measure on $G(\Af)$ such
that open compact subgroups of $G(\Af)$ have rational volume.}
If $\pi\in\Pi_G$, then $\pi_f^K$ is a finite-dimensional representation
of $\Hf_{G,K}\otimes_\Q\C$.

\begin{cor}
Let $(G,h)$ be as in Theorem~\ref{thm_Matsushima}, and let $K$ be an open
compact subgroup of $G(\Af)$. Then we have an isomorphism of
graded $\Hf_{G,K}\otimes_\Q\C$-modules
\[\H^*(M_K(G,h)(\C),\C)\simeq\bigoplus_{\pi\in\Pi_G}
\pi_f^K\otimes\H^*(\ggoth,K_\infty;\pi_\infty)^{m(\pi)}.\]

\end{cor}

\subsection{\'Etale cohomology of canonical models: the Kottwitz conjecture}

Let $\ell$ be a prime number.

If $(G,h)$ is a Shimura datum, then the projective system $(M_K(G,h)(\C))_K$
with its action of $G(\Af)$ has
a model over the reflex field $F=F(G,h)$.
So the $\ell$-adic \'etale
cohomology $\H^*_\et(M_K(G,h)_{\overline{F}},\oQ_\ell)$ has commuting actions
of $\Hf_{G,K}$ and $\Gal(\overline{F}/F)$.
For every isomorphism $\oQ_\ell\simeq\C$, we have comparison isomorphisms
\[\H^*_\et(M_K(G,h)_{\overline{F}},\oQ_\ell)\simeq\H^*(M_K(G,h)(\C),\C)\]
equivariant for the action of $\Hf_{G,K}$, because this action comes from
the geometric Hecke correspondences
(see Subsection~\ref{subsection_lss}).
So, when $G_\der$ is of $\Q$-rank $0$,
Matsushima's formula tells us that the action of $\Hf_{G,K}\otimes_\Q
\oQ_\ell$ on the cohomology groups $\H^i_\et(M_K(G,h)_{\overline{F}},\oQ_\ell)$
is semi-simple, that the only representations of $\Hf_{G,K}$ that appear are
the $\pi_f^K$ for $\pi\in\Pi_G$, and that the corresponding $\pi_f^K$-multiplicity space
$\H^i_K(\pi_f):=\Hom_{\Hf_{G,K}}(\pi_f^K,\H^i_\et(M_K(G,h)_{\overline{F}},\oQ_\ell))$
is of dimension
\[\sum_{\pi'\in\Pi_G\mid\pi'_f\simeq\pi_f}m(\pi')\dim\H^i(\ggoth,K_\infty;\pi'_\infty).
\]

We would like to calculate the action of $\Gal(\overline{K}/K)$ on
$\H^i_K(\pi_f)$. 
Kottwitz has a very precise conjecture describing
$\H^i_K(\pi_f)$, and we want to state that conjecture in the simplest case.
We need some preparation.

\vspace{.3cm}

\subsubsection{The Langlands group of $F$}
\label{ssection_Langlands_group}

The Langlands group $\Lf_F$ of $F$ is a conjectural
group scheme over $\C$ whose
irreducible representations on $n$-dimensional vector spaces should
classify the cuspidal automorphic representations of $\GL_n(\Ade_F)$.
Remember that we defined discrete automorphic representations of
$\GL_n(\Ade_F)=\GL_n(F\otimes_\Q\Ade)$ in~\ref{sssection_dar};
roughly speaking, cuspidal automorphic representations are discrete
automorphic representations that don't arise from an automorphic representation
of a Levi subgroup of $\GL_n$ via parabolic induction. Langlands
conjectures that there is a bijection $\pi\mapsto\phi_\pi$ from the
set of cuspidal automorphic representations $\pi$ of $\GL_n(\Ade_F)$ to
the set of equivalence classes of representations $\phi_\pi:\Lf_F\ra
\GL_n(\C)$, and that this correspondence is determined by local
compatibilities. We call $\phi_\pi$ the \emph{Langlands parameter} of $\pi$.

More precisely, if $\pi$ is a discrete automorphic
representation of $\GL_n(\Ade_F)$, then we can write
$\pi$ as a restricted tensor product $\bigotimes_v'\pi_v$ over all places
of $F$, where $\pi_v$ is an irreducible admissible representation of
$\GL_n(F_v)$. On the other hand, for each place $v$ of $F$, we have the
(non conjectural) Langlands groups $\Lf_{F_v}$ of $F_v$, with a
(conjectural) embedding $\Lf_{F_v}\subset\Lf_F$, and the
(non conjectural) local Langlands correspondence relates irreducible
admissible representations of $\GL_n(F_v)$ and $n$-dimensional representations
of $\Lf_{F_v}$ (``local Langlands parameters'').
The local Langlands correspondence over $\R$ and $\C$ was proved by Langlands
in~\cite{Lan-LL} (for all groups, not just $\GL_n$).
The local Langlands correspondence over $p$-adic local fields was proved independently
by Harris-Taylor (\cite{HaTa}), Henniart (\cite{Henniart}) and Scholze
(\cite{Sc-LL}).
For our purposes it is enough to understand the unramified local
Langlands correspondence (over $p$-adic fields),
which is just given by the Satake isomorphism;
see Dick Gross's notes~\cite{Gross-Satake} for an introduction
to the Satake isomorphism.
If $\pi$ is a cuspidal automorphic representation and we write
$\pi=\bigotimes_v'\pi_v$, then we expect that, for every place $v$ of
$F$, the restriction $\phi_{\pi\mid\Lf_{F_v}}$ corresponds to $\pi_v$ by the
local Langlands correspondence, and that this uniquely determines $\phi_\pi$.
In fact, it should be enough to know $\phi_\pi$ at the finite places $v$ such
that $\pi_v$ is unramified (hence the corresponding representation of
$\Lf_{F_v}$ is given by the Satake isomorphism).

It is also expected that $\Lf_F$ canonically surjects to $\Gf_{F,\C}$,
where $\Gf_F$ is the motivic Galois group of $F$, a group scheme over
$\Q$ defined as the Tannakian group of the conjectural category of
mixed motives over $F$.
\footnote{To make this precise, we need a fiber functor. We fix an
embedding of $F$ into $\C$ and take the fiber functor given by the
corresponding Betti realization with $\Q$-coefficients.}
The irreducible representations of $\Lf_F$ factoring through
$\Gf_{F,\C}$ are supposed to correspond to automorphic representations
satisfying a certain condition at the infinite places of $F$,
called \emph{algebraic automorphic representations}. 
On the other hand, for every prime number $\ell$, the \'etale
$\ell$-adic realization functor defines a continuous morphism of groups
$\Gal(\overline{F}/F)\ra\Gf_F(\oQ_\ell)$, that is supposed
to be injective (by the conservativity conjecture) and have dense
image (by the Tate conjecture). 

\begin{remark}
The philosophy behind the Langlands group of a number field, and
its relation to the motivic Galois group, are explained much better
in Clozel's paper \cite{Clo-MFA}. The extension of Clozel's ideas to
groups other than $\GL_n$ was worked out in the paper~\cite{Buzzard-Gee}
by Buzzard and Gee.

\end{remark}

Now we come back to the case of a connected reductive group $G$.
We need to define the $L$-group of $G$; a good reference for this
is Part~I of~\cite{Bo-ALF} or Section~1 of \cite{K-CTT}.
Let $T$ a be maximal torus of $G_{\oQ}$ and $B$ be a Borel subgroup
of $G_{\oQ}$ containing $T$.
The \emph{based root datum} of $G$ is the family
$\mathfrak{R}(B,T):=(X^*,\Phi,\Delta,X_*,\Phi^\vee,\Delta^\vee)$,
where $X^*=X^*(T_{\oQ})$, $\Phi\subset X^*$ is the set of roots of
$T_{\oQ}$ in $G_{\oQ}$, $\Delta$ is the set of simple roots corresponding to $B$,
$X_*=X_*(T_{\oQ})$ and $\Delta^\vee\subset\Phi^\vee\subset X_*$ are the set of simple coroots
and coroots. The based root datum is independent of the choice of
$(T,B)$ in the following sense: if $T'$ is a maximal torus of $G_{\oQ}$ and
$B'\supset T'$ is a Borel subgroup of $G_{\oQ}$, then there exists
$g\in G(\oQ)$ such that $B'=gTg^{-1}$ and $T'=gTg^{-1}$, and the isomorphism
$\mathfrak{R}(B,T)\iso\mathfrak{R}(B',T')$ induced by $g$ is independent of the
choice of $g$. So we just write $\mathfrak{R}(G)$ instead of $\mathfrak{R}(B,T)$.
In particular, we get a morphism of groups from
$\Aut(G_{\oQ})$ to the group of automorphisms of the based root datum of $G$, which turns
out to be surjective with kernel the group of inner automorphisms; any
pinning of $G$ gives a splitting of the surjective morphism $\Aut(G_{\oQ})\ra\Aut(\mathfrak{R}(G))$.

The \emph{dual group} $\widehat{G}$ of $G$ is the complex
connected reductive group with based root datum $(X_*,\Phi^\vee,\Delta^\vee,X^*,\Phi,\Delta)$.

Until now, we just needed $G$ to be a connected reductive group over $\oQ$. But
as $G$ is defined over $\Q$, we get an action of $\Gal(\oQ/\Q)$ on its based
root datum $\mathfrak{R}(G)$, and
if we fix a pinning of $\widehat{G}$, then this defines an action of
$\Gal(\oQ/\Q)$ on $\widehat{G}$. The \emph{$L$-group of $G$}
is ${}^L G=\widehat{G}\rtimes W_\Q$, where $W_\Q\subset\Gal(\oQ/\Q)$ is
the Weil group of $\Q$.
\footnote{For the Weil group, see for example Tate's Corvallis notes~
\cite{Tate}.}

\begin{remark}
If $\ell$ is a prime number, we could define $\widehat{G}$ to be the
connected reductive group
over $\oQ_\ell$ with root datum
$(X_*,\Phi^\vee,X^*,\Phi)$, and we would get a group ${}^LG$ over
$\oQ_\ell$. We write $\widehat{G}(\oQ_\ell)$ and ${}^LG(\oQ_\ell)$ 
for the resulting groups, when we want to distinguish them from the
complex versions.
Which form of the $L$-group we use depends on the context:
for Langlands parameters defined on $\Lf_F$, we use the complex form, and
for Langlands parameters defined on $\Gal(\overline{F}/F)$, we use the
$\ell$-adic form. 

\end{remark}

\begin{example}
If $G$ is of type $A$, $D$, $E$, $F$ or $G$, then $\widehat{G}$ is of the
same type as $G$. If $G$ is of type $B_n$ (resp. $C_n$), then $\widehat{G}$ is
of type $C_n$ (resp. $B_n$). We can also relate other properties of
$G$ and $\widehat{G}$: for example, the derived group $G_\der$ is simply
connected if and only if $Z(\widehat{G})$ is connected, and in that
case $Z(\widehat{G})$ is the dual group of $G/G_\der$. 

Here are some examples of dual groups:
\begin{itemize}
\item $\widehat{\GL_n}=\widehat{\U(p,q)}=\GL_n(\C)$ if $p+q=n$;
\item $\widehat{\SL_n}=\widehat{\SU(p,q)}=\PGL_n(\C)$ if $p+q=n$;
\item $\widehat{\PGL_n}=\SL_n(\C)$;
\item $\widehat{\Sp_{2n}}=\SO(2n+1,\C)$, $\widehat{\GSp_{2n}}=\GSpin_{2n+1}(\C)$;
\item $\widehat{\GU(p,q)}=\GL_1(\C)\times\GL_n(\C)$.

\end{itemize}

We see that different groups can have isomorphic dual groups. In fact, if
$G'$ is an inner form of $G$, then the $L$-groups ${}^LG$ and ${}^LG'$ are
isomorphic. However, as $\U(p,q)$ is not an inner form of $\GL_n$, the
actions of $\Gal(\oQ/\Q)$ on their dual groups are not the same, and we
get non-isomorphic $L$-groups. If $G= \GL_n$ (or more generally if
$G$ is split over $\Q$), then $\Gal(\oQ/\Q)$ acts trivially on
$\widehat{G}$. On the other hand, if $G=\U(p,q)$ and if $E$ is the
imaginary quadratic extension of $\Q$ that we used to define $G$, then
$G$ splits over $E$, so
$\Gal(\oQ/\Q)$ acts on $\widehat{G}=\GL_{p+q}(\C)$ via its quotient
$\Gal(E/\Q)$, and the nontrivial element of $\Gal(E/\Q)$ acts
as a non-inner automorphism of $\GL_{p+q}(\C)$, i.e. a conjugate of the
automorphism $g\mapsto{}^tg^{-1}$ (which conjugate depends on the choice of
the pinning).

\end{example}

Coming back to the Langlands correspondence for $G$, there are several complications:
\begin{itemize}
\item A cuspidal automorphic representation $\pi$ of $G(\Ade)$ should now
have a Langlands parameter $\phi_\pi$ with values not in $\GL_n(\C)$ but in
${}^L G$;

\item There is a still a characterization of algebraic automorphic
representations (conjecturally corresponding to the parameters that
factors through $\Gf_{F,\C}$), but it is more complicated, see
Buzzard and Gee's paper \cite{Buzzard-Gee};

\item Distinct cuspidal automorphic representations can have the same
Langlands parameter. We say that they are in the same $L$-packet.

\end{itemize}

If $\pi$ is an algebraic cuspidal automorphic representation of $G$,
then we expect $\phi_\pi:\Gf_{\Q,\C}\ra {}^L G$ to be defined over a finite
extension $L$ of $\Q$ in $\C$. Choosing a finite place $\lambda$ of
$L$ over a prime number $\ell$, we get a morphism from
$\Gf_\Q(\oQ_\ell)$ into the $\ell$-adic version of ${}^L G$, and this
gives a morphism $\sigma_\pi:\Gal(\oQ/\Q)\ra{}^L G(\oQ_\ell)$, also
called the Langlands parameter of $\pi$, and
whose value on
the Frobenius elements at almost all prime numbers $p$ is predicted by
the Satake parameter of $\pi_p$. Now the conjecture only involves
well-defined objects, and we can actually try to prove it !

\begin{remark}\label{rmk_Arthur_param}
If we are very brave and want to classify all discrete automorphic
representations of $G(\Ade)$, then there is an extension of the
Langlands conjecture due to Arthur. Now a discrete automorphic representation
$\pi$ should have a parameter $\psi_\pi:\Lf_F\times\SL_2(\C)\ra{}^L G$,
satisfying a long list of properties (in particular compatibility with
a local version of the Arthur conjectures), and there is a somewhat
explicit formula to calculate the multiplicity $m(\pi)$. For a quick
review of Arthur's conjectures, see Section~8 of Kottwitz's paper
\cite{K-SV}. 
Warning: if $\pi$ is cuspidal, then we can recover
$\phi_\pi$ from $\psi_\pi$ and vice versa, but the two parameters are
not equal.

\end{remark}

\subsubsection{The Kottwitz conjecture}
\label{sssection_KC}

Let $(G,h)$ be a Shimura datum, and let $F=F(G,h)$ be its reflex
field. We assume that $G_\der$ is of $\Q$-rank $0$, so that the
Shimura varieties $M_K(G,h)(\C)$ are compact.

The conjugacy class of the cocharacter $\mu:=h_\C\circ r:\GL_{1,\C}
\ra G_\C$ is defined over $F$ (see~\ref{sssection_rf}), hence
it defines a finite-dimensional representation $r_\mu$ of
${}^L G_F(\oQ_\ell):=\widehat{G}(\oQ_\ell)\rtimes W_F$ in the following way
(see Lemma~2.1.2 of Kottwitz's paper \cite{K-TOI}):
choose a maximal torus $\widehat{T}$ of $\widehat{G}$ and a Borel
subgroup $\widehat{B}$ containing $\widehat{T}$ that are part of a splitting
fixed by $W_F$. 
The cocharacter $h_\C\circ r$ corresponds to a unique dominant character
$\mu$ of $\widehat{T}$, an we denote by $V_\mu$ the corresponding
highest weight representation of $\widehat{G}(\oQ_\ell)$. The action of
$\widehat{G}(\oQ_\ell)$ on $V_\mu$ extends to a unique action $r_\mu$ of 
${}^LG_F(\oQ_\ell)$ such that
$W_F$ acts trivially on the highest weight subspace.

\begin{examples}
\begin{itemize}
\item[(1)] If $G=\GSp_{2d}$ and $(G,h)$ is the Shimura datum of
Example~\ref{ex_SD_GSp}, then $\widehat{G}=\GSpin_{2d+1}(\C)$ and
$r_\mu:\widehat{G}\ra\GL_{2^d}(\C)$ is the spin representation.

If $d=1$, then $\GSp_{2d}=\GL_2$ and $r_\mu$ is the standard representation
if $\widehat{\GL_2}=\GL_2(\C)$.

If $d=2$, then we have an exceptional isomorphism
$\GSpin_5(\C)\simeq\GSp_4(\C)$, and $r_\mu$ is isomorphic to the standard
representation of $\GSp_4(\C)$. 

\item[(2)] If $G=\GU(p,q)$ and $(G,h)$ is the Shimura datum of
Example~\ref{ex_GU}, then $\widehat{G}\simeq\GL_1(\C)\times\GL_{p+q}(\C)$ and
$r_\mu$ is, up to twists by characters, the $q$th exterior power of the
standard representation of $\GL_{p+q}(\C)$.

\end{itemize}
\end{examples}

Let $K$ be an open compact subgroup of $G(\Af)$, 
let $d=\dim M_K(G,h)$,
and let $\pi$ be a cuspidal automorphic
representation of $G(\Ade)$ such that the $\pi_f^K$-multiplicity space
$\H^i_K(\pi_f)$ of the $i$th cohomology group
$\H^i(M_K(G,h)(\C),\C)$ is nonzero for at least one
$i\in\Z$. Then $\pi$ should be algebraic, so its Langlands parameter
should give rise to a Galois representation $\sigma_\pi:\Gal(\oQ/\Q)\ra
{}^L G(\oQ_\ell)$ as above.

Here we give a simplified version of the conjecture (from Section~1
of Kottwitz's paper~\cite{K-SSV}); the general version of the conjecture,
also due to Kottwzitz, can be found in Section~10 of~\cite{K-SV}, and we
will discuss it in~\ref{ssection_KC}.

\begin{conjecture}\label{conj_Kottwitz}
Under some assumptions on the Shimura datum (see Remark~\ref{remark_KC}(1)),
there is an explicitly defined integer $a(\pi_f)$ such that
we have an equality of virtual
representations of $\Gal(\overline{F}/F)$:
\[\sum_{i=0}^{2d}(-1)^i[\H^i_K(\pi_f)(d/2)]=a(\pi_f)[(r_\mu\circ\sigma_\pi)].\]
Moreover, the integers $i$ such that $\H^i_K(\pi_f^K)\not=0$ all have the
same parity.

\end{conjecture}

\begin{remark}\label{remark_KC}
\begin{itemize}
\item[(1)]
The conjecture as stated is in a very naive form and false for general
Shimura varieties.
In fact, we only expect it to be true when $G_\der$ has $\Q$-rank $0$
and
$G$ has no endoscopy, i.e. satisfies Proposition~\ref{prop_no_end}
(the correct technical condition is that $G$ should admit no
nontrivial elliptic endoscopic triple). 

As the Shimura varieties are then projective (and smooth for $K$ small enough), 
by the Weil conjectures (proved by Deligne, see~\cite{De-WeilI}), 
the representation
$\H^i_K(\pi_f)$ of $\Gal(\overline{F}/F)$ is pure of weight $i$ for every
$i\in\Z$, so we can separate the degrees in the formula of
Conjecture~\ref{conj_Kottwitz} by using Frobenius weights.

\item[(2)] While we cannot hope to prove the Kottwitz conjecture
without first constructing the Langlands parameter of $\sigma_\pi$, we do
know what the image by $\sigma_\pi$ of the Frobenius element $\Frob_\wp$ at a
place $\wp$ of $F$ over almost every
prime number $p$ should be, so we can try to prove that
$\Frob_\wp$ has the correct characteristic polynomial on the
$\H^i_K(\pi_f)$.
For the simple Shimura
varieties of Subsection~\ref{ssection_SSV}, Kottwitz proved this consequence
of his conjecture in \cite[Theorem~1]{K-SSV}.

\item[(3)] If we know the local Langlands correspondence for $G$, we
can also try to check that the restriction to all local Galois
groups of the representations
$\H^i_K(\pi_f)$ are as predicted by the Kottwitz conjecture.
This is a much harder problem and we won't discuss it here; see for example
the book \cite{HaTa} of Harris-Taylor for an example among many of this
kind of calculation.

\end{itemize}
\end{remark}

\subsection{Applications of the Kottwitz conjecture}

\subsubsection{The zeta function of a Shimura variety}
\label{sssection_zeta}

\begin{definition}
Let $X$ be a smooth proper variety over a finite field $\Fi_q$.
The \emph{Hasse-Weil zeta function} of $X$ is the following
formal power series in $q^{-s}$:
\[Z(X,s)=\exp(\sum_{n\geq 1}\frac{\card(X(\Fi_{q^n}))}{n}q^{-sn}).\]

\end{definition}

Using the Grothendieck-Lefschetz fixed point formula (cf.
Theorem~\ref{thm_GL}), we get the following result.

\begin{theorem}[Grothendieck, see Theorem 3.1 of~\cite{RFT}]
\label{thm_ZF}
Let $\Frob_q\in\Gal(\overline{\Fi}_q/\Fi_q)$ be the geometric Frobenius (the
inverse of the arithmetic Frobenius $a\mapsto a^q$). Then
\[Z(X,s)=\prod_{i=0}^{2\dim(X)}\det(1-q^{-s}\Frob_q,\H^i_\et(X_{\overline{\Fi}_q},
\Q_\ell))^{(-1)^{i+1}}.\]

\end{theorem}

In particular, the formal power series $Z(X,s)$ is actually a rational
function in $q^{-s}$. (This last corollary was already known by work
of Dwork, see~\cite{Dwork}.)

\begin{remark}
We can define the zeta function of any algebraic variety over
$\Fi_q$ by the same formula, and Theorem~\ref{thm_ZF} still holds providing
we use \'etale cohomology with proper supports.

\end{remark}

Now let $X$ be a proper smooth algebraic variety over a number field $F$.
For all but finite many finite places $\wp$ of $F$, the variety $X$ has
a proper smooth model $\Xf$ over $\Of_{F,\wp}$ (we say that $X$ has
\emph{good reduction} at $\wp$), and we set
\[\zeta_{X,\wp}(s)=\zeta_{\Xf_{\kappa(\wp)}}(\card(\kappa(\wp))^{-s}),\]
where $\kappa(\wp)=\Of_F/\wp$ is the residue field of $\wp$.
By Theorem~\ref{thm_ZF} and the specialization theorem for \'etale
cohomology,
\footnote{Which follows
from the proper and smooth base change theorems, see
\cite{SGA43} Expos\'es XII and XVI.
}
this does not depend on the choice of the model.

If $v$ is a finite place of $F$ where $X$ does not have good reduction or
an infinite place, we will not give the definition of $\zeta_{X,v}(s)$;
we will just say that $\zeta_{X,v}(s)$ is a rational function of
$\card(\kappa(v))^{-s}$ if $v$ is finite and a product of $\Gamma$ functions
if $v$ is infinite. See Serre's notes~\cite{Serre-zeta} for the definitions
and for the conjectures on Hasse-Weil zeta functions.

\begin{definition}
The \emph{Hasse-Weil zeta function} of $X$ is the infinite product
\[\zeta_X(s)=\prod_{v\ \mathrm{place\ of}\ F}\zeta_{X,v}(s).\]

\end{definition}

\begin{example}
If $X=\Spec\Q$, then $\zeta_X$ is the Riemann zeta function.

\end{example}

A priori this product only makes sense for $\Re(s)$ big enough. The
\emph{Hasse-Weil conjecture} predicts that $\zeta_X(s)$ has a meromorphic
continuation to $\C$ and a functional equation similar to the one
of the Riemann zeta function. 

This conjecture seems to be out of reach
in general, but for Shimura varieties we can approach it using
the Kottwitz conjecture. The general idea goes as follows:
\begin{itemize}
\item[(1)] Essentially by Theorem~\ref{thm_ZF}, we have an equality
\[\zeta_X(s)=\prod_{i=0}^{2\dim(X)}L(\H^i_\et(X_{\overline{F}},\Q_\ell),s)^{(-1)^{i+1}},\]
where, for every continuous representation $\rho$ of $\Gal(\overline{F}/F)$,
we denote by $L(\rho,s)$ the \emph{$L$-function} of $\rho$.
(See Section~3 of Tate's Corvallis notes~\cite{Tate}.)

\item[(2)] Suppose that $(G,h)$ is a Shimura datum such that
$G_\der$ has $\Q$-rank $0$, that $F=F(G,h)$ and
$X=M_K(G,h)$.
For every $i$, we have up to semi-simplification
\[\H^i_\et(M_K(G,h)_{\overline{F}},\oQ_\ell)\simeq\bigoplus_{\pi\in\Pi_G}
\H^i_K(\pi_f)^{\dim(\pi_f^K)}\]
as representations of $\Gal(\overline{F}/F)$.

\item[(3)] The Kottwitz conjecture predicts that $\H^i_K(\pi_f)$ is a sum
of copies of $r_\mu\circ\sigma_\pi$, where $\sigma_\pi:\Gal(\overline{F}/F)\ra
{}^LG(\oQ_\ell)$ is the Langlands parameter of $\pi$ and $r_\mu$ is the
algebraic representation of ${}^LG$ defined in~\ref{sssection_KC}.
But the local compatibility between $\pi$ and $\sigma_\pi$ implies
immediately that
\[L(r_\mu\circ\sigma_\pi,s)=L(\pi,s,r_\mu),\]
where the $L$-function $L(\pi,s,r_\mu)$ is defined in
Borel's survey \cite{Bo-ALF}.

\item[(4)] In theory we understand the analytic properties of
$L$-functions of automorphic representations better, so we get
some information on the zeta function of $M_K(G,h)$.

\end{itemize}

In practice the automorphic $L$-functions that appear are usually
not standard $L$-functions and so our understanding of them is still
limited. However, this method can still go through when $r_\mu$ is
the standard representation of a classical group, such as in the case
of modular curves (when $G=\GL_2$) or Picard modular surfaces (when
$G=\GU(2,1)$), or more generally when $G_\R$ is a unitary group of
signature $(n-1,1)\times(n,0)^{r-1}$, as can be the case for
Kottwitz's simple Shimura varieties (see~\ref{subsubsection_SSV}).

\subsubsection{The global Langlands correspondence}
\label{sssection_Langlands_corr}

Let $(G,h)$ be a Shimura datum, and let $F=F(G,h)$.
In a way, the Kottwitz conjecture
says that the cohomology $\varinjlim_K\H^*_\et(M_K(G,h)_{\overline{F}},\oQ_\ell)$
realizes the global Langlands correspondence for those automorphic
representations of $G(\Ade)$ that contribute to Matsushima's formula for
$M_K(G,h)(\C)$.
So we could try to use this cohomology to construct the global Langlands
correspondence in that case, and then use results like the main theorem
of Kottwitz's \cite{K-SSV} to check that this does satisfy the desired
compatibility with the local correspondence (at least in the unramified
case). 

For $\GL_2$ and its inner forms, this is older than Kottwitz's conjecture;
see for example the papers~\cite{Ku-Shi} and~\cite{Deligne-FML}
of Kuga-Shimura and Deligne, respectively. For $\GL_n$ with
$n\geq 3$, we have no Shimura variety, so we have to use other groups,
and we run into several problems:
\begin{itemize}
\item[(1)] The representation of $\Gal(\overline{F}/F)$
that appears in the cohomology of
$M_K(G,h)$ is not $\sigma_\pi$ (this would not even make sense, as
$\sigma_\pi$ is a morphism into ${}^LG(\oQ_\ell)$) but
$r_\mu\circ\sigma_\pi$;

\item[(2)] There are multiplicities (the integer $a(\pi_f)$ in the
Kottwitz conjecture);

\item[(3)] We want to construct the Langlands correspondence for
$\GL_n$, not some strange unitary or symplectic group;

\item[(4)] We want to get a representation of $\Gal(\oQ/\Q)$, not
$\Gal(\overline{F}/F)$;

\item[(5)] This will only ever work for ``cohomological''
automorphic representations, i.e. those $\pi$ that
appear in the cohomology of Shimura varieties. This is a condition on
$\pi_\infty$: roughly, we need it to have nontrivial 
$(\ggoth,K_\infty)$-cohomology.

\end{itemize}

All of these can be somewhat addressed, at some cost. For point (1), we
can choose the group such that $\widehat{G}$ is classical and $r_\mu$
is the standard representation. This will for example be the case if
$G=\GU(n-1,1)$, although in practice we will rather want $G$ to be a more
complicated unitary group (defined by a CM extension of $\Q$ of
degree $2r>2$) of signature $(n-1,1)\times(n,0)^{r-1}$ at infinity, so
that we can get a simple Shimura variety. Of course, this solution
puts even greater restrictions on the groups that we can use, so it
seems that we are making problems (3) and (4) worse. For problem (2), we
can sometimes calculate the multiplicities, and they tend to be equal to $1$
for nice unitary groups.

Problem (3) can be attacked using the \emph{Langlands functoriality principle}.
The idea is that, if discrete automorphic representations of $G(\Ade)$
are parametrized by morphisms $\Lf_\Q\times\SL_2(\C)\ra{}^L G$, then, if
$H$ is another connected reductive algebraic group over $\Q$ and if
we have a morphism ${}^LG\ra{}^L H$, then we should be able to
``transfer'' discrete automorphic representations from $G$ to $H$.
As with the global Langlands correspondence or the Arthur conjectures, 
this principle can be made
very precise at ``good'' primes (i.e. primes where both groups and all
automorphic representations we consider are unramified) using the Satake
isomorphism, so we can pin down the conjectural transfer using a local-global
compatibility principle. Of course, things are not so simple: the conjectural
Arthur parametrization is not bijective in general so we can only
expect to transfer $L$-packets, and the heuristic does not tell us how to
actually construct the transfer.

One favorable case is when $\widehat{G}$ is the set of fixed points
of an automorphism of $\widehat{H}$, because then
we can use the (twisted) Arthur-Selberg trace formula to construct
the transfer map, although this is very technically difficult;
see~\ref{sssection_Arthur_TF} for a very simple instance
of the (untwisted) trace formula.
This is for example the case if $G=\Sp_{2d}$ and $H=\GL_{2d+1}$
(see Arthur's book~\cite{Arthur-fonc}),
if $G$ is the group $\U(p,q)$ constructed using a quadratic imaginary
extension $E$ of $\Q$ and $H=\GL_{p+q}(E)$, seen as an algebraic
group over $\Q$
(see the works by Mok \cite{Mok} and Kaletha-Minguez-Shin-White~\cite{KMSW}),
or if $G=\GL_n$ and $H=\GL_n(E)$ (seen as an algebraic group over $\Q$)
for $E/\Q$ a cyclic extension
(see the book \cite{ArClo} of Arthur-Clozel).
Using trace formula techniques, we can transfer discrete automorphic
representations of $\GL_{p+q}(\Ade)$ to $\GL_{p+q}(\Ade_E)$, then
back down to $\U(p,q)(\Ade)$, and then use the cohomology of Shimura
varieties for this unitary group to construct Galois representations.
One caveat is that this only works for
representations of $\GL_{p+q}(\Ade_E)$ that are conjugate self-dual
(because their parameter should be stable by the involution of
$\widehat{\GL_{p+q}(E)}$ whose set of fixed points is equal to
$\widehat{\U(p,q)}$).
This kind of construction will also work for more general unitary
groups defined using other CM fields or division algebras over them, and
in fact we want to use some of them rather than $\U(p,q)$ because of problem (1).
Again, this is difficult and there are many technical problems, including
annoying congruence conditions on $p+q$. The first results that we are
aware of on the construction of parts of the Langlands correspondence
using the cohomology of Shimura varieties, outside of the case
of $\GL_2$,
are due to Clozel (see \cite{Clo-RG}, which
rests on the results of Kottwitz from \cite{K-SSV}). More recent
and more powerful results can be found in the papers \cite{Shin} 
of Shin and \cite{ScShin} of Scholze-Shin; these feature
unitary Shimura varieties that are not simple, so they use a more
complicated form of the Kottwitz conjecture (see Subsection~\ref{ssection_KC}),
as well as the fundamental lemma.
See Dat's Bourbaki notes~\cite{Dat-LF} for a statement of the fundamental
lemma and some context; the fundamental lemma is known by work of
Laumon-Ngo (\cite{LauNgo}),
Ngo (\cite{Ngo})
and Waldspurger (\cite{Wald1}, \cite{Wald2}).

Problem (4) can be addressed by ``gluing'' Galois representations
constructed using different Shimura varieties. More precisely, suppose
that you have a discrete automorphic representation $\pi$
of $\GL_n(\Ade)$, and that you know that it is self-dual and that $\pi_\infty$
satisfies all the required conditions for $\pi$ to transfer to cohomological
automorphic representations of unitary groups. By varying the CM field $E_i$
(and the corresponding unitary group), we get a family of representations
$\sigma_i$ of $\Gal(\overline{E}_i/E_i)$, and $\sigma_i$ should be
the restriction of $\sigma_\pi$. We can check that these representations
are compatible by looking at what happens at ``nice'' prime numbers, and then
glue them if we chose the $E_i$ disjoint enough.

Problem (5) is already an issue for $\GL_2$ and requires new techniques,
as in the paper~\cite{Deligne-Serre} of Deligne and Serre.

An additional problem, mentioned in the discussion of (3), is
that we can only ever transfer self-dual automorphic representations of
$\GL_n$ to groups that have Shimura varieties, but we will not try
to explain the solution to this problem here. 
See for example the papers \cite{HC} of
Chenevier-Harris, \cite{HLTT} of Harris-Lan-Taylor-Thorne,
\cite{Sc-torsion} of Scholze, or \cite{Boxer} of Boxer.

\subsection{Proving the Kottwitz conjecture}
\label{proof_KC}

We will present the original approach to the Kottwitz conjecture,
due to Ihara, Langlands and Kottwitz. We will not talk about the
more refined approaches through Igusa varieties or through Scholze's
methods (see for example the book \cite{HaTa}, or the
papers \cite{Shin} and \cite{ScShin} already cited above), 
that also allow us to understand the cohomology at ramified primes.

The situation is the following: we have a Shimura datum $(G,h)$ such that
$G_\der$ is of $\Q$-rank $0$ (so that the Shimura varieties $M_K(G,h)$
are projective), and a discrete automorphic representation $\pi$ of
$G(\Ade)$. We fix a small enough open compact subgroup $K$ of $G(\Af)$, and
we write $F=F(G,h)$. 
We are trying to understand the representation of $\Gal(\overline{F}/F)$ on
the $\pi_f^K$-multiplicity spaces
$\H^i_K(\pi_f)$ in the cohomology groups
$\H^i_\et(M_K(G,h)_{\overline{F}},\oQ_\ell)$.

\subsubsection{The specialization theorem}

We restrict our attention to the action of the local Galois group at finite
places $\wp$ of $F$ that are ``nice enough'', i.e. such that
$M_K(G,h)$ has a proper smooth integral model $\Mcal_K$ over $\Of_{F,\wp}$.
Remember from~\ref{sssection_abelian_type}
that, if $(G,h)$ is of abelian type, then
we have control (in terms of $G$ and $K$) over the places $\wp$ where we have
a smooth integral model; when $(G,h)$ is of Hodge type and $G_\der$ has $\Q$-rank $0$, we also
know that this integral model is projective by the work of Madapusi
Pera (see \cite{Mada-Hodge}). 
However,
in general we only know that we have smooth projective integral models 
for all but finitely many $\wp$.

Fix $\wp$ as in the previous paragraph, and let $\kappa(\wp)$ be the
residue field of $\wp$. Then we have an exact sequence
\[1\ra I_\wp\ra\Gal(\overline{F}_\wp/F_\wp)\ra\Gal(\overline{\kappa(\wp)}/
\kappa(\wp))\ra 1,\]
where $I_\wp$ is the \emph{inertia group at $\wp$}. Moreover, as $\kappa(\wp)$
is a finite field, its absolute Galois group is topologically generated
by the geometric Frobenius $\Frob_\wp$, which is the inverse of the
arithmetic Frobenius $a\mapsto a^{\card(\kappa(\wp))}$.

The specialization theorem for \'etale cohomology 
(which follows
from the proper and smooth base change theorems, see
\cite{SGA43} Expos\'es XII and XVI) tells us that the
representations $\H^i_\et(M_K(F,h)_{\overline{F}_\wp},\oQ_\ell)$ of
$\Gal(\overline{F}_\wp/F_\wp)$ are \emph{unramified}, i.e. that
$I_\wp$ acts trivially on them, and that we have isomorphisms of
representations of $\Gal(\overline{\kappa(\wp)}/\kappa(\wp))$:
\begin{equation}\label{eq_specialization}
\H^i_\et(M_K(F,h)_{\overline{F}_\wp},\oQ_\ell)\simeq
\H^i_\et(\Mcal_{K,\overline{\kappa(\wp)}},\oQ_\ell).\end{equation}

If we have a Hecke correspondence defined by $g\in G(\Af)$ and
$K'\subset K\cap gKg^{-1}$, then these isomorphisms will be compatible
with the corresponding Hecke operator, provided that $M_{K'}(G,h)$ also
has a proper smooth model over $\Of_{F,\wp}$ and that the
Hecke correspondence extends to the models over $\Of_{F,\wp}$.

So we can now work over the finite field $\kappa(\wp)$.

\vspace{.3cm}

\subsubsection{The Grothendieck-Lefschetz fixed point formula and
Deligne's conjecture}

Suppose that we were only trying to understand the representation
of $\Gal(\overline{F}_\wp/F_\wp)$ on
$\H^*_\et(M_K(G,h)_{\overline{F}_\wp},\oQ_\ell)$ (and not on the multiplicity
spaces $\H^i_K(\pi_f^K)$), and that we only cared about semisimplications. Then, by
the isomorphism \eqref{eq_specialization}, it would suffice to calculate
the characteristic polynomial of $\Frob_\wp$ acting on
$\H^*_\et(\Mcal_{K,\overline{\kappa(\wp)}},\oQ_\ell)$. We can do this
thanks to the Grothendieck-Lefschetz fixed point formula (already
mentioned before Theorem~\ref{thm_ZF}):

\begin{theorem}[Grothendieck, cf. Th\'eorème 3.2 of \cite{RFT}]
\label{thm_GL}
Let $\Fi_q$ be a finite field, $\overline{\Fi}_q$ be an algebraic closure
of $\Fi_q$, 
$\ell$ be a prime number different from the characteristic of $\Fi_q$
and $\Frob_q$ be the geometric Frobenius automorphic of $\overline{\Fi}_q$.
Then, for every separated $\Fi_q$-scheme of finite
type $X$ and every positive integer $r$, we have
\[\sum_{i=0}^{2\dim(X)}(-1)^i\Tr(\Frob_q^r,\H^i_{\et,c}(X_{\overline{\Fi}_q},\Q_\ell))=
\card(X(\Fi_{q^r})),\]
where $\Fi_{q^r}$ is the unique extension of $\Fi_q$ of degree $r$ in
$\overline{\Fi}_q$.

\end{theorem}

If we want instead to 
understand the action of $\Gal(\overline{F}_\wp/F_\wp)$ on the Hecke
isotypic components, we need to calculate the traces of Hecke operators
multiplied by powers of $\Frob_\wp$ on $\H^*_\et(M_K(G,h)_{\overline{F}_\wp},
\oQ_\ell)$.
For this, we use
a generalization of the Grothendieck-Lefschetz
fixed point formula. We will not state the
most general version here, but just the consequence that we need.

We use the notation of Theorem~\ref{thm_GL}. Let $X,X'$ be separated
$\Fi_q$-schemes of finite type, and let $a,b:X'\ra X$ be finite morphisms. 
Suppose that the
trace morphism $\Tr_b:b_*b^*\ra\id$ exists.
\footnote{This holds for example if $b$ is flat or if $X$, $X'$ are
normal; this is due to Deligne, see
SGA~4 XVII~6.2 in \cite{SGA43}.}
Let $u$ be the endomorphism of $\H^*_{\et,c}(X_{\overline{\Fi_q}},\Q_\ell)$ that is
the composition of the pullback by $a$ map $\H^*_{\et,c}(X_{\overline{\Fi_q}},\Q_\ell)
\ra\H^*_{\et,c}(X'_{\overline{\Fi_q}},\Q_\ell)$ and of the map
$\H^*_{\et,c}(X'_{\overline{\Fi_q}},\Q_\ell)\ra\H^*_{\et,c}(X_{\overline{\Fi_q}},\Q_\ell)$ 
induced by $\Tr_b$.
Finally, we denote by $F_X:X\ra X$ the Frobenius morphism (which is
identity on the underlying topological spaces and raises functions to the
$q$th power).

\begin{theorem}\label{thm_Deligne_conj}
For any big enough positive integer $r$, we have
\[\sum_{i=0}^{2\dim(X)}(-1)^i\Tr(\Frob_q^r\cdot u,\H^i_{\et,c}(X_{\overline{\Fi}_q},\Q_\ell))=
\card(\{x'\in X'(\overline{\Fi}_q)\mid a(x')=F_X^r\circ b(x')\}).\]

\end{theorem}

Although it is still often refereed to as ``Deligne's conjecture'',
this statement is a theorem: if $X,X'$ are reductions of Shimura varieties
and $a,b$ are Hecke operators, it was proved by Pink in \cite{Pink-DC}.
In the general case, it was proved independently by Fujiwara (\cite{Fu})
and Varshavsky (\cite{Var}).

So we now need to understand the set of points of the Shimura varieties
$M_K(G,h)$ and $M_{K'}(G,h)$ (or rather of their integral models) over
$\overline{\kappa(\wp)}$, as well as the action of the Frobenius morphism and
of Hecke correspondences on these points.

\vspace{.3cm}

\subsubsection{The Langlands-Rapoport conjecture}

The Langlands-Rapoport conjecture gives a purely group-theoretical
description of the set of points of a Shimura variety over the algebraic
closure of the residue field at a good place $\wp$ of the reflex field, as
well as a description of the action of the Frobenius at $\wp$ and of
Hecke operators on this set. We only give a rough
statement here. We actually present a corollary of the original
conjecture, that is, a description of the points of a Shimura variety over
a finite field (not its algebraic closure); the derivation is explained
in Section~5 of Milne's paper~\cite{Milne-LR}.
As a consequence, our
notation is not the same as in the paper of Langlands and Rapoport or the book of
Reimann (see the next paragraph), in
fact we follow Kottwitz's notation in~\cite{K-SV}. 
We will also restrict to the case of groups with simply connected
derived subgroups.

The original paper of Langlands and Rapoport on their conjecture is~\cite{LanRap};
another good source is Milne's survey paper~\cite{Milne-LR}. A more compact
presentation of the conjecture can be found in Appendix~B of Reimann's
book~\cite{Reim}, where the author also fixes some technical mistakes in
the previous references.

Let $(G,h)$ be a Shimura datum with reflex field $F$, let $p$ be a prime
number and $\wp$ a place of $F$ above $p$. Suppose that $G$ and $F$ are
unramified at $p$, and let $K=K^pK_p$ be a level with $K_p\subset G(\Q_p)$ 
hyperspecial.
Then we expect the Shimura variety $M_K(G,h)$ to have a ``nice''
model $\Mcal_K$ over $\Of_{F,\wp}$, and the Langlands-Rapoport conjecture
gives a description of $\Mcal_K(\Fi_q)$, where $\Fi_q$ is an extension
of $\Of/\wp$. The description
has the rough shape
\begin{equation}\label{LRconj}
\Mcal_K(\Fi_q)=\coprod_\varphi I_\varphi(\Q)\sous(X^p(\varphi)\times
X_p(\varphi)),
\end{equation}
where $\varphi$ is in a certain set of parameters, $I_\varphi$ is an
algebraic group over $\Q$ (the centralizer of $\varphi$),
$X^p(\varphi)$ involves the finite adeles outside of $p$ and
$X_p(\varphi)$ is a purely $p$-adic object. The conjecture also includes
a description of the actions of the Frobenius morphism
and of the Hecke operators.

It would take too long to explain what the parameters $\varphi$ are, but
we can say that they give rise to triples $(\gamma_0,\gamma,\delta)$,
where:
\begin{itemize}
\item $\gamma_0$ is a semisimple element of $G(\Q)$, given up to
$G(\oQ)$-conjugacy;
\item $\gamma=(\gamma_\ell)_{\ell\not=p}$ 
is an element of $G(\Af^p)$, given up to $G(\Af^p)$-conjugacy,
and such that $\gamma_\ell$ and $\gamma_0$ are conjugated under
$G(\oQ_\ell)$ for every $\ell\not=p$;
\item $\delta$ is an element of $G(L)$, where $L$ is an unramified
extension of degree $r=[\Fi_q:\Of/\wp]$ of $F_\wp$, such that, if we
denote by $\sigma\in\Gal(L/\Q_p)$ the lift of the (arithmetic) Frobenius, 
then $N(\delta)=\delta\sigma(\delta)\ldots\sigma^{r-1}(\delta)$ is
$G(\oQ_p)$-conjugate to $\gamma_0$.

\end{itemize}
There are some more conditions on the triple $(\gamma_0,\gamma,\delta)$,
see for example Section~2 of \cite{K-SV}. 
One of these conditions is that a certain invariant $\alpha(\gamma_0;\gamma,\delta)$
defined by Kottwitz should vanish, and it implies that there exists an inner
form $I$ of $G_{\gamma_0}$, the centralizer of $\gamma_0$ in $G$,
such that $I_{\Q_v}\simeq G_{\Q_v,\gamma_v}$ for $v\not=p,\infty$, 
$I_{\Q_p}$ is isomorphic to the twisted centralizer of $\delta$ (see page 169 of~\cite{K-SV})
and $I_\R$ is anisotropic modulo the center of $G$.
We then take $I_\varphi=I$,
\[X^p(\varphi)=\{g\in G(\Af^p)/K^p\mid g^{-1}\gamma g\in K^p\}\]
and
\[X_p(\varphi)=\{g\in G(L)/G(\Of_L)\mid g^{-1}\delta\sigma(g)\in
G(\Of_L)\mu_h(\varpi_L)G(\Of_L)\},\]
where we extended $G$ to a reductive group scheme over $\Z_p$,
$\varpi_L$ is a uniformizer of $L$, and
$\mu_h$ is the morphism $h_\C\circ r:\GL_{1,\C}\ra G_\C$, seen as
a conjugacy class of morphisms $\GL_{1,\oQ_p}\ra G_{\oQ_p}$ that is stable
by $\Gal(\oQ_p/L)$; as $G$ is quasi-split over $L$ for $r$ big enough,
we may assume up to taking $r$ big enough that $\mu_h$ is defined over
$L$ (see Lemma~1.1.3 of~\cite{K-TOI}).

If the Shimura datum $(G,h)$ is PEL of type $A$ or $C$, so that
$\Mcal_K$ has a modular description as in Definition~\ref{def_PEL_moduli},
then the triples $(\gamma_0,\gamma,\delta)$ should parametrize
the $\Q$-isogeny classes of triples $(A,\lambda,\iota)$ as in the moduli
problem. This parametrization rests on Honda-Tate theory, which classifies
abelian varieties $A$ over finite fields using their Frobenius, seen as
a central element of $\End(A)\otimes_\Z\Q$. See Part III of
\cite{K-SV} for the case of Siegel modular varieties, and
\cite{K-NP} for the PEL cases of type $A$ and $C$.

The Langlands-Rapoport conjecture is not known for general Shimura varieties,
because we do not even have integral models for general Shimura varieties.
It is also not known for Shimura varieties of abelian type,
except in some particular cases such as the case of quaternionic
Shimura varieties (see the paper \cite{Reim} of Reimann),
but there are some weaker versions that are proved and suffice for the
application to the Kottwitz conjecture.
Here is my understanding of the current situation:
\begin{itemize}
\item If $(G,h)$ is of PEL type and $G_\der$ is simply connected (i.e.
$G$ is of type $A$ or $C$), then Kottwitz formulated a (strictly)
weaker form of the
Langlands-Rapoport conjecture in \cite{K-SV} and proved this reformulation
in \cite{K-NP}. This version is enough to provide a counting-point
formula that can be compared to the Arthur-Selberg trace formula (see below).

\item Dong Uk Lee (see \cite{DUL})
generalized Kottwitz's results for Shimura varieties of Hodge type.
The methods, and even the formulation itself,
do not seem to be generalizable to Shimura varieties
of abelian type.

\item Before that, Kisin has proved
his own weaker version of the Langlands-Rapoport conjecture in \cite{Kisin-LR}
for Shimura varieties of abelian type. This version of the conjecture
is not obviously weaker or stronger than Kottwitz's version. One of its
advantages is that it makes it easier to go from Shimura varieties of
Hodge type to Shimura varieties of abelian type, and one drawback is
that it does not imply a counting-point formula that we can use with the
Arthur-Selberg trace formula.

\item In \cite{KSZ}, Kisin, Shin and Zhu, working with Shimura varieties
of abelian type, formulated and proved a version
of the Langlands-Rapoport conjecture that is weaker than the original
conjecture but stronger than Kottwitz's and Kisin's previous versions, 
and that is enough for
the application to the Kottwitz conjecture. (In fact they did more than that, as
they also ``stabilized'' the expression that they got; see
Remark~\ref{rmk_stab}(1).)

\end{itemize}

The upshot is that, if $g\in G(\Af)$ has a trivial component at
the prime number $p$ under $\wp$ and $f^\infty=\ungras_{KgK}\in\Hf_K$, 
we get a formula for the
trace of $\Frob_\wp^r\cdot f^\infty$ on the $\ell$-adic cohomology of
$M_K(G,h)$ involving terms such as orbital integrals for
$f^\infty$ (i.e. integrals of $f$ over $G(\Af)$-conjugacy orbits
of elements of $G(\Q)$) and twisted orbital integrals of
a function (depending on $r$) at $p$; a priori we only get this for
$r$ big enough, but then the identity can be extended to all
non-negative $r$.
\footnote{This does require some work. We need to see that both sides are finite sums
of terms of the form $a\cdot b^r$ with $a,b\in\C$, which is not obvious for the
twisted orbital integrals.}
This is the kind of input that
we can plug into the geometric side of the Arthur-Selberg trace formula,
see~\ref{sssection_Arthur_TF}. 
The spectral side of the trace formula will then give us an expression
that can be massaged into what we want, that is, in cases when the simplest
form of the Kottwitz conjecture applies, the trace of
$\Frob_\wp^r\cdot f^\infty$ on the virtual representation
\[\sum_{i\geq 0}(-1)^i\sum_{\pi\in\Pi_G}a(\pi_f)\pi_f^K\otimes
(r_\mu\circ\sigma_\pi)\]
of $\Hf_K\times\Gal(\overline{F}/F)$.
So we win.

\begin{remark}\label{rmk_stab}
\begin{itemize}
\item[(1)] The sentence ``This is the kind of input that
we can plug into the geometric side of the Arthur-Selberg trace formula''
in the previous paragraph is sweeping 
a lot of difficulties under the rug. If we are looking at cases 
where $G$ has no endoscopy (such as the simple Shimura varieties of
Kottwitz), then the formula given by the Langlands-Rapoport conjecture
is not too far from the geometric side of the trace formula
(see Section~4 of \cite{K-SSV}). In general, we must first perform a
complicated process known as ``stabilization'', which uses difficult
results such as the fundamental lemma, cf.~\ref{sssection_Langlands_corr}.
See Section~4 of Kottwitz's paper \cite{K-SV} for an explanation of stabilization in the simpler case when $G_\der$ is simply connected, and the
book \cite{KSZ} of Kisin-Shin-Zhu for the case of general Shimura data.

\item[(2)] Note that we never used Matsushima's formula in our outline
of the proof of the Kottwitz conjecture. In fact, though it serves
as a guide, Matsushima's formula is not logically necessary to the proof.

\end{itemize}
\end{remark}

\subsubsection{The Arthur-Selberg trace formula for a compact quotient}
\label{sssection_Arthur_TF}

We consider the situation of~\ref{sssection_dar},
so $G$ is a connected reductive group over $\Q$, $S_G$ is the maximal
$\Q$-split torus in the center of $G$ and $A_G=S_G(\R)^0$.
We write $L^2_G$ for $L^2(G(\Q)\sous G(\Ade)/A_G,\C)$, with the action
of $G(\Ade)$ given by right translations.

We also assume that $G_\der$ is of $\Q$-rank $0$, so that
the quotient $G(\Q)\sous G(\Ade)/A_G$ is compact.
Then the representation $L^2_G$ is semi-simple, and we have
\begin{equation}\label{eq_L2}
L^2_G\simeq\widehat{\bigoplus_{\pi\in\Pi_G}}\pi^{m(\pi)}.
\end{equation}

Fix a Haar measure on $G(\Ade)$.
If $f$ is a smooth function with compact support on $G(\Ade)$, then
we write $R(f)$ for the action of $f$ on $L^2_G$ by right convolution.
Thanks to our hypothesis on $G$, the operator $R(f)$ is of trace class,
and the goal of the Arthur-Selberg trace formula is to give two expressions
of its trace.

\begin{theorem}[See Section~1 of \cite{A-ITF} or \cite{Boo}]
We have
\begin{align*}
\Tr(R(f)) & = \sum_{\pi\in\Pi_G}m(\pi)\Tr(\pi(f)) \\
& = \sum_{\gamma\in G(\Q)/\sim}\vol(G_\gamma(\Q)\sous G_\gamma(\Ade)/A_G)
\int_{G_\gamma(\Ade)\sous G(\Ade)}f(x^{-1}\gamma x)dx,
\end{align*}
where, in the first formula, $\pi(f)$ is the operator
$\int_G f(x)\pi(x)dx$ acting on the space of $\pi$ and,
in the second formula, the sum is over all elements of
$G(\Q)$ modulo conjugation and, if $\gamma\in G(\Q)$, we denote
by $G_\gamma$ the centralizer of $\gamma$ in $G$.

\end{theorem}

Note that we need to choose Haar measures on the groups $G_\gamma(\Ade)$
to make sense of the second formula for $R(f)$
(we also use the counting measure on
$G_\gamma(\Q)$), but that the result does not depend on that choice.

The first formula for $R(f)$ is called the \emph{spectral side}. It follows
from the isomorphism of~\eqref{eq_L2}.

The second formula for $R(f)$ is called the \emph{geometric side}.
We can deduce it by noting that $R(f)$ is an integral operator with
kernel
\[K(x,y)=\sum_{\gamma\in G(\Q)}f(x^{-1}\gamma y),\]
so the trace of $R(f)$ must be equal to $\int_G K(x,x)dx$. See
Section~1 of \cite{A-ITF} for more details.

\subsection{The general Kottwitz conjecture}
\label{ssection_KC}

Let us finally explain briefly how the Kottwitz conjecture generalizes
to arbitrary Shimura varieties. There are two sources of difficulty:
\begin{itemize}
\item The fact that general Shimura varieties are not compact, which 
introduces complications in the geometry: 
we will have to deal with compactifications and choose which cohomology
theory to use.

\item Endoscopy, which will make the Arthur-Selberg trace formula harder
to analyze.

\end{itemize}

\begin{remark}
Although we are concentrating on the case of constant coefficients for
simplicity, we
could also consider cohomology with coefficients in a local system defined
by a representation of the group $G$, as in Section~5 of Pink's 
paper~\cite{Pink-l-adic}. In fact, it is necessary to work in this
more general situation to make some inductive arguments work, but it does
not complicate the calculation of cohomology in an essential way.

\end{remark}

\subsubsection{Compactifications}

There are many ways to compactify locally symmetric spaces, see for example
the book~\cite{BJ} of Borel and Ji for a detailed presentation.
However, in these notes we are interested in compactifications of Shimura
varieties that have the structure of an algebraic variety, and then
there are two choices:
\begin{itemize}
\item[(1)] The Baily-Borel compactification, also called minimal Satake
compactification, Satake-Baily-Borel compactification or just
minimal compactification. Its main advantage is that it is canonical, so
for example Hecke correspondences will extend to the Baily-Borel
compactifications of Shimura varieties, and that its boundary is stratified
by lower-dimensional Shimura varieties.\footnote{Almost: the boundary strata
  are actually quotients of Shimura varieties by finite groups.}
Its main drawback is that it is
very singular in general; for example, the boundary of the Baily-Borel
compactification of the Siegel modular variety $\Mcal_{d,n}$ has codimension
$d$.

\item[(2)] The toroidal compactifications, which map to the Baily-Borel
compactification. As the plural indicates, they
are not canonical and depend on extra data, so we cannot extend the action
of the Hecke algebra to the cohomology of any fixed toroidal compactification.
On the other hand, we can always find toroidal compactifications which are
smooth projective and whose boundary is a divisor with normal crossings, and
the boundary of a toroidal compactification is stratified by so-called
\emph{mixed Shimura varieties}, which are a generalization of Shimura
varieties to certain non-reductive groups.
Even if we were only interested in the Baily-Borel compactification, the
toroidal compactifications would be a very useful tool to study it; they
are also interesting in their own right.

\end{itemize}

Sections~I.4, I.5 and~III.4 of the book~\cite{BJ} contain a
review of Satake and Baily-Borel compactifications of locally
symmetric spaces as topological spaces and complex analytic spaces.
The fact that the Baily-Borel compactification is a projective algebraic
variety over $\C$ was proved by Satake for Siegel modular varieties, and
by Baily and Borel in general, see the paper~\cite{BB}. As for toroidal
compactifications of complex Shimura varieties, they are studied in the
book~\cite{AMRT} of Ash, Mumford, Rapoport and Tai. In his
thesis~\cite{Pink-these}, Pink constructed canonical models of the
toroidal compactifications over the reflex field, and deduced that
the Baily-Borel compactification also descends to a projective algebraic
variety over that field. As in the case of Shimura varieties, the story of
integral models is not yet complete,\footnote{And as in the case of Shimura
  varieties, we only discuss places of good reduction here.}
and as for canonical models, the crucial part is to construct models
of the toroidal compactifications (we then get a model of the Baily-Borel
compactification as a byproduct).
The case of modular curves is treated in the paper~\cite{DeRa} of
Deligne-Rapoport, and that of Siegel modular varieties in the book~\cite{CF}
of Chai and Faltings. The next case to be handled systematically was that
of PEL type Shimura varieties, in Lan's thesis~\cite{Lan-PEL}. Then
Madapusi Pera treated the case of Hodge type Shimura varieties in his
paper~\cite{Mada-Hodge}. The case of abelian type Shimura varieties seems
to be within reach, but as far as I know it has not been written yet.

\vspace{.3cm}

\subsubsection{Intersection cohomology and Zucker's conjecture}

It is perfectly possible to try calculating compact support cohomology
of noncompact Shimura varieties with the methods oulined in~\ref{proof_KC},
this is for example what Laumon does in \cite{Lau1}, \cite{Lau2}.
The issue is that the Arthur-Selberg trave formula will be harder to use,
because there is no simple spectral description of the cohomology such
as the one given by Matsushima's formula in the compact case.
Indeed, the generalization of Matsushima's formula, which is a theorem
of Borel and Casselman (see~\cite{BC}) gives a spectral description of the
$L^2$ cohomology of the Shimura variety, so it would be simpler to study that.
The next problem is that $L^2$ cohomology has a very non-algebraic
definition (it is the cohomology of the complex of $L^2$ differential forms,
for a certain metric on the locally symmetric space); it has an action of
the Hecke algebra, but it is not clear
at all a priori how to make the absolute Galois group of the reflex field
act on it.

Fortunately, this problem is solved by Zucker's conjecture, which says that
there should be a Hecke-equivariant isomorphism between $L^2$ cohomology
of a Shimura variety and intersection cohomology of its Baily-Borel
compactification. Intersection cohomology is a cohomology theory adapted
to singular spaces.
It was first introduced by Goresky and MacPherson
(see~\cite{GM1}) as a modification of Betti
cohomology for stratified spaces, and then further studied by them in~\cite{GM2}
using a sheaf-theoretic interpretation due to Deligne. Then, in the
book~\cite{BBD}, Beilinson, Bernstein, Deligne and Gabber put
intersection cohomology into the wider context of perverse sheaves:
if $X$ is a possibly singular space, its intersection cohomology is the
cohomology of a particular complex of sheaves on it, called the
\emph{intersection complex}, which is a simple object in the abelian
category of perverse sheaves (for the middle perversity). 
We have been deliberately vague about the
precise nature of $X$, as one of the advantages of the sheaf-theoretic
point of view is that it can be adapted to many different situations;
for example, $X$ could be an algebraic variety over a field and the
intersection complex could be an \'etale complex with $\Q_\ell$-coefficients; 
or $X$ could be an algebraic variety over $\C$ and the intersection complex
could be a complex of $\Q$-vector spaces for the analytic topology, or
a mixed Hodge module. 
\footnote{We are avoiding the cases of non-algebraic $X$ or integral
coefficients, as there is no unique ``middle perversity'' there. But the
theory still makes sense, it is just slightly more complicated.}
In any case, 
the intersection cohomology of $X$ satisfies some of the properties that
we would expect from the cohomology of a nonsingular variety:
\begin{itemize}
\item Poincar\'e duality if $X$ is proper;
\item the hard Lefschetz theorem if $X$ is projective;
\item purity if $X$ is proper (in the appropriate sense, e.g.
in the sense of Deligne (see~\cite{De-WeilI}) if $X$ is a variety over
a finite field or as a Hodge structure if $X$ is a variety over $\C$ and
we are using mixed Hodge modules).

\end{itemize}
It also satisfies comparison theorems between the different cohomology theories.
(See for example Chapter~6 of~\cite{BBD} for the comparison between the
\'etale and Betti theories.)

\vspace{.3cm}

\subsubsection{The Borel-Casselman theorem}
\label{sssection_BC}

We come back to Shimura varieties. Let $(G,h)$ be a Shimura datum,
$K\subset G(\Af)$ be a small enough level, $X$ be the Shimura variety
$M_K(G,h)$, $X^*$
be the Baily-Borel compactification of $X$ and $IH^*(X^*)$ be its
intersection cohomology. We have seen that Hecke
correspondences
extend to Baily-Borel compactifications. Although intersection cohomology
is not as functorial as usual cohomology, we can still make
the Hecke algebra at level $K$ act on $IH^*(X^*)$ (because
Hecke operators are finite morphisms between Baily-Borel compactifications
and become finite \'etale when restricted to Shimura varieties). Zucker's
conjecture (see ~\cite[(6.20)]{ZuC}) says that, if we use the Betti version of
$IH^*(X^*)$ with complex coefficients, then there is a
Hecke-equivariant isomorphism between $IH^*(X^*)$ and $L^2$
cohomology of $X$. Also, Zucker's conjecture is actually a theorem, as it has
been proved (at least) four times: by Looijenga (\cite{Loo}),
Looijenga-Rapoport (\cite{LR}), Saper-Stern (\cite{SS}) and
Saper (\cite{Saper}).
If we combine it with the generalization of Matsushima's formula
due to Borel and Casselman (\cite{BC}), we finally get the following
theorem (note the similarity with Theorem~\ref{thm_Matsushima}):

\begin{theorem}\label{thm_Matsushima_general}
Let $(G,h)$ be a Shimura datum. Then we
have a $G(\Af)$-equivariant isomorphism of graded $\C$-vector spaces
\[\varinjlim_K IH^*(M_K(G,h)^*(\C))\simeq\bigoplus_{\pi\in\Pi_G}
\pi_f\otimes\H^*(\ggoth,K_\infty;\pi_\infty)^{m(\pi)},\]
where $G(\Af)$ acts on the factors $\pi_f$ on the right hand side,
$M_K(G,h)^*$ is the Baily-Borel compactification of $M_K(G,h)$ and we
are using the Betti version of intersection cohomology (with constant
coefficients for simplicity).

\end{theorem}

\subsubsection{The Kottwitz conjecture}

Let $F$ be the reflex field of $(G,h)$. If we use the \'etale $\ell$-adic
version of intersection cohomology, then we get cohomology groups
$\varinjlim_K IH^*(M_K(G,h)^*_{\overline{F}},\oQ_\ell)$ equipped with
commuting actions of $G(\Af)$ and $\Gal(\overline{F}/F)$, which become
isomorphic as $G(\Af)$-modules to the cohomology groups of
Theorem~\ref{thm_Matsushima_general} if we choose any isomorphism
$\oQ_\ell\simeq\C$. We introduce the following virtual representation
of $G(\Af)\times\Gal(\oQ_\ell/\Q_\ell)$:
\[IH=\sum_{i\geq 0}(-1)^i[\varinjlim_K IH^i(M_K(G,h)^*_{\overline{F}},\oQ_\ell)],\]
where the brackets denote the class of a representation in the group
of virtual representations.
The sum is actually finite, because the intersection cohomology group
vanishes for $i\geq 2\dim(M^K(G,h))$, and this dimension is independent of
the level $K$. We have a decomposition
\[IH=\sum_{\pi_f}[\pi_f]\otimes\sigma(\pi_f),\]
where $\pi_f$ goes through the set of equivalence classes of
irreducible representations of $G(\Af)$ and each $\sigma(\pi_f)$ is a
virtual representation of $\Gal(\overline{F}/F)$.
By Theorem~\ref{thm_Matsushima_general}, we know that the only $\pi_f$ that
can occur are the finite parts of discrete automorphic representations of
$G$, and we have a formula for the dimension of $\sigma(\pi_f)$.
The Kottwitz conjecture, stated on page 201 of~\cite{K-SV},
gives a precise description of $\sigma(\pi_f)$.

We will not state a precise form of the conjecture, as it would
require introducing too much notation, but we will try to give its flavor.
The most naive guess is that we should have something like
\[\sigma(\pi_f)=\pm\sum_{\pi'\in\Pi_G,\ 
\pi'_f\simeq\pi_f}m(\pi')[(r_\mu\circ\sigma_{\pi'})],\]
where $r_\mu$ and $\sigma_\pi$ are as in Conjecture~\ref{conj_Kottwitz},
but this is not true because of endoscopy. What we actually expect is as
follows. As mentioned in Remark~\ref{rmk_Arthur_param}, discrete
automorphic representations of $G$ should be classified by \emph{Arthur
parameters} $\psi:\Lf_\Q\times\SL_2(\C)\ra{}^LG$. The parameters corresponding
to the representations appearing in the cohomology of Shimura varieties
are expected to factor through the motivic Galois group of $\Q$ and
to be defined over a finite extension of $\Q$ in $\C$, 
and there is an explicit
recipe giving $\sigma_\pi$ from the parameter corresponding to
$\pi$; we denote by $V_\psi$ the space of the representation
$r_\mu$ with the action of $\Gal(\overline{F}/F)$ coming from $\sigma_\pi$.
\footnote{We are cheating here: if we want an action of $\Gal(\overline{F}/F)$
and not of the motivic Galois group of $F$, we have to take as coefficients
an $\ell$-adic completion of the field of definition of $\psi$.}
Let $S_\psi$ be the centralizer of $\psi$ in $\widehat{G}$; this group, or
rather its quotient $S_\psi/S_\psi^0 Z(\widehat{G})^{\Gal(\oQ/\Q)}$, should
control the structure of the Arthur packet of representations associated
to $\psi$ and the multiplicities of these representations. We can decompose
$V_\psi$ as a direct sum
\[V_\psi=\bigoplus_{\nu}V_{\psi,\nu},\]
where $\nu$ runs over characters of $S_\psi$ and $V_{\psi,\nu}$ is the part
where $S_\psi$ acts by $\nu^{-1}$. Then the Kottwitz conjecture says that
\[\sigma(\pi_f)=\sum_\psi\sum_\nu n(\pi_f,\nu)[V_{\psi,\nu}],\]
where $\psi$ runs over Arthur parameters whose packets contains a
representation with finite part $\pi_f$, and $n(\pi_f,\nu)$ is the product of
a sign and of a multiplicity (which is the multiplicity of $\nu$ in
some finite-dimensional representation of $S_\psi$ depending on $\pi_f$).

\vspace{.3cm}

\subsubsection{Proving the Kottwitz conjecture}

The general strategy of~\ref{proof_KC} still applies with some adaptations:
\begin{itemize}
\item[(1)] The specialization theorem still holds for intersection cohomology,
as long as we have toroidal compactifications that are proper and smooth
over the base ring (which is the case for Hodge type Shimura varieties
at places of good reduction). This follows from the specialization theorem
for the cohomology of the toroidal compactifications and from the
decomposition theorem (see Corollaires~5.4.2 and~5.4.6 of~\cite{BBD}),
which says that intersection cohomology of the Baily-Borel compactification
is a direct factor of this cohomology.

\item[(2)] Deligne's conjecture makes sense (and is known) for general
complexes and cohomological correspondences, not just for constant
coefficients. It reduces the calculation of the trace of a Hecke
operator twisted by a power of the Frobenius on intersection cohomology
to counting fixed points (in the Baily-Borel compactification) and calculating
the stalks of the intersection complex at these fixed points.

\item[(3)] The Langlands-Rapoport conjecture makes no distinction between
compact and noncompact Shimura varieties. It also theoretically
gives a description
of points of the Baily-Borel compactification over a finite field, as
this compactification is stratified by Shimura varieties (or explicit
finite quotients thereof).

\item[(3')] The intersection complex can be quite difficult to calculate,
but in the case of Shimura varieties there is an explicit description
given by the theory of weighted cohomology of Goresky-Harder-MacPherson
(see~\cite{GHM}); in particular, it is a complex of local systems on each
stratum of the Baily-Borel compactification. Once the Langlands-Rapoport
conjecture (or an appropriate weakening) is known, this reduces the
calculation of the cohomology to a complicated combinatorial problem.

\item[(4)] As for the Arthur-Selberg trace formula, we now need to use
the stable trace formula for a noncompact quotient, fortunately in a
somewhat simple form (see Section~29 of Arthur's overview~\cite{A-ITF}
for an introduction to the stable trace formula, and Arthur's
paper~\cite{A-L2} for the simplifications that occur when calculating the
cohomology of Shimura varieties).
The stabilization process that we must apply to the point-counting
formula coming from (3) and (3') is also more complicated, because we need
to stabilize the terms coming from the boundary. 
\footnote{Some of these complications already arise in the compact case
with endoscopy. Indeed, even when the group $G$ has no nontrivial
parabolic subgroups, its endoscopic groups almost certainly do because
they are quasi-split, and
this produces ``boundary terms'' in their stable trace formulas 
on the right-hand side of identity~\eqref{eq_KC}. These terms must
cancel each other out, as there are no corresponding terms on the
left-hand side of~\eqref{eq_KC}.}
For an example of this process, see Yihang Zhu's thesis~\cite{Zhu-phd}.

In the end, we expect to get the following formula for $h\in\Hf_K$ a
Hecke operator trivial at $p$, $\Frob_\wp$ the geometric Frobebius
at a place $\wp$ of $F$ above $p$ and $j\in\Nat$:
\begin{equation}\label{eq_KC}
\Tr(h\times\Frob_\wp,IH^*(M_K(G,h)_{\overline{\Fi_p}},\oQ_\ell))=
\sum_H \iota(G,H)ST_H((h^p)^H\varphi_{p,j}\varphi_\infty),
\end{equation}
where $H$ runs over elliptic endoscopic groups of $G$ (to be more precise,
it should run over elliptic endoscopic triples), $\iota(G,H)$ is a positive
rational number, $ST_H$ is the stable trace formula for $H$,
$(h^p)^H$ is the so-called \emph{endoscopic transfer} to $H(\Af^p)$ of
the prime-to-$p$ part of $h$, $\varphi_{p,j}$ is a function in the spherical
Hecke algebra of $G(\Q_p)$ that only depends on $j$ (not on $h$), and
$\varphi_\infty$ is a function on $G(\R)$ that depends only on the coefficients
of the cohomology (here we are taking constant coefficients, but we
could take coefficients in any algebraic representation of $G$).

\item[(5)] The last step
is to transform the right-hand side of the previous
formula into an expression of the form
$\sum_{\pi_f}\Tr(h,\pi_f)\Tr(\Frob_\wp^j,\sigma(\pi_f))$.
Kottwitz explains a way to do this in Section~10 of~\cite{K-SV}, but his
method requires Arthur's conjectural classification of discrete automorphic
representations, which is only available for classical and unitary groups
(see~\cite{Arthur-fonc}, \cite{Mok} and~\cite{KMSW}) so far, while
the groups of Shimura varieties of PEL type (i.e. the best understood Shimura
varieties) tend to be \emph{general} symplectic, orthogonal
or unitary groups (i.e. $\GSp$, $\GO$ and $\GU$); 
there are some results for general symplectic groups,
see Bin Xu's thesis~\cite{Xu1} and his paper~\cite{Xu2}, though as far as
I understand they are not sufficient to completely prove the Kottwitz
conjecture.
It also requires the agreement of the various constructions of the
classification of representations of real groups, that was proved for
classical and unitary groups in papers of Arancibia-Moeglin-Renard
(\cite{ArMo}), Arancibia (\cite{Ar}),
Arancibia-Mezo (\cite{ArMe}) and Adams-Arancibia-Mezo
(\cite{AArMe}).

\end{itemize}

\bibliography{Shimura_varieties}

\end{document}